\documentclass{article}
\usepackage{amssymb}
\usepackage{amsmath}
\usepackage{amsthm}
\usepackage{amsfonts}
\usepackage{eucal}
\usepackage{xspace}
\usepackage{color}
\usepackage{a4wide}
\usepackage{enumerate}
\usepackage[titletoc,toc,title]{appendix}
\usepackage[colorlinks=true,citecolor=red,linkcolor=blue,urlcolor=RubineRed,pdfpagetransition=Blinds,pdftoolbar=false,pdfmenubar=false]{hyperref}
\usepackage{color}
\usepackage{xcolor}
\definecolor{vi}{RGB}{200,0,200} 
\providecommand{\U}[1]{\protect\rule{.1in}{.1in}}
\providecommand{\U}[1]{\protect\rule{.1in}{.1in}}
\providecommand{\U}[1]{\protect\rule{.1in}{.1in}}
\everymath{\displaystyle}
\newtheorem{theorem}{Theorem}[section]
\newtheorem{lemma}{Lemma}[section]
\newtheorem{proposition}{Proposition}[section]

\theoremstyle{definition}
\newtheorem{remark}{Remark}[section]
\newtheorem{definition}{Definition}[section]
\newcommand{\field}[1]{\ensuremath{\mathbb{#1}}}
\newcommand{\C}{\field{C}\xspace}

\newcommand{\R}{\field{R}\xspace}
\newcommand{\N}{\field{N}\xspace}
\newcommand{\intdoble}{\displaystyle \int \!\!\!\! \int}

\begin{document}


	\title{Boundary controllability of a one-dimensional phase-field system with one control force}
		
\author{Manuel \textsc{Gonz\'alez-Burgos}\thanks{Dpto.~Ecuaciones Diferenciales y An\'alisis Num\'erico and Instituto de Matem\'aticas de la Universidad de Sevilla (IMUS), Facultad de Matem\'aticas, Universidad de Sevilla, C/ Tarfia S/N, 41012 Sevilla, Spain. Supported by grant MTM2016-76990-P, Ministry of Economy and Competitiveness
(Spain). E-mail: \texttt{manoloburgos@us.es}}, Gilcenio R.~\textsc{Sousa-Neto}\thanks{Centro Acad\^emico do Agreste, NICIT, Universidade Federal de Pernambuco, 55002-970, Caruaru, PE, Brazil (\texttt{gilceniorodrigues@gmail.com}). Partially supported by CAPES (Brazil) and MathAmSud COSIP.}}
	\date{}
\maketitle


	\begin{abstract}
In this paper, we present some controllability results for linear and nonlinear phase-field systems of Caginalp type considered in a bounded interval of $\R$ when the scalar control force acts on the temperature equation of the system by means of the Dirichlet condition on one of the endpoints of the interval. In order to prove the linear result we use the moment method providing an estimate of the cost of fast controls. Using this estimate and following the methodology developed in~\cite{tucsnack}, we prove a local exact boundary controllability result to constant trajectories of the nonlinear phase-field system. To the authors' knowledge, this is the first nonlinear boundary controllability result in the framework of non-scalar parabolic systems, framework in which some ``hyperbolic'' behaviors could arise. 
%
	\end{abstract}




\section{Introduction}\label{s1}
\setcounter{equation}{0}

This work deals with the boundary controllability properties of a phase-field system of Caginalp type (see~\cite{cag}) which is a model describing the transition between the solid and liquid phases in solidification/melting processes of a material occupying an interval:
	\begin{equation}
	\left\{
	\begin{array}{ll}
\displaystyle \tilde{\theta}_t - \xi\tilde{\theta}_{xx} + \dfrac{1}{2}\rho\xi\tilde{\phi}_{xx} + \dfrac{\rho}{\tau}\tilde{\theta} = f_1(\tilde{\phi})	& \mbox{in } Q_T :=(0,\pi) \times (0,T), 
	\\
	\noalign{\smallskip}
\displaystyle \tilde{\phi}_t - \xi\tilde{\phi}_{xx} - \dfrac{2}{\tau}\tilde{\theta} = f_2(\tilde{\phi}) & \mbox{in }  Q_T, 
	\\
	\noalign{\smallskip}
\displaystyle \tilde{\theta}(0,\cdot) = v,\ \tilde{\phi}(0,\cdot) = c,\ \tilde{\theta}(\pi,\cdot)=0 , \ \tilde{\phi}(\pi,\cdot) = c & \mbox{on }  (0,T),	 
	\\
	\noalign{\smallskip}
\displaystyle \tilde{\theta}(\cdot,0) = \tilde{\theta}_0, \  \tilde{\phi}(\cdot,0)=\tilde{\phi}_0 & \mbox{in }  (0,\pi).		\end{array}
	\right. \label{PFSy}
	\end{equation}
Here, $T>0$ is some final time, $\tilde{\theta} = \tilde{\theta}(x,t)$ denotes the temperature of the material, $\tilde{\phi} = \tilde{\phi} (x,t)$ is the phase-field function used to identify the solidification level of the material, $c\in \{ - 1,0,1\}$ and the functions $f_1$ and $f_2$ are the nonlinear terms which come from the derivative of the classical regular double-well potential $W$ and are defined by
	\[ 
\displaystyle f_1(\tilde{\phi}) = -\frac{\rho}{4\tau}\left( \tilde{\phi}-\tilde{\phi}^3 \right) \quad \mbox{and} \quad f_2(\tilde{\phi}) = \frac{1}{2\tau} \left( \tilde{\phi}-\tilde{\phi}^3 \right). 
	\]
Besides, $\rho > 0$ is the latent heat, $\tau > 0$ represents the relaxation time and $\xi > 0$ is the thermal diffusivity. Finally, $v \in L^2(0,T)$ is the control force, which is exerted at point $x=0$ by means of the boundary Dirichlet condition, and the initial data $\tilde \theta_0$, $\tilde \phi_0$ are given functions.


The phase function $\tilde{\phi}$ describes the phase transition of the material (solid or liquid) in such a way that $\tilde \phi = 1$ means that the material is in solid state and $\tilde \phi = - 1$ in liquid state. Observe that the temperature $\tilde \theta$ of the material could be zero and this could occur with the material in solid or liquid phase. On the other hand, the phase-field variable $\tilde{\phi}$ does not have a direct physical meaning. This is the reason way we control the temperature $\tilde \theta$ which, in fact, is the unique variable with physical meaning. 

The objective of this paper is to prove a null controllability result at time $T$ for the temperature variable $\tilde \theta$ of system~\eqref{PFSy}. If we consider the transition region associated to the temperature, i.e., the set
	$$
\Gamma (t) := \left\{ x \in (0, \pi) : \tilde \theta (x, t) = 0 \right\},
	$$
then, the problem under consideration consists of proving that there exists a control $v$ such that the transition region associated to the temperature $\tilde \theta $ satisfies $\Gamma (T) = (0, \pi)$. It is interesting to underline that in this case the material could be in solid phase ($\tilde \phi (\cdot, T)=1$), liquid phase ($\tilde \phi (\cdot, T)=-1$) or in an intermediate phase (mushy) which corresponds to $\tilde \phi (\cdot, T)=0$. In this work we are interested in showing the null controllability result at time $T$ for the temperature $\tilde \theta$ but keeping the material in solid state, $c=1$, or liquid state, $c=-1$, at time $T$, that is to say, proving that there exists a control $v\in L^2 (0,T)$ such that system~\eqref{PFSy} has a solution $\tilde y = (\tilde \theta, \tilde \phi)$ (in an appropriate space) such that
	\begin{equation}\label{exact}
\tilde \theta (\cdot, T) = 0 \quad \hbox{and} \quad \tilde \phi (\cdot, T) = c \quad \hbox{in } (0, \pi).
	\end{equation}
We give a complementary analysis and results in the Appendix~\ref{APPX2_}, where we deal with the case where $c=0$.

\smallskip

As said before, the objective of this work is to study the controllability properties of system~\eqref{PFSy}. Let us observe that we are exerting only one control force on the system (a boundary control) but we want to control the corresponding state $\tilde y = (\tilde \theta, \tilde \phi) $ which has two components. In fact, the second equation in~\eqref{PFSy} is indirectly controlled by means of the term $-2 \tilde \theta/ \tau$. On the other hand, \eqref{PFSy} is a nonlinear system with nonlinearities with a super-linear behavior at infinity. Therefore, we can expect a local controllability result at time $T$ for this system, that is to say, an exact controllability result to the trajectory $(0,c)$ when the initial datum $(\tilde \theta_0, \tilde \phi_0)$ is sufficiently close to $(0,c)$ in an appropriate norm (see for instance~\cite{FC-Z,D-FC-GB-Z} for similar results in the scalar parabolic framework).

%
%
%

System~\eqref{PFSy} is a particular class of more general $n \times n$ nonlinear parabolic control systems of the form:
	\begin{equation}\label{genralsystems}
	\left\{ 
	\begin{array}{ll}
y_{t}-D\Delta y+A y= F(y) + Bv1_{\omega } & \hbox{in }Q_{T}:=\Omega \times (0,T), \\ 
	\noalign{\smallskip} 
y=Cu1_{\Gamma _{0}},\quad & \hbox{on } \Sigma_{T}:=\partial \Omega \times (0,T), \\ 
	\noalign{\smallskip} 
y(\cdot ,0)=y_{0} & \hbox{in }\Omega ,
	\end{array}
	\right.
	\end{equation}
where $\omega $ and $\Gamma _{0}$ are, respectively, open subsets of the smooth bounded domain $\Omega \subset \R^{N}$ and of its boundary $\partial \Omega $, $D \in \mathcal{L} ( \R^{n} ) $, with $n\geq 1$, is a positive definite matrix, $B ,C \in \mathcal{L}(\R \xspace^m,\;\R^n)$, with $m\leq n$, and $A=\left( a_{ij}\right) _{1\leq i,j\leq n}\in \mathcal{L} ( \R^{n} )$ are given matrices. In~\eqref{genralsystems}, $F \in C^0 (\R^n; \R^n)$ is a nonlinear given function. Unlike the scalar case, even in the linear case $F \equiv 0$, new difficulties arise in the study of the controllability properties of~\eqref{genralsystems}. When $m<n$, the issue for this system is to control the whole components of the system with a control function acting, locally in space or on a part of the boundary, only on some of them. We refer to~\cite{ACML4} for a review of results for the controllability problem of system~\eqref{genralsystems}.

The controllability properties of system~\eqref{PFSy} has been analyzed before in the $N$-dimensional case ($N \ge 1$) when a distributed control supported in an open subset of the domain is exerted on the system. The first local controllability results for a nonlinear phase-field system controlled by one distributed control force are proved in~\cite{AK-B-D-K} under certain restrictions on the dimension $N$. In~\cite{GB-PG}, the authors introduce a new approach to deal with the distributed null controllability of some linear coupled parabolic systems that makes possible to generalize the results in~\cite{AK-B-D-K} to more general phase-field systems such as~\eqref{PFSy}. Finally, in~\cite{FC-GB-dT} the authors study the controllability of some (linear and semilinear) non-diagonalizable parabolic systems of PDEs and provide some Kalman rank conditions which characterize the controllability properties in the linear case. In these previous works the null controllability result for the linear and nonlinear problem uses in a fundamental way global Carleman inequalities for scalar parabolic problems. To our knowledge, this is the first time that the boundary controllability properties of a nonlinear phase-field system are analyzed.

It is important to underline that in this work we are considering a boundary controllability problem for a non-scalar parabolic system. As said before, in the study of these boundary controllability problems new phenomena and technical difficulties arise. Let us briefly describe them.  In the linear case ($F \equiv 0$), it is well-known (see~\cite{MGBBoundC,AK-B-GB-dT-2011,MGBminT}) that the distributed ($C \equiv 0$) and boundary ($B \equiv 0$) controllability properties of system~\eqref{genralsystems} are different and not equivalent. In fact, the boundary controllability of system~\eqref{genralsystems} could present ``hyperbolic" behaviors such as the non-equivalence between the approximate and null controllability or the existence of a minimal time of controllability, i.e., the existence of $T_0 \in [0, \infty]$ such that the system is null controllable at time $T$ if $ T > T_0$ and it is not if $T < T_0$ (see~\cite{MGBminT},~\cite{FAMLnewminim} and the references therein for more details). On the other hand, global Carleman inequalities seem not to be too useful when we deal with boundary controllability properties of non-scalar parabolic systems (see~\cite{AK-B-GB-dT-2011}) and this creates a new technical difficulty: we want to obtain a nonlinear boundary controllability result without having global Carleman estimates for the corresponding adjoint systems to linearized versions of system~\eqref{PFSy}. 

%
%

As noted above, the decision of exerting the control force on the temperature variable $\tilde{\theta}$ was taken because it is the only variable in~\eqref{PFSy} with physical meaning. Indeed, the values of $\tilde{\phi}$ determine the material phase and, consequently, imposing a boundary control for $\tilde{\phi}$ would mean that specific phases for the material on the boundary are maintained throughout the solidification process (which is not an usual situation in practice). On the other hand, exerting the control on the boundary in the temperature variable can be seen as having a regulable external source which heats/cools down the material at point $x=0$. From the physical point of view, this boundary control is more interesting than a distributed control supported on an open subset of the domain (internal source).

The main objective of this work is to provide an exact controllability result of system~\eqref{PFSy} to the constant trajectory $(0, c)$ with $c=\pm 1$ (the case $c=0$ follows the same structure of the case $c=\pm 1$, so it will be dealt with in  Appendix~\ref{APPX2_}). Observe that the nonlinearities $f_1$ and $f_2$ in~\eqref{PFSy} can be written as  
	\[
	\left\{
	\begin{array}{l}
\displaystyle f_1(\phi) = -\frac{\rho}{4\tau}(\phi - \phi^3) = \frac{\rho}{2\tau}(\phi - c) \pm \frac{3\rho}{4\tau}(\phi - c)^2 + \frac{\rho}{4\tau}(\phi - c)^3, \\
	\noalign{\smallskip}
\displaystyle f_2(\phi) = \frac{1}{2 \tau}(\phi - \phi^3) = -\frac{1}{\tau}(\phi - c) \mp \frac{3}{2\tau}(\phi - c)^2 - \frac{1}{2\tau}(\phi - c)^3.
	\end{array}
	\right.
	\]
and therefore, performing the change of variable $(\theta, \phi) = (\tilde \theta, \tilde{\phi} - c)$,  system~\eqref{PFSy} becomes
	\begin{equation}\label{bisPFSy}
	\left\{
	\begin{array}{ll}
\theta_t - \xi\theta_{xx} + \dfrac{1}{2}\rho\xi\phi_{xx} - \dfrac{\rho}{2\tau}\phi + \dfrac{\rho}{\tau}\theta  = g_1(\phi) & \mbox{in } Q_T,  \\
	\noalign{\smallskip}
\phi_t - \xi\phi_{xx} + \dfrac{1}{\tau}\phi - \dfrac{2}{\tau}\theta  = g_2(\phi)	 & \mbox{in }  Q_T, \\
	\noalign{\smallskip}
\theta(0,\cdot) = v,\quad  \phi(0,\cdot) = \theta(\pi,\cdot)=\phi(\pi,\cdot) = 0 & \mbox{on } (0,T), \\
	\noalign{\smallskip}
\theta(\cdot,0) = \theta_0, \quad \phi(\cdot,0) =\phi_0 & \mbox{in } (0,\pi),
	\end{array}
	\right.
	\end{equation}
where $(\theta_0, \phi_0) = (\tilde{\theta}_0, \tilde{\phi}_0 - c)$ and the functions $g_1$ and $g_2$ are given by
	\begin{equation}\label{g1g2}
g_1(\phi) = \pm \frac{3\rho}{4\tau}\phi^2 + \frac{\rho}{4\tau}\phi^3  \quad \mbox{and} \quad g_2(\phi) = \mp \frac{3}{2\tau}\phi^2 - \frac{1}{2\tau}\phi^3.
	\end{equation}

With the previous change of variables in mind, the exact controllability to the trajectory $(0,c)$ of system~\eqref{PFSy} at time $T>0$ is equivalent to the null controllability at the same time $T$ of system~\eqref{bisPFSy}. In order to prove the null controllability at time $T >0$ of system~\eqref{bisPFSy}, we will rewrite the controllability problem as a fixed-point problem for a convenient operator in appropriate spaces. To perform this fixed-point strategy, we will first study the controllability properties of the following autonomous linear system:
	\begin{equation}\label{PFSylinear}
	\left\{
	\begin{array}{ll}
\theta_t - \xi\theta_{xx} + \dfrac{1}{2}\rho\xi\phi_{xx} - \dfrac{\rho}{2\tau}\phi + \dfrac{\rho}{\tau}\theta  =0 & \mbox{in } Q_T,  \\
	\noalign{\smallskip}
\phi_t - \xi\phi_{xx} + \dfrac{1}{\tau}\phi - \dfrac{2}{\tau}\theta  = 0	 & \mbox{in }  Q_T, \\
	\noalign{\smallskip}
\theta(0,\cdot) = v,\quad  \phi(0,\cdot) = \theta(\pi,\cdot)=\phi(\pi,\cdot) = 0 & \mbox{on } (0,T), \\
	\noalign{\smallskip}
\theta(\cdot,0) = \theta_0, \quad \phi(\cdot,0) =\phi_0 & \mbox{in } (0,\pi),
	\end{array}
	\right.
	\end{equation}
which is a linearization of system~\eqref{bisPFSy} around the equilibrium $(0,0)$. System~\eqref{PFSylinear} can also be written in the vectorial form
	\begin{equation}\label{vectPFSy}
	\left\{
	\begin{array}{ll}
y_t - D y_{xx} + A y = f & \mbox{in } Q_T,\\
	\noalign{\smallskip}
y(0,\cdot) = Bv, \quad y(\pi,\cdot)=0 & \mbox{on } (0,T), \\
	\noalign{\smallskip}
y (\cdot,0) = y_0, & \mbox{in }  (0,\pi), \\
	\end{array}
	\right. 
	\end{equation}
with $y_0 = ( \theta_0, \phi_0)$, $f =( 0, 0)$ and
	\begin{equation}\label{ADB}
D = \left(
	\begin{array}{cc}
\xi & -\dfrac{1}{2} \rho \xi \\
	\noalign{\smallskip}
0 &  \xi
	\end{array}
	\right),
\quad
A =
	\left(
	\begin{array}{cc}
\dfrac{\rho}{\tau} & -\dfrac{\rho}{2\tau} \\
	\noalign{\smallskip}
-\dfrac{2}{\tau} & \dfrac{1}{\tau}
	\end{array}
	\right),
\quad
B =
	\left(
	\begin{array}{cc}
1 \\ 0
	\end{array}
	\right).
	\end{equation}

We will see that, for every $v \in L^2 (0,T) $, $f \in L^2 (Q_T ;\R^{2})$ and $y_{0}\in H^{-1}(0,\pi ;\R^{2})$, system~\eqref{vectPFSy} possesses a unique solution defined by transposition which satisfies 
	\begin{equation*}
y \in L^{2} (Q_T ;\R^{2}  )\cap C^{0}\left([0,T]; H^{-1}  (0,\pi ;\R^{2} ) \right) ,
	\end{equation*}
and depends continuously on the data $v$, $f$ and $y_{0}$. Observe that the previous regularity permits to pose the boundary controllability of system~\eqref{PFSylinear} in the space $H^{-1}  (0,\pi ;\R^{2} )$.

Let us present our first main result: the boundary approximate controllability at time $T>0$ of system~\eqref{PFSylinear}. One has:

\begin{theorem} \label{CAproximada_}
Let us consider $\xi$, $\rho$ and $\tau$ three positive real numbers and let us fix $T>0$. Then, system~\eqref{PFSylinear} is approximately controllable in $H^{-1} (0, \pi; \R^2)$ at time $T $ if and only if one has 
	\begin{equation}\label{H2}
\xi^2\tau^2(\ell^2 - k^2)^2 - 2\xi\rho\tau(\ell^2 + k^2)- 2 \rho-1 \neq 0, \quad \forall k,\ell\geq1, \quad \ell >k.
	\end{equation}
\end{theorem}
%


\begin{remark}
Condition~\eqref{H2} characterizes the approximate controllability property of system~\eqref{PFSylinear}. Thus,~\eqref{H2} is a necessary condition for the null controllability of this system at time  $T> 0$. Observe that this condition is independent of the final time $T$. We will also see that condition~\eqref{H2} is equivalent to the following property (see Proposition~\ref{propert}): \textit{``The eigenvalues of the vectorial operators
	\begin{equation}\label{LL*}
L = -D \partial_{xx} + A \quad \mbox{and} \quad L^* = -D^* \partial_{xx} + A^*,
	\end{equation}
with domains $D(L)=D(L^*)=H^2(0,\pi;\mathbb{R}^2) \cap H_0^1(0,\pi;\mathbb{R}^2)$, have geometric multiplicity  equal to one''}. Thus, condition~\eqref{H2} is a Fattorini-Hautus criterium for the boundary approximate controllability of the linear system~\eqref{PFSylinear} (see~\cite{Fattorini}).
\end{remark}

In this work, we will also analyze the null controllability properties of system~\eqref{PFSylinear}. In this sense, one has:

\begin{theorem} \label{CNull_}
Let us us fix $T>0$ and consider $\xi$, $\rho$ and $\tau$, positive real numbers satisfying~\eqref{H2} and
	\begin{equation}\label{H1}
\xi \neq \dfrac{1}{j^2}\dfrac{\rho}{\tau},\quad \forall j\geq1.
	\end{equation}
Then, system~\eqref{PFSylinear} is exactly controllable to zero in $H^{-1} (0, \pi; \R^2)$ at time $T>0$. Moreover, there exist  two positive constants $C_0$ and $M$, only depending on $\xi$, $\rho$ and $\tau$, such that for any $T  > 0$, there is a bounded linear operator 
	$$ 
\mathcal{C}_T^{(0)}: H^{-1}(0,\pi;\mathbb{R}^2) \rightarrow L^2(0,T)
	$$
satisfying
	\begin{equation}\label{CT0}
\| \mathcal{C}_T^{(0)} \|_{\mathcal{L}(H^{-1}(0,\pi;\mathbb{R}^2),L^2(0,T))} \leq C_0 \, e^{M/T },
	\end{equation}
and such that the solution 
	$$
y =( \theta , \phi) \in L^2 (Q_T; \R^2) \cap C^0([0,T];H^{-1}(0,\pi;\mathbb{R}^2)) 
	$$ 
of system~\eqref{PFSylinear} associated to $y_0 = ( \theta_0 , \phi_0) \in H^{-1}(0,\pi;\mathbb{R}^2)$ and $v=\mathcal{C}_T^{(0)}(y_0)$ satisfies $y(\cdot , T)=0$. 
	\end{theorem}
%


\begin{remark}
From the results stated in~\cite{AK-B-D-K} and~\cite{GB-PG}, it is well known that the linear system~\eqref{PFSylinear} is approximate and null controllable at any time $T>0$ and any positive $\xi$, $\rho$ and $\tau$, when the scalar control $v \in L^2(Q_T)$ acts on the temperature equation of~\eqref{PFSy} as a right-hand side source supported on an open subset $\omega$ of the domain. 
These distributed controllability results are independent of condition~\eqref{H2} and only use the cascade structure of system~\eqref{PFSylinear}.
Nevertheless, this cascade structure is not enough when one deals with the boundary controllability of non-scalar problems (see for example~\cite{MGBBoundC}, \cite{AK-B-GB-dT-2011}, \cite{MGBminT}, ... ).
Again, the approximate and null controllability results stated in Theorems~\ref{CAproximada_} and~\ref{CNull_} show the different nature of the controllability problem of scalar or non-scalar parabolic systems.
\end{remark}


\begin{remark}\label{r1.3}
Theorem~\ref{CNull_} also provides an estimate of the cost of the control for system~\eqref{PFSylinear} that drives the system from an initial datum $y_0 = (\theta_0, \phi_0) \in H^{-1} (0, \pi; \R^2)$ to the equilibrium at time $T>0$. To be precise, under assumption~\eqref{H2} and~\eqref{H1}, Theorem~\ref{CNull_} implies that the set
	$$
\mathcal{Z}_T(y_0) := \{ v \in L^2(0,T) : y = (\theta, \phi) \hbox{ solution of }\eqref{PFSylinear}  \hbox{ associated to $y_0$ satisfies } y(\cdot, T) = 0\},
	$$
is nonempty for any $T>0$ and any $y_0 = (\theta_0 , \phi_0) \in H^{-1} (0, \pi; \R^2)$. We can then define the control cost for system~\eqref{PFSylinear} as
	\begin{equation*}
\mathcal{K}(T)= \sup_{\|y_0\|=1} \left( \inf_{v\in\mathcal{Z}_{T}(y_0)} \| v \|_{ L^2(0,T)} \right), \quad \forall T>0.
	\end{equation*}

Observe that as a direct consequence of Theorem~\ref{CNull_} and inequality~\eqref{CT0}, we can obtain the following estimate of this cost for system~\eqref{PFSylinear} at time $T>0$:
	\begin{equation}\label{cost}
\mathcal{K}(T) \leq C_0 \, e^{\frac MT},\quad  \forall T > 0,
	\end{equation}
where $C_0$ and $M$ are positive constants only depending on the parameters in system~\eqref{PFSylinear} (see~\cite{Seidman} and \cite{FC-E} for similar results in the scalar parabolic framework).
\end{remark}


\begin{remark}
As said before, condition~\eqref{H2} is equivalent to the simplicity of the spectrum of $L$ and $L^*$. We will see in Proposition~\ref{propert} that condition~\eqref{H1} implies a stronger property of the spectra of $L$ and $L^*$: If we denote $\{ \Lambda_k \}_{k \ge 1} \subset (0, \infty)$ the sequence of eigenvalues of the operator $L$, with $\Lambda_k \le \Lambda_{k+1}$ for any $k \ge 1$, then, there exist $\delta > 0$ and an integer $ q \ge 1$ such that
	\begin{equation}\label{separability}
\left| \Lambda_{k} - \Lambda_n  \right| \ge \delta \left| k^2 - n^2 \right| , \quad \forall k,n \in \N, \quad | k - n | \ge q .
	\end{equation}
This gap condition for the spectrum of $L$ is crucial to prove the null controllability at any positive time $T$ of system~\eqref{PFSylinear} with control cost satisfying the estimate~\eqref{cost} for positive constants $C_0$ and $M$ only depending on $\xi$, $\rho$ and $\tau$ (for similar results, see~\cite{FaRu} and~\cite{Mil04}).

In the case in which assumption~\eqref{H1} does not hold, that is to say, if for some integer $j \ge 1$ one has
	\begin{equation*}
\xi = \dfrac{1}{j^2}\dfrac{\rho}{\tau},
	\end{equation*}
then, the eigenvalues of $L$ (and $L^*$) concentrate (see Remark~\ref{r3.1}) and the gap condition~\eqref{separability} is not valid. In fact, one has
	$$
\displaystyle \inf_{k,\ell\geq1, k \not= \ell}|\Lambda_{k} - \Lambda_{\ell} | = 0 . 
	$$
In~\cite{MGBminT}, the authors proved that when the eigenvalues $\{ \Lambda_k \}_{k \ge 1}$ of the operator $L = - D \partial_{xx} + A$ concentrate, the controllability problem for system~\eqref{vectPFSy} ($ f \equiv 0 $) has a minimal time $T_0 \in [0, \infty]$ of null controllability which is related to the condensation index of the sequence. Even in the case $T_0 = 0$ (and therefore, system~\eqref{PFSylinear} is null controllable for any $T > 0$), without the separability assumption~\eqref{separability}, providing an estimate of the control cost $\mathcal K (T)$ with respect to $T >0$ is an open problem.
\end{remark}


\begin{remark}
In order to prove Theorem~\ref{CNull_} we will use the moment method, introduced in~\cite{FaRu} in the framework of  the boundary controllability of the one-dimensional scalar heat equation. To this end, we will carry out an analysis of the properties of the eigenvalues of $L$ which will imply inequality~\eqref{CT0} and estimate~\eqref{cost} for the control cost of system~\eqref{PFSylinear}. These two inequalities are essential in the proof of the controllability property of the nonlinear system~\eqref{PFSy}.
\end{remark}

Let us now present the local exact controllability result to the trajectory $(0,c)$ ($c = \pm 1$) for the nonlinear system~\eqref{PFSy}. This is our third main result. One has:

%
	\begin{theorem}\label{NullC}
Let us consider $\xi,\tau$ and $\rho$ three  positive numbers satisfying~\eqref{H2} and~\eqref{H1}, 
and let us fix $T>0$ and $c=-1$ or $c=1$. Then, there exists $\varepsilon>0$ such that, for any $(\tilde{\theta}_0,\tilde{\phi}_0)\in H^{-1}(0,\pi)\times (c + H_0^1(0,\pi))$ fulfilling
	\begin{equation}\label{epsi}
\|\tilde{\theta}_0\|_{H^{-1}} + \|\tilde{\phi}_0 - c\|_{H_0^1} \le \varepsilon,
	\end{equation}
there exists $v\in L^2(0,T)$ for which system \eqref{PFSy} has a unique solution 
	$$
(\tilde\theta,\tilde{\phi})\in \left[L^2(Q_T)\cap C^0([0,T]; H^{-1}(0,\pi; \R^2)) \right] \times C^0(\overline{Q}_T)
	$$
which satisfies~\eqref{exact}.
	\end{theorem}
%
%

Theorem~\ref{NullC} establishes a local exact boundary controllability result at time $T$ for the nonlinear system~\eqref{PFSy} when the parameters $\xi$, $\rho$ and $\tau$ satisfy~\eqref{H2} and~\eqref{H1}. Similar distributed controllability results\footnote{The distributed control acts as a source in the temperature equation.} have been proved in the $N$-dimensional case, without any assumption on the parameters, using the cascade structure of the system (see~\cite{AK-B-D-K} and~\cite{GB-PG}). As in the linear case~\eqref{PFSylinear}, this cascade structure is not enough for dealing with the boundary controllability of system~\eqref{PFSy}.

We end the presentation of our main results with some remarks. 


\begin{remark}
Following~\cite{tucsnack}, Theorem~\ref{NullC} will be proved using a point-fixed strategy. The key point in its proof will be a boundary null controllability result for the non-homogeneous system~\eqref{vectPFSy} when the function $f$ is in an appropriate weighted-$L^2$ space.
In turn, this null controllability result for~\eqref{vectPFSy} will use in a crucial way the estimates~\eqref{CT0} and~\eqref{cost}.
\end{remark}


\begin{remark}
The main results established in this paper only deal with the boundary controllability of linear or nonlinear systems in space dimension one. This restriction is mainly due to the fact that in its proofs we will use the moment method.
In general, the boundary controllability of parabolic systems in higher space dimension remains widely open and only some partial answers are known in the linear setting. 
To our knowledge, the only results on this issue are those of~\cite{AB}, \cite{AB-L} and~\cite{BBGB}. 
In the two first articles, the results for parabolic systems are deduced from the study of the boundary control problem of two coupled wave equations using transmutation techniques. 
As a result they rely on some geometric constraints on the control domain. In~\cite{BBGB}, the author characterize the boundary null-controllability of  system~\eqref{genralsystems} in the linear case ($B \equiv 0$ and $F \equiv 0$) when $\Omega$ is a cylindrical domains of the form $\Omega =(0,\pi)\times \Omega_2$  ($\Omega_2$ is a smooth domain of $\R^{N-1}$, $N>1$) and $\Gamma_0 := \{ 0 \} \times \omega_2$ ($\omega_2 $ is an open subset of $\Omega_2$) without imposing geometric constraints on $\omega_2 $. 
It is important to highlight that these results use that the diffusion matrix $D $ is a multiple of the identity matrix. 
The boundary controllability of systems~\eqref{PFSy} and~\eqref{PFSylinear} in the $N$-dimensional case is completely open.
\end{remark}

The rest of the paper is organized as follows: In Section~\ref{s2}, we give some existence and uniqueness results for the linearized versions of the phase-field system~\eqref{PFSy} and we recall some known results on existence and bounds on biorthogonal families to complex exponentials. Section~\ref{s3} is devoted to studying the spectral properties of the parabolic operators $L$ and $L^*$ given in~\eqref{LL*}. In Section~\ref{s4} we prove the controllability results for the linear problem~\eqref{PFSylinear}: In Subsection~\ref{s4.1} we prove the approximate controllability result at time $T$ for system~\eqref{PFSylinear} (Theorem~\ref{CAproximada_}) and in Subsection~\ref{s4.2} the corresponding null controllability result (Theorem~\ref{CNull_}). Theorem~\ref{NullC} is proved in Section~\ref{s5}. Before (Subsection~\ref{s5.1}), we prove a null controllability result for the non-homogeneous system~\eqref{vectPFSy} when $f$ belongs to appropriate spaces. As a consequence, we provide a proof of Theorem~\ref{NullC} in Subsection~\ref{s5.2}. We finish this paper with two appendices. In Appendix~\ref{APPX1_}, we prove the existence and uniqueness result for the linearized system~\eqref{vectPFSy} and for its backward formulation. In Appendix~\ref{APPX2_} we give some additional results on the null controllability of the phase-field system~\eqref{PFSy}, that is to say, we deal with the case $c=0$ (see~\eqref{exact}).


\section{Preliminary results}\label{s2}

\setcounter{equation}{0}

In this paper we will use the following notations for norms. If $X$ is a Banach space, the norms of the spaces $L^2 (0,T; X)$ and $C^0 ([0,T]; X)$ will be respectively denoted by $\| \cdot \|_{L^2 (X)}$ and $\| \cdot \|_{C^0 (X)}$. We will also work with the spaces $L^2(0,\pi; \R^2)$, $H_0^1(0,\pi; \R^2)$ and $H^{-1}(0,\pi; \R^d)$, with norms denoted by $\| \cdot \|_{L^2}$, $\| \cdot \|_{H_0^1}$ and $\| \cdot \|_{H^{-1}}$. On the other hand, we will use $\langle \cdot\,, \cdot \rangle$ as notation for the usual duality pairing between $H^{-1} (0,1; \R^2)$ and $H_0^1 (0,1; \R^2)$.

Finally, throughout the paper $C$ will stand for a generic positive constant that only depends on the coefficients $\xi$, $\tau$ and $\rho$ in system~\eqref{PFSy}, whose value may change from one line to another. Frequently, we will use the notation $C_T$ when it is convenient to specify the dependence of the generic constant with respect to the final time $T$.

In this section we will give some results related to the existence, uniqueness and continuous dependence with respect to the data of the linear problem~\eqref{vectPFSy}. To this aim, let us consider the linear backwards in time problem:
	\begin{equation}\label{AdjvectPFSy}
	\left\{
	\begin{array}{ll}
- \varphi_t - D^* \varphi_{xx} + A^* \varphi = g & \mbox{in } Q_T,\\
	\noalign{\smallskip}
\varphi(0,\cdot)=\varphi(\pi,\cdot)=0 & \mbox{on } (0,T), \\
	\noalign{\smallskip}
\varphi(\cdot, T) = \varphi_0 & \mbox{in }  (0,\pi),\\
	\end{array}
	\right. 
	\end{equation}
where $D$ and $A$ are given in~\eqref{ADB} and $\varphi_0$ and $g$ are functions in appropriate spaces.

Let us start with a first result on existence and uniqueness of strong solutions to system~\eqref{AdjvectPFSy}. One has:
\begin{proposition}\label{adjexist}
Let us assume that $\varphi_0 \in H_0^1(0,\pi;\mathbb{R}^2)$ and $g\in L^2(Q_T;\mathbb{R}^2)$. Then, system~\eqref{AdjvectPFSy} has a unique strong solution 
	$$
\varphi\in C^0([0,T];H_0^1(0,\pi;\mathbb{R}^2))\cap L^2(0,T;H^2(0,\pi;\mathbb{R}^2) \cap H_0^1(0,\pi;\mathbb{R}^2)).
	$$ 
In addition, there exists a positive constant $C$, only depending on $D$ and $A$, such that
	\begin{equation}\label{acot1}
\|\varphi\|_{C^0(H_0^1)} +  \|\varphi\|_{L^2(H^2\cap H_0^1)} \leq e^{CT} \left(\|g\|_{L^2(L^2)} + \|\varphi_0\|_{H_0^1} \right).
	\end{equation}
\end{proposition}

In view of Proposition~\ref{adjexist}, we can define solution by transposition to system~\eqref{vectPFSy}.

%
%
\begin{definition}\label{deft_}
Let $y_0\in H^{-1}(0,\pi;\mathbb{R}^2)$, $v\in L^2(0,T)$ and $f\in L^2(Q_T;\mathbb{R}^2)$  be given. It will be said that $y \in L^2(Q_T;\mathbb{R}^2)$ is a solution by transposition to~\eqref{vectPFSy} if, for each $g\in L^2(Q_T;\mathbb{R}^2)$, one has
	\begin{equation}\label{transp}
\displaystyle\intdoble_{Q_T}y\cdot g \, dx\, dt = \left\langle y_0,\varphi(\cdot,0) \right\rangle - \int_0^T B^*D^* \varphi_x(0,t)v(t) \, dt + \intdoble_{Q_T}	f \cdot \varphi \, dx\, dt, 
	\end{equation}
where $\varphi \in C^0([0,T];H_0^1(0,\pi;\mathbb{R}^2))\cap L^2(0,T;H^2 (0,\pi;\mathbb{R}^2)\cap H_0^1(0,\pi;\mathbb{R}^2))$ is the solution of \eqref{AdjvectPFSy} associated to $g$ and $\varphi_0=0$ (recall that $\langle \cdot\,, \cdot \rangle$ stands for the usual duality pairing between $H^{-1} (0,1; \R^2)$ and $H_0^1 (0,1; \R^2)$).
\end{definition}

With this definition we have: 
%
%
\begin{proposition}\label{linearexist}
Let us assume that $y_0 = (\theta_0, \phi_0) \in H^{-1}(0,\pi;\mathbb{R}^2)$, $v\in L^2(0,T)$ and $f\in L^2(Q_T;\mathbb{R}^2)$. Then, system~\eqref{vectPFSy} admits a unique solution by transposition $ y = (\theta, \phi)$ that satisfies
	\begin{equation*}
	\left\{
	\begin{array}{l}
\displaystyle y \in L^2(Q_T;\mathbb{R}^2)\cap C^0([0,T];H^{-1}(0,\pi;\mathbb{R}^2)), \quad y_t \in L^2 (0,T; ( H^2(0,\pi;\mathbb{R}^2) \cap H_0^1(0,\pi;\mathbb{R}^2) )' ), \\
	\noalign{\smallskip}
y_t - D y_{xx} + A y = f  \mbox{ in } L^2 (0,T; ( H^2(0,\pi;\mathbb{R}^2) \cap H_0^1(0,\pi;\mathbb{R}^2) )' ),\\
	\noalign{\smallskip}
y (\cdot , 0 ) = y_0 \mbox{ in }H^{-1}(0,\pi;\mathbb{R}^2) ,\\
	\noalign{\smallskip}
	\end{array}
	\right.
	\end{equation*}
and
	\begin{equation}\label{acot2}
\| y \|_{L^2(L^2)} + \| y \|_{C^0(H^{-1})} + \| y_t \|_{L^2 ((H^2 \cap H_0^1)')} \leq Ce^{CT} \left( \| y_0 \|_{H^{-1}} + \|v\|_{L^2 (0,T)} + \|f\|_{L^2(L^2)} \right) ,
	\end{equation}
for a constant $C>0$ only depending on the parameters $\xi$, $\rho$ and $\tau$ in system~\eqref{vectPFSy}. Moreover
\begin{enumerate}[(a)]
\item If $\phi_0\in L^2(0,\pi)$, then $\phi\in L^2(0,T;H_0^1(0,\pi))\cap C^0([0,T];L^2(0,\pi))$ and, for a new constant $C>0$ (only depending on $\xi$, $\rho$ and $\tau$), one has
	\begin{equation}\label{acot3}
\|\phi\|_{L^2(H_0^1)} + \|\phi\|_{C^0(L^2)} \leq C \left( \|y\|_{L^2(L^2)} + \| \phi_0 \|_{L^2} + \|f\|_{L^2(L^2)} \right).
	\end{equation}
\item If $\phi_0\in H_0^1(0,\pi)$, then $\phi\in L^2(0,T;H^2(0,\pi) \cap H_0^1(0,\pi))\cap C^0([0,T];H_0^1(0,\pi))$ and, for a new constant $C>0$ (only depending on $\xi$, $\rho$ and $\tau$), one has
	\begin{equation}\label{acot4}
\|\phi\|_{L^2(H^2\cap H_0^1)} + \|\phi\|_{C^0(H_0^1)} \leq C \left(  \|y\|_{L^2(L^2)} + \| \phi_0 \|_{H_0^1} + \|f\|_{L^2(L^2)} \right),
	\end{equation}
and, in particular, $y = (\theta,\phi)\in L^2(Q_T)\times C^0(\overline Q_T)$.
\end{enumerate}
	\end{proposition}
%
%
%
One can prove Propositions~\ref{adjexist} and~\ref{linearexist} using, for instance, the well-known Galerkin method. For the sake of completeness we present an idea of the proof of this two propositions in Appendix~\ref{APPX1_}.

Observe that, when $g = 0$, the backward problem~\eqref{AdjvectPFSy} is the corresponding adjoint system to~\eqref{PFSylinear}:
	\begin{equation}\label{adjoint}
	\left\{
	\begin{array}{ll}
- \varphi_t - D^* \varphi_{xx} + A^* \varphi = 0 & \mbox{in } Q_T,\\
	\noalign{\smallskip}
\varphi(0,\cdot)=\varphi(\pi,\cdot)=0 & \mbox{on } (0,T), \\
	\noalign{\smallskip}
\varphi(\cdot, T) = \varphi_0 & \mbox{in }  (0,\pi).\\
	\end{array}
	\right. 
	\end{equation}

The controllability properties of system~\eqref{PFSylinear} can be characterized in terms of appropriate properties of the solutions to~\eqref{adjoint}. In order to provide these characterizations, we need a new result which relates the solutions of systems~\eqref{PFSylinear} and~\eqref{adjoint}. One has:

%
%
	\begin{proposition} \label{premomprob}
Let us consider $y_0=(\theta_0,\phi_0)\in H^{-1}(0,\pi;\mathbb{R}^2)$ and $v\in L^2(0,T)$. Then, the solution $y= (\theta,\phi)$ of system~\eqref{PFSylinear} associated to $y_0$ and $v$, and the solution $\varphi$ of the adjoint system~\eqref{adjoint} associated to $\varphi_0 \in H_0^{1}(0, \pi; \R^2)$ satisfy
	\begin{equation}\label{pmpeq}
\displaystyle \int_0^T B^*D^*\varphi_x(0,t) v(t) \, dt = \langle y(\cdot,T), \varphi_0 \rangle - \langle y_0, \varphi(\cdot,0) \rangle . 
	\end{equation}
	\end{proposition}
%
%

%
\begin{proof}
	The proof is a consequence of Proposition~\ref{linearexist}. Observe that is enough to prove that~\eqref{pmpeq} holds under the regularity assumption $y_0\in C_0^1(0,\pi;\mathbb{R}^2)$ and $v \in C_0^1([0, \pi])$. Indeed, using density arguments,  the estimates of Proposition~\ref{linearexist} and the linearity of~\eqref{PFSylinear}, it follows that the identity~\eqref{pmpeq} is valid for all $y_0\in H^{-1}(0,\pi;\mathbb{R}^2)$ and $v \in L^2(0,T)$.
	
On the other hand, when $y_0\in C_0^1(0,\pi;\mathbb{R}^2)$, $v \in C^1([0, \pi])$ and $\varphi_0 \in H_0^1(0,\pi;\mathbb{R}^2)$, after some integrations by parts, it is not difficult to prove that the corresponding solution $y$ of~\eqref{PFSylinear} and $\varphi$, solution of the adjoint system~\eqref{adjoint}, satisfy equality~\eqref{pmpeq}. This ends the proof.
%
%
\end{proof}

One important consequence of the previous result is the characterization of the approximate and null controllability properties of the linear system~\eqref{PFSylinear} in terms of suitable properties of the solutions of the adjoint system~\eqref{adjoint}. One has:

\begin{theorem}\label{equivcontrol}
Let us consider $T>0$. Then, 
\begin{enumerate}
\item System~\eqref{PFSylinear} is approximately controllable at time $T>0$ if and only if the following unique continuation property holds:
\begin{quote}
\begin{center}
``Let $\varphi_0 \in H_0^1(0,\pi;\mathbb{R}^2)$ be given and let $\varphi$ be the corresponding solution of the adjoint problem~\eqref{adjoint}. Then, if $B^*D^*\varphi_x(0,t)=0$ on $(0,T)$, one has $\varphi_0=0$ in $(0,\pi)$.''
\end{center}
\end{quote}
\item System~\eqref{PFSylinear} is null controllable at time $T$ if and only if there exists a constant $C_T > 0$ such that, for any $\varphi_0 = (\theta_0, \phi_0) \in H_0^{1}(0, \pi; \R^2)$, the corresponding solution of~\eqref{adjoint} satisfies the observability inequality
	$$
\| \varphi (\cdot, T) \|_{H_0^1}^2 \le C_T \int_0^T \left| B^*D^*\varphi_x(0,t) \right|^2 \, dt.
	$$
\end{enumerate}
\end{theorem}

This result is well known. For a proof see, for instance~\cite{C},~\cite{TW} and~\cite{Z}.

\begin{remark}
The constant $C_T$ appearing in the observability inequality for the adjoint system~\eqref{adjoint} is closely related to the cost $\mathcal K (T)$ for system~\eqref{PFSylinear} (see Remark~\ref{r1.3}). To be precise, if the observability inequality holds, then $\mathcal Z(T) \not= \emptyset$, for any $y_0 = (\theta_0 , \phi_0) \in H^{-1} (0, \pi; \R^2)$, and
	$$
\mathcal K (T) \le \sqrt{C_T} .
	$$
On the other hand, assume that $\mathcal Z(T) \not= \emptyset$, for any $y_0 = (\theta_0 , \phi_0) \in H^{-1} (0, \pi; \R^2)$, and define $\mathcal K (T)$ as in Remark~\ref{r1.3}. Then, the previous observability inequality for~\eqref{adjoint} holds with $C_T = \mathcal K (T)^2$.

For a proof of the previous properties, see for example~\cite{C} (see Theorem 2.44, p.~56),~\cite{Z} or~\cite{TW}.
\end{remark}

%

We will finish this section giving two known results which will be used later. They are related to the existence and bounds of biorthogonal families to real exponentials. One has:

%
%
\begin{lemma}\label{bioex}
Let us consider a sequence $\{\Lambda_k\}_{k\geq1} \subset \mathbb{R}_+$ satisfying $\Lambda_k\neq\Lambda_n$, for any $k,n\in\mathbb{N}$ with $k\neq n$, and
	\begin{equation}\label{serconv}
\displaystyle \sum_{k\geq1} \dfrac{1}{\Lambda_k} < \infty.
	\end{equation}
Then, there exists a family $\{q_k\}_{k\geq1}\subset L^2(0,T)$ biorthogonal  to $\{ e^{-\Lambda_kt} \}_{k\geq1}$, i.e., a family $\{q_k\}_{k\geq1}$ in $ L^2(0,T)$ such that 
	\begin{equation*}
\displaystyle\int_0^T  q_k(t)e^{-\Lambda_jt} dt = \delta_{kj},\quad \forall k,j\geq1.
	\end{equation*}
	\end{lemma}
%

We also have:

\begin{lemma}\label{1.5M}
Let us consider a sequence $\{\Lambda_k\}_{k\geq1} \subset \mathbb{R}_+$ such that $\Lambda_k\neq\Lambda_n$, for any $k,n\in\mathbb{N}$ with $k\neq n$. Let us also assume that there exist an integer $q \ge 1$ and positive constants $p$, $\delta$ and $\alpha$ such that
	\begin{equation}\label{56}
	\left\{
	\begin{array}{l}
\left| \Lambda_k - \Lambda_n \right| \geq \delta \left| k^2-n^2 \right| ,\quad \forall k,n\in \N, \ \left| k-n \right|\geq q, \\
	\noalign{\smallskip}
\inf_{k\neq n,\ \left| k-n \right|<q} \left| \Lambda_k - \Lambda_n \right|>0,
	\end{array}
	\right.
	\end{equation}
and
	\begin{equation}\label{56_2}
\left| p\sqrt{r}-\mathcal{N}(r) \right|\leq\alpha,\quad \forall r>0. 
	\end{equation}
(In~\eqref{56_2}, $\mathcal{N}(r)$ is the counting function associated to $\{\Lambda_k\}_{k\geq1}$, defined by $\mathcal{N}(r)= \#\{ k :  \Lambda_k\leq r \}$). Then, there exists $\widetilde T_0>0$ such that, for any $T\in(0,\widetilde T_0)$, we can find a family $\{q_k\}_{k\geq1}\subset L^2(0,T)$ biorthogonal to $\{ e^{-\Lambda_kt} \}_{k\geq1}$ which in addition satisfies
	\begin{equation*}
\|q_k\|_{L^2(0,T)} \leq C e^{C\sqrt{\Lambda_k} + \frac{C}{T}},\quad \forall k \geq 1,
	\end{equation*}
for a positive constant $C$ independent of $T$.
\end{lemma}

A proof of Lemma~\ref{bioex} can be found in~\cite{FaRu} and~\cite{AK-B-GB-dT-2011}. On the other hand, Lemma~\ref{1.5M} is a particular case of a more general result proved in~\cite{BBGB} (see Theorem~1.5 in pages~2974--2975).

%

\section{Spectral properties of the operators $L$ and $L^*$}\label{s3}
\setcounter{equation}{0}
Let us consider the vectorial operators $L$ and $L^*$ given in~\eqref{LL*}, with domains 
	$$
D(L)=D(L^*)=H^2(0,\pi;\mathbb{R}^2) \cap H_0^1(0,\pi;\mathbb{R}^2).
	$$ 
This section will be devoted to giving some spectral properties of the operators $L$ and $L^*$ which will be used below. Recall that the matrices $D$ and $A$ are given in~\eqref{ADB}.

In what follows, for simplicity, we will use the notation
	\begin{equation}\label{rkquant}
r_k:= \sqrt{ \dfrac{\xi\rho}{\tau} k^2 + \left(\dfrac{\rho+1}{2\tau}\right)^2 } , \quad \forall k \ge 1.
	\end{equation}
On the other hand, it is well-known that the operator $- \partial_{xx}$ with homogeneous Dirichlet boundary conditions admits a sequence of positive eigenvalues, given by $\{k^2 \}_{k\geq1}$, and a sequence of normalized eigenfunctions $\{\eta_k\}_{k\geq1}$, which is a Hilbert basis of $L^2(0,\pi)$, given by
	\begin{equation}\label{sin}
\eta_k (x) = \sqrt{\dfrac{2}{\pi}}\sin kx, \quad x \in (0,\pi).
	\end{equation}
With the previous notation, we have the following result:

%
%
	\begin{proposition}\label{Lspec}
Let us consider the operators $L$ and $L^*$ given in~\eqref{LL*} (the matrices $D$ and $A$ are given in~\eqref{ADB}). Then, 
\begin{enumerate}
\item The spectra of $L$ and $L^*$ are given by $\sigma(L) = \sigma (L^*)= \{\lambda_k^{(1)}, \lambda_k^{(2)} \}_{k\geq1}$ with
	\begin{equation}\label{lamk12}
\lambda_k^{(1)} = \xi k^2 + \dfrac{\rho + 1}{2\tau} - r_k ,\quad \lambda_k^{(2)} = \xi k^2 + \dfrac{\rho + 1}{2\tau} + r_k, \quad \forall k \ge 1,			
	\end{equation}
where $r_k$ is given in~\eqref{rkquant}. 

\item For each $k\geq1$, the eigenspaces of ${L}$ (resp., ${L}^*$) corresponding to ${\lambda}_k^{(1)}$ and ${\lambda}_k^{(2)}$ are respectively generated by

	\begin{equation}\label{pk12}
\Psi_k^{(1)} = \frac{1}{4\sqrt{\tau r_k}}\left(
	\begin{array}{c}
1-\rho+2\tau r_k	 \\
	\smallskip
4
	\end{array}
	\right)
\eta_k,
\quad
\Psi_k^{(2)} = \frac{1}{4\sqrt{\tau r_k}}
	\left(
	\begin{array}{c}
1-\rho-2\tau r_k	\vspace{0.2cm}\\
	\smallskip
4
	\end{array}
	\right)
\eta_k, 
	\end{equation}
(resp., 
	\begin{equation}\label{pk1_22}
\Phi_k^{(1)} = \frac{1}{4\sqrt{\tau r_k}}
	\left(
	\begin{array}{c}
4 \\
	\smallskip
\rho-1+2\tau r_k
	\end{array}
	\right)
\eta_k,
\quad
\Phi_k^{(2)} = \frac{-1}{4\sqrt{\tau r_k}}
	\left(
	\begin{array}{c}
4 \\
	\smallskip
\rho-1-2\tau r_k
	\end{array}
	\right)
	\eta_k  \hbox{)}.
\end{equation}
\end{enumerate}
\end{proposition}
%
%
	\begin{proof}
We will prove the result for the operator $L$. The same reasoning provides the proof for its adjoint $L^*$.

We look for a complex $\lambda$ and a function $\varphi \in H^2(0, \pi; \C^2) \cap H_0^1(0, \pi; \C^2 )$ such that $ \varphi \not\equiv 0$ and $ L ( \varphi) = \lambda \varphi $. Using that the function $\eta_k$ is the normalized eigenfunction of the Dirichlet-Laplace operator in $(0, \pi)$ associated to the eigenvalue $k^2$, we can find $\varphi$ as
	$$
\varphi (x) = \sum_{n \ge 1} a_n \eta_n (x), \quad \forall x \in (0, \pi),
	$$
where $\{ a_n \}_{n \ge 1} \subset \C^2$ and, for some $k \ge 1$, $a_k \not= 0$. From the identity $ L ( \varphi) = \lambda \varphi $ we deduce
	$$
\sum_{n \ge 1} \left( n^2 D + A - \lambda I_2 \right)  a_n \eta_n (x) = 0, \quad \forall x \in (0, \pi),
	$$
($I_2 \in \mathcal{L}(\C^2)$ is the identity matrix). From this identity, it is clear that the eigenvalues of the operator $L$ correspond to the eigenvalues of the matrices
	$$
k^2 D + A, \quad \forall k \ge 1.
	$$
and an associated eigenfunction of $L$ is given choosing $a_n = z_k \delta_{kn} $, for any $n \ge 1$, where $z_k \in \mathbb{C}^2$ is an associated eigenvector of $k^2 D + A$, that is to say, $ \Psi_k (\cdot) = z_k \eta_k (\cdot)$.

Taking into account the expression of the characteristic polynomial of $k^2 D + A$:
	$$
p(x) = x^2 - \left( 2\xi k^2 + \dfrac{\rho+1}{\tau} \right)x + \xi^2 k^4 +\dfrac{\xi}{\tau}k^2, \quad k \ge 1,
	$$
a direct computation provides the formulae~\eqref{lamk12} and~\eqref{pk12} as eigenvalues and associated eigenfunctions of the operator $L$. This finishes the proof. 
	\end{proof}
%

%

Let us now analyze some properties of the eigenvalues and eigenfunctions of the operators $L$ and $L^*$. These properties will be used below. We start with some properties of the sequences $ \{\lambda_k^{(1)}\}_{k\geq1}$ and $\{\lambda_k^{(2)}\}_{k\geq1}$. One has

	\begin{proposition}\label{propert}
	Under the assumptions of Proposition~\ref{Lspec}, the following properties hold:
\begin{enumerate}
\item[(P1)] $ \{\lambda_k^{(1)}\}_{k\geq1}$ and $\{\lambda_k^{(2)}\}_{k\geq1}$ (see~\eqref{lamk12}) are increasing sequences satisfying
	\begin{equation*}
0 < \lambda_{k}^{(1)} < \lambda_k^{(2)},\quad \forall k\geq1.		
	\end{equation*}
\item[(P2)] The spectrum of $L$ and $L^*$ is simple, i.e., $\lambda_{k}^{(2)} \neq \lambda_{\ell}^{(1)}$, for all $k,\ell\geq 1$ if and only if the parameters $\xi$, $\rho$ and $\tau$ satisfy condition~\eqref{H2}.

\item[(P3)] Assume that the parameters $\xi$, $\rho$ and $\tau$ satisfy~\eqref{H1}, i.e., there exists $j\geq0$  such that 
	\begin{equation}\label{xi}
\frac{1}{(j+1)^2}\frac{\rho}{\tau}<\xi<\frac{1}{j^2}\frac{\rho}{\tau}.
	\end{equation}
Then, there exists an integer $k_0=k_0(\xi,\rho,\tau, j) \geq 1$ and a constant $C= C(\xi,\rho,\tau, j) >0$  such that
	\begin{equation}\label{orden}
\lambda_{k+j}^{(1)} < \lambda_{k}^{(2)} < \lambda_{k+1+j}^{(1)} < \lambda_{k+1}^{(2)}<\cdots , \ \forall k\geq k_0, \ \hbox{and} \ \min_{k \ge k_0} \left\{ \lambda_k^{(2)} - \lambda_{k+j}^{(1)} , \lambda_{k+j+1}^{(1)} - \lambda_k^{(2)} \right\} > C. 
	\end{equation}
\item[(P4)]  Assume now that the parameters $\xi$, $\rho$ and $\tau$ satisfy~\eqref{H2} and~\eqref{H1}. Then, one has: 
	\begin{equation}\label{differ}
\displaystyle \inf_{k,\ell\geq1}|\lambda_{k}^{(2)}-\lambda_{\ell}^{(1)}|>0 ,
	\end{equation}
and there exists a positive integer $k_1 \in \N$, depending on $\xi$, $\rho$ and $\tau$, such that
	\begin{equation}\label{nueva}
\min \left\{ \left| \lambda_k^{(1)} - \lambda_\ell^{(1)} \right| , \left| \lambda_k^{(2)} - \lambda_\ell^{(2)} \right|, \left| \lambda_k^{(2)} - \lambda_\ell^{(1)} \right|   \right\} \ge \frac{\xi}{2} |k^2 - \ell^2| , \quad \forall k, \ell \ge 1, \quad |k - \ell | \ge k_1 .
	\end{equation}
		\end{enumerate}
	\end{proposition}

\begin{proof}
Let us start proving property~$(P1)$. From the expressions of $\lambda_k^{(1)}$ and $\lambda_k^{(2)}$ (see~\eqref{lamk12}), we directly get $\lambda_{k}^{(1)} < \lambda_k^{(2)}$ for any $k \ge 1$. On the other hand, using the inequality
	\begin{equation*}
r_k = \sqrt{ \dfrac{\xi\rho}{\tau}k^2 + \left(\dfrac{\rho+1}{2\tau}\right)^2} < \sqrt{\xi^2 k^4 + 2\xi k^2\dfrac{\rho+1}{2\tau} + \left(\dfrac{\rho+1}{2\tau}\right)^2} = \xi k^2 + \dfrac{\rho+1}{2\tau},\quad \forall k\geq 1,
	\end{equation*}
we also deduce $0 < \lambda_k^{(1)}$ for any $k \ge 1$.

Let us now prove that $\{\lambda_k^{(1)}\}_{k\geq1}$ and $\{\lambda_k^{(2)}\}_{k\geq1}$ are increasing sequences. Indeed, 
	\begin{equation*}
	\begin{array}{l}
\displaystyle \lambda_{k+1}^{(1)} - \lambda_k^{(1)} = \xi(2k+1) + \sqrt{ \dfrac{\xi\rho}{\tau}k^2 + \left(\dfrac{\rho+1}{2\tau}\right)^2} - \sqrt{ \dfrac{\xi\rho}{\tau}(k+1)^2 + \left(\dfrac{\rho+1}{2\tau}\right)^2} \\
	\noalign{\smallskip}
\displaystyle \phantom{\displaystyle \lambda_{k+1}^{(1)} - \lambda_k^{(1)}} = \xi(2k+1) - \dfrac{\xi\rho}{\tau} \dfrac{2k+1} { \sqrt{ \frac{\xi\rho}{\tau}k^2 + \left(\frac{\rho+1}{2\tau}\right)^2} + \sqrt{ \frac{\xi\rho}{\tau}(k+1)^2 + \left(\frac{\rho+1}{2\tau}\right)^2} } \\
	\noalign{\smallskip}
\displaystyle \phantom{\displaystyle \lambda_{k+1}^{(1)} - \lambda_k^{(1)}} = \xi(2k+1) \left[ 1 - \dfrac{\rho}{\tau}\dfrac{1}{\sqrt{ \frac{\xi\rho}{\tau}k^2 + \left(\frac{\rho+1}{2\tau}\right)^2} + \sqrt{ \frac{\xi\rho}{\tau}(k+1)^2 + \left(\frac{\rho+1}{2\tau}\right)^2}} \right] \rightarrow  \infty,
	\end{array}
	\end{equation*}
as $k \rightarrow  \infty$. Moreover,
	\begin{equation*}
				\sqrt{\dfrac{\xi\rho}{\tau}k^2 + \left(\dfrac{\rho+1}{2\tau}\right)^2}
				+ 
				\sqrt{ \dfrac{\xi\rho}{\tau}(k+1)^2 + \left(\dfrac{\rho+1}{2\tau}\right)^2}
			\geq
				\dfrac{\rho+1}{2\tau} + \dfrac{\rho+1}{2\tau}
			>
				\dfrac{\rho}{\tau},
	\end{equation*}
which implies $\lambda_{k+1}^{(1)} - \lambda_k^{(1)} > 0$, for any $k\geq1$. Thus, $\{\lambda_k^{(1)}\}_{k\geq1}$ is a positive increasing sequence. Clearly $\{\lambda_k^{(2)}\}_{k\geq1}$ is also a positive increasing sequence and $\lambda_{k+1}^{(2)} - \lambda_k^{(2)}\rightarrow  \infty$, as $k\rightarrow \infty $. This proves property~$(P1)$.

\smallskip

Let us now see property $(P2)$. Using property~$(P1)$, one has that, for any integers $k, \ell \ge 1$ with $\ell \le k$, clearly $ \lambda_{\ell}^{(1)} \le \lambda_{k}^{(1)} < \lambda_{k}^{(2)}$. Therefore, in order to prove the equivalence we can assume
that $\ell > k$. We have
	\begin{equation*}
\lambda_{\ell}^{(1)} - \lambda_{k}^{(2)}  = \dfrac{\xi\rho}{\tau}(\ell^2-k^2)\left( \dfrac{\tau}{\rho} - \dfrac{1}{r_\ell - r_k} \right).
	\end{equation*}
Thus, $\lambda_{k}^{(2)} \neq \lambda_{\ell}^{(1)}$ for any $k, \ell \ge 1$, with $\ell > k$, if and only if 
	\begin{equation*}
r_\ell^2 \neq \left(r_k + \dfrac{\rho}{\tau}\right)^2, \quad \forall k, \ell \ge 1, \quad \ell >k.
	\end{equation*}
From the expression of $r_k$ (see~\eqref{rkquant}) we readily deduce $2r_k> \frac{\rho}{\tau}$ and $\xi \tau (\ell^2-k^2) -\rho  + 2\tau r_k > 0$ ($\ell > k$). So,
	\begin{equation*}
	\left\{
	\begin{array}{l}
\displaystyle r_\ell^2 - \left(r_k + \dfrac{\rho}{\tau}\right)^2
			=
\dfrac{\rho}{\tau} \left[ \left( \xi(\ell ^2-k^2) - \dfrac{\rho}{\tau} \right) - 2r_k \right] 
			=
\dfrac{\rho}{\tau^2} \dfrac{ \left( \xi\tau(\ell^2-k^2) - \rho \right)^2 - 4\tau^2r_k^2}{\xi \tau (\ell^2 - k^2) -\rho  + 2\tau r_k}
			\\
	\noalign{\medskip}
\displaystyle \phantom{r_\ell^2 - \left(r_k + \dfrac{\rho}{\tau}\right)^2}
			=
\dfrac{\rho}{\tau^2} \dfrac{ \xi^2\tau^2(\ell^2-k^2)^2 -2\xi\tau\rho(\ell^2-k^2)  - \rho^2 - 4\xi\tau\rho k^2 - 2\rho - 1}{\xi \tau (\ell^2 - k^2) -\rho  + 2\tau r_k}
			\\
	\noalign{\medskip}
\displaystyle \phantom{r_\ell^2 - \left(r_k + \dfrac{\rho}{\tau}\right)^2}
			=
				\dfrac{\rho}{\tau^2}\dfrac{ \xi^2\tau^2(\ell^2-k^2)^2 -2\xi\tau\rho(\ell^2+k^2) - 2\rho - 1}{\xi \tau (\ell^2 - k^2) -\rho  + 2\tau r_k},
	\end{array}
	\right.
	\end{equation*}
and we get that $\lambda_{k}^{(2)} \neq \lambda_{\ell}^{(1)}$ for any $k, \ell \ge 1$, with $\ell > k$, if and only if  condition~\eqref{H2} holds. This finishes the proof of property~$(P2)$.

\smallskip

In order to prove property~$(P3)$, we are going to use the expressions
	\begin{equation}\label{bis-lamb}
\lambda_k^{(1)} = \xi k^2 + \dfrac{\rho + 1}{2\tau} - \sqrt{\frac{\xi \rho}{\tau}} \, k - \frac{\epsilon_k}{k} ,\quad \lambda_k^{(2)} = \xi k^2 + \dfrac{\rho + 1}{2\tau} + \sqrt{\frac{\xi \rho}{\tau}} \, k + \frac{\epsilon_k}{k} , \quad \forall k \ge 1,			
	\end{equation}
which can be easily deduced from the expressions of $\lambda_k^{(i)}$, $i=1,2$, and $r_k$ (see~\eqref{lamk12} and~\eqref{rkquant}). In~\eqref{bis-lamb}, $\{ \epsilon_k \}_{k \ge 1}$ is a new positive increasing sequence satisfying
	\begin{equation}\label{new2}
\lim_{k \to \infty} \epsilon_k = \frac 12 \left( \dfrac{\rho + 1}{2\tau} \right)^2 \sqrt{\frac{\tau }{\xi \rho}} \, .
	\end{equation}

Using~\eqref{bis-lamb}, we will prove that, for any $i \ge 1$, the difference $\lambda_{k+i}^{(1)} - \lambda_k^{(2)}$ behaves at infinity as
	\begin{equation}\label{summi}
\lim_{k\rightarrow  \infty} \frac{\lambda_{k+i}^{(1)} - \lambda_k^{(2)}}{\xi i(2k+i)} = 1 - \sqrt{\dfrac{1}{i^2}\dfrac{\rho}{\xi\tau}} \not= 0.
	\end{equation}
Indeed, a simple computation gives
	\begin{equation}\label{new}
\lambda_{k+i}^{(1)} - \lambda_k^{(2)}
				=
\xi i(2k+i) - \sqrt{\frac{\xi \rho}{\tau}} \, \left( 2k +i \right) - \frac{\epsilon_{k + i}}{k + i} -\frac{\epsilon_k}{k}
				=
\xi i(2k+i)\left[ 1 - \sqrt{\dfrac{1}{i^2}\dfrac{\rho}{\xi\tau}} - \widetilde \epsilon_k^{(i)} \right],
	\end{equation}
where $\{ \widetilde \epsilon_k^{(i)}\}_{k \ge 1}$ is a sequence converging to zero. From assumption~\eqref{H1} we can conclude~\eqref{summi}. 

We will obtain the proof of property $(P3)$ from~\eqref{summi}. Observe that assumption~\eqref{H1} implies that the parameters $\xi$, $\rho$ and $\tau$ satisfies~\eqref{xi} for an appropriate integer $j \ge 0$. Therefore, if $j = 0$, then, $\xi > \frac{\rho}{\tau}$ and~\eqref{summi} implies $\lim_{k \to \infty} \left(\lambda_{k+1}^{(1)} - \lambda_k^{(2)} \right) = \infty$. On the other hand, one has $\lim_{k \to \infty} \left(\lambda_{k}^{(2)} - \lambda_k^{(1)} \right) = \lim_{k \to \infty} 2r_k = \infty$. Thus, there exists an integer $k_0 \geq 1$ and a constant $C>0$ such that~\eqref{orden} holds for $j=0$.

If $j \ge 1$, again, the property~\eqref{summi} implies
	\begin{equation*}
			 	\lim_{k \to \infty} \left(\lambda_{k+i}^{(1)} - \lambda_k^{(2)} \right) = - \infty,\quad \mbox{if}\ i\leq j
			 \quad \mbox{and} \quad
			 	\lim_{k \to \infty} \left(\lambda_{k+i}^{(1)} - \lambda_k^{(2)} \right) = \infty,\quad \mbox{if}\ i\geq j+1. \\
	\end{equation*}
We can also conclude the existence of an integer $k_0 \geq 1$ and a positive constant $C$ such that~\eqref{orden} holds. This shows property~$(P3)$.

\smallskip

Let us finalize the proof showing property~$(P4)$. First, inequality~\eqref{differ} is a direct consequence of property~$(P2)$ and~\eqref{orden}.  Secondly, if we take $k, \ell \ge 1$, from~\eqref{bis-lamb}, one deduces:
	$$
	\left\{
	\begin{array}{l}
\displaystyle \left| \lambda_k^{(1)} - \lambda_\ell^{(1)} \right| = \xi \left| k^2 - \ell^2 \right| \left|  1 - \sqrt{\frac{\xi \rho}{\tau}} \frac{1}{k + \ell } - \frac{\epsilon_k}{k (k^2 - \ell^2)} + \frac{\epsilon_\ell}{\ell (k^2 - \ell^2)} \right| \\
	\noalign{\smallskip}
\displaystyle \phantom{\displaystyle \left| \lambda_k^{(1)} - \lambda_\ell^{(1)} \right| } \ge \xi \left| k^2 - \ell^2 \right| \left[ 1 - \left( \sqrt{\frac{\xi \rho}{\tau}} + |\epsilon_k| + | \epsilon_\ell | \right) \frac{1}{k+ \ell } \right].
	\end{array}
	\right.
	$$
Observe that $\{ \epsilon_k \}_{k \ge 1}$ is a convergent sequence and $k + \ell \ge | k - \ell |$, for any $k, \ell \in \N$. Hence, there exists a integer $q_1 \ge 1$ (depending on the parameters of system~\eqref{PFSylinear}) such that
	$$
\displaystyle \left| \lambda_k^{(1)} - \lambda_\ell^{(1)} \right| \ge \frac \xi2 \left| k^2 - \ell^2 \right|, \quad \forall k, \ell \in \N, \quad |k - \ell | \ge q_1. 
	$$
A similar inequality can be deduced for a new $q_2 \in \N$ if we change $\left| \lambda_k^{(1)} - \lambda_\ell^{(1)} \right|$ by $\left| \lambda_k^{(2)} - \lambda_\ell^{(2)} \right|$.

Finally, if we repeat the previous reasoning, we can write
	$$
	\left\{
	\begin{array}{l}
\displaystyle \left| \lambda_k^{(2)} - \lambda_\ell^{(1)} \right| = \xi \left| k^2 - \ell^2 \right| \left|  1 + \sqrt{\frac{\xi \rho}{\tau}} \frac{1}{k - \ell } + \frac{\epsilon_k}{k (k^2 - \ell^2)} + \frac{\epsilon_\ell}{\ell (k^2 - \ell^2)} \right| \\
	\noalign{\smallskip}
\displaystyle \phantom{\displaystyle \left| \lambda_k^{(1)} - \lambda_\ell^{(1)} \right| } \ge \xi \left| k^2 - \ell^2 \right| \left[ 1 - \left( \sqrt{\frac{\xi \rho}{\tau}} + \frac 12 |\epsilon_k| + \frac 12 | \epsilon_\ell | \right) \frac{1}{\left| k - \ell \right|} \right].
	\end{array}
	\right.
	$$
Again, from this inequality we conclude the existence of $q_3 = q_3 (\xi, \rho, \tau) \in \N $ such that  
	$$
\displaystyle \left| \lambda_k^{(2)} - \lambda_\ell^{(1)} \right| \ge \frac \xi2 \left| k^2 - \ell^2 \right|, \quad \forall k, \ell \in \N, \quad |k - \ell | \ge q_3. 
	$$
This proves inequality~\eqref{nueva} if we take $k_1 = \max \{ q_1, q_2, q_3\}$. This completes the proof of $(P4)$ and the proof of the result.
\end{proof}

\begin{remark}\label{r3.0}
From the previous proof we can give more information about condition~\eqref{H2}. To be precise, let us see that, in fact, this condition only has to be checked for a finite number of positive integers $k$ and $\ell$, with $k < \ell$. To this end, let us consider $ j \ge 0$ such that
	\begin{equation}\label{xi-bis}
\frac{1}{(j+1)^2} \frac{\rho}{\tau} < \xi \le \frac{1}{j^2}\frac{\rho}{\tau}, 
	\end{equation}
(if $\frac{\rho}{\tau} < \xi$, then $j=0$). Taking into account that $\{ \epsilon_k \}_{k \ge 1}$ is a positive increasing sequence (see~\eqref{bis-lamb}), identity~\eqref{new} for $i = j$ implies 
	$$
\lambda_{k +j}^{(1)} < \lambda_k^{(2)} \quad \hbox{and} \quad \lambda_{k +j+ 1}^{(1)} < \lambda_{k+1}^{(2)}, \quad \forall k \ge 1.
	$$
On the other hand, using again~\eqref{new} for $i = j$ and the expression of $\{ \widetilde \epsilon_k^{(i)}\}_{k \ge 1}$, we can write
	$$
\lim_{k \to \infty} \left(\lambda_{k+j + 1}^{(1)} - \lambda_k^{(2)} \right) = \infty. 
	$$
Therefore, there exists $k_0 \ge 1$ (only depending on $\xi$, $\rho$ and $\tau$) such that 
	$$
\lambda_{k+j}^{(1)} < \lambda_{k}^{(2)} < \lambda_{k+1+j}^{(1)} < \lambda_{k+1}^{(2)}<\cdots , \quad \forall k\geq k_0.
	$$
In particular, we can assure that $\lambda_{k}^{(2)} \neq \lambda_{\ell}^{(1)}$, for all $k \ge k_0$ and $\ell\geq k_0 + j$. 

In fact, the previous reasoning provides a stronger property for the spectrum of $L$ and $L^*$: \textit{if for some $k, \ell \ge 1$ one has $\lambda_{k}^{(2)} = \lambda_{\ell}^{(1)}$, then $k < k_0$ and $\ell = k + j +1$, with $j \ge 0$ given in~\eqref{xi-bis}}.  That is to say, condition~\eqref{H2}  only has to be checked for $k < k_0$ and $\ell = k + j + 1$.

We can also provide an estimate of $k_0$. Indeed, coming back to formula~\eqref{new} with $i = j + 1$, we infer the following expression of the sequence $\{ \widetilde \epsilon_k^{(j+1)}\}_{k \ge 1}$:
	\begin{equation*}
	\begin{split}
\widetilde \epsilon_k^{(j+1)} &= \frac{1}{\xi ( j + 1) (2k + j + 1)} \left[ \frac{\epsilon_{k + j + 1}}{k + j + 1} + \frac{\epsilon_{k } }{k} \right] \le \frac{\epsilon_{k + j + 1} + \epsilon_{k }  }{\xi ( j + 1) k (2k + j + 1)}  \\
	\noalign{\smallskip}
& < \left( \dfrac{\rho + 1}{2\tau} \right)^2 \sqrt{\frac{\tau }{\xi \rho}} \frac{ 1 }{\xi ( j + 1) k (2k + j + 1)} .
	\end{split}
	\end{equation*}
In the previous inequality we have used that $\{ \epsilon_k \}_{k \ge 1}$ is a positive increasing sequence satisfying~\eqref{new2}. Thus, using once again~\eqref{new}, it is easy to see that if we take $k_0 \ge 1$ such that
	$$
\left( \dfrac{\rho + 1}{2\tau} \right)^2 \sqrt{\frac{\tau }{\xi \rho}} \frac{ 1 }{\xi ( j + 1) k_0 (2k_0 + j + 1)} \le \frac 12 \left( 1 - \sqrt{\dfrac{1}{(j + 1)^2}\dfrac{\rho}{\xi\tau}}  \right),
	$$
we have
	$$
\lambda_{k}^{(2)} < \lambda_{k+1+j}^{(1)} , \quad \forall k\geq k_0 \, .
	$$

\end{remark}

\begin{remark}\label{r3.1}
If condition~\eqref{H1} does not hold, i.e., if for some integer $j \ge 1$ one has
	\begin{equation*}
\xi = \dfrac{1}{j^2}\dfrac{\rho}{\tau},
	\end{equation*}
then, the gap condition~\eqref{differ} is not valid. Indeed, from~\eqref{bis-lamb} we deduce
	$$
\lambda_{k+j}^{(1)} - \lambda_k^{(2)} = - \left( \frac{\epsilon_{k+j}}{k+j} + \frac{\epsilon_k}{k} \right) , \quad \forall k \ge 1.
	$$
In particular ($\{ \epsilon_k \}_{k \ge 1}$ is a positive sequence), $\lambda_{k+j}^{(1)} < \lambda_{k}^{(2)}$ for any $k \ge 1$ and
	$$
\lim_{k \to \infty} \left(  \lambda_{k+j}^{(1)} - \lambda_k^{(2)} \right) = 0.
	$$
In this case, we can rearrange the sequence $\{\lambda_k^{(1)}, \lambda_k^{(2)} \}_{k\geq1}$ as follows: there exists an integer $k_0 \ge 1$ such that
	$$
\lambda_{k - 1}^{(2)} < \lambda_{k+j}^{(1)} < \lambda_{k}^{(2)}, \quad \forall k \ge k_0.
	$$
The previous inequality can be directly deduced from~\eqref{summi}.
\end{remark}

Let us now check that the sequence of eigenvalues of $L$ and $L^*$ fulfills the conditions in Lemma~\ref{1.5M}. We will do it in the next result:

\begin{proposition}\label{corh1h2}
Let us assume that the parameters $\xi$, $\rho$ and $\tau$ satisfy \eqref{H2}. Then, the sequence $\{\lambda_k^{(1)}, \lambda_k^{(2)} \}_{k\geq1}$, given by~\eqref{lamk12}, can be rearranged into an increasing sequence $\Lambda= \{ \Lambda_k\}_{k\geq1}$ that satisfies~\eqref{serconv} and $\Lambda_k\neq\Lambda_n$, for all $k,n \in \N$ with $k\neq n$. In addition, if~\eqref{H1} holds, the sequence $\{\Lambda_k\}_{k\geq1}$  also satisfies \eqref{56} and \eqref{56_2}.
\end{proposition}

\begin{proof}
As a consequence of property $(P1)$ in Proposition~\ref{propert}, we deduce that the sequence of eiegenvalues $\{\lambda_k^{(1)}, \lambda_k^{(2)} \}_{k\geq1}$ can be rearranged into a positive increasing sequence $\Lambda= \{ \Lambda_k\}_{k\geq1}$ that satisfies \eqref{serconv}. Under assumption~\eqref{H2}, we can also apply property~$(P2)$ of the same proposition and conclude that the elements of the sequence $\Lambda $ are pairwise different.

Let us now assume that, in addition, the parameters $\xi$, $\rho$ and $\tau$ also fulfill condition~\eqref{H1}. In this case, we can give an explicit rearrangement of the sequence $\{\lambda_k^{(1)}, \lambda_k^{(2)} \}_{k\geq1}$. Indeed, if $j \ge 0$ is such that the parameters satisfy~\eqref{xi}, property~$(P3)$ in Proposition~\ref{propert} provides an integer $k_0 \ge 1$ for which one has~\eqref{orden}. Thus, if $1 \le k \le 2 k_0 + j -2$, we define $\Lambda_k$ such that
	$$
\{ \Lambda_k \}_{1 \le k \le 2 k_0 + j -2} \equiv \{ \lambda_k^{(1)} \}_{1 \le k \le k_0+j -1} \cup \{ \lambda_k^{(2)} \}_{1 \le k \le k_0 -1}   \hbox{ and }  \Lambda_k < \Lambda_{k+1}, \quad \forall k: 1 \le k \le 2 k_0 + j - 3.
	$$
From the  $(2k_0 + j - 1)$-th term, we define
	\begin{equation}\label{Lambda}
	\left\{
	\begin{array}{l}
\Lambda_{2k_0 + j + 2k-1} = \lambda^{(1)}_{k_0 + j + k}, \quad \forall k \ge 0, \\
	\noalign{\smallskip}
\Lambda_{2k_0 + j + 2k} = \lambda^{(2)}_{k_0 + k}, \quad \forall k \ge 0. \\
	\end{array}
	\right.	
	\end{equation}
Clearly, $\Lambda = \{ \Lambda_k \}_{ k \ge 1}$ is an increasing sequence and $\{ \Lambda_k \}_{ k \ge 1} = \{\lambda_k^{(1)}, \lambda_k^{(2)} \}_{k\geq1}$. Furthermore, thanks to~\eqref{differ} in Proposition~\ref{propert}, the sequence $\Lambda$ also satisfies the second inequality in~\eqref{56} for every $q \ge 1$.

Our next task will be to prove the first inequality of~\eqref{56} for appropriate $q \ge 1$ and $\delta >0$. It is interesting to underline that it is enough to prove the existence of $ q \in \N$ and $\widetilde \delta >0$ such that one has
	\begin{equation}\label{56_3}
\left| \Lambda_k - \Lambda_n \right| \ge \widetilde \delta \left| k^2 - n^2 \right|, \quad \forall k,n \ge q, \quad \left| k - n\right| \ge  q.
	\end{equation}
Indeed, let us see that the first inequality in~\eqref{56} is valid for $q \ge 1$ and a new positive constant $\delta$. Observe that we can assume that $k \ge n \ge 1$. Hence, it is sufficient to prove~\eqref{56} with $k \ge n \ge 1$, with $n \le q -1 $ and $k - n \ge q$. First, it is clear that if in addition $k \le 2q$, thanks to~\eqref{H2}, we can conclude inequality~\eqref{56} for an appropriate positive constant $\delta_0$.

Let us now take $k \ge n \ge 1$, with $n \le q -1 $ and $k \ge 2q$ (and therefore, $k - n \ge q$). From~\eqref{56_3} and using $k \ge q +n \ge q+ 1$, $1 \le n \le q-1$ and $k - q \ge q$, we have %
	$$
	\left\{
	\begin{array}{l}
\displaystyle \left| \Lambda_k - \Lambda_n \right| = \Lambda_k - \Lambda_n \ge \Lambda_k - \Lambda_q \ge \widetilde \delta \left| k^2 - q^2 \right| = \widetilde \delta \left| k^2 - n^2 \right| \left[ 1- \frac{q^2 - n^2}{k^2 - n^2} \right] \\
	\noalign{\smallskip}
\displaystyle \phantom{\displaystyle \left| \Lambda_k - \Lambda_n \right|} \ge \widetilde \delta \left| k^2 - n^2 \right| \left[ 1- \frac{q^2 - n^2}{(q + 1)^2 - n^2} \right] \ge \widetilde \delta \left[ 1- \frac{q^2 - 1}{(q + 1)^2 - 1} \right]  \left| k^2 - n^2 \right| \\
	\noalign{\smallskip}
\displaystyle \phantom{\displaystyle \left| \Lambda_k - \Lambda_n \right|} = \frac{\widetilde \delta \left( 2 q +1 \right)}{(q + 1)^2 - 1} \left| k^2 - n^2 \right| .
	\end{array}
	\right.
	$$
Summarizing, assuming~\eqref{56_3}, we have deduced the first inequality in~\eqref{56} for $q \ge 1$ and 
	$$
\delta = \min\left\{ \delta_0, \,  \frac{\widetilde \delta \left( 2 q +1 \right)}{(q + 1)^2 - 1} \right\}> 0.
	$$

Let us show~\eqref{56_3} for suitable $\widetilde \delta > 0$ and $q \in \N$. To this aim, we will use the properties~\eqref{orden} and~\eqref{nueva}, in Proposition~\ref{propert}, and the expression of $\Lambda_k$ for $k \ge 2 k_0 + j -1$ (see~\eqref{Lambda}; recall that $j \ge 0$ is such that the parameters $\xi$, $\rho$ and $\tau$ satisfy~\eqref{xi}). We will work with $q \in \N$ given by
	\begin{equation}\label{q}
q \ge \max \left\{ 2k_0 + j -1, 2k_1 + 2 j +1,  6j + 3 \right\}.
	\end{equation}
Thus, if $k, n \in \N$ are such that $k ,n \ge q$ and $\left| k - n \right| \ge q$, then $\Lambda_k$ and $\Lambda_n$ are given by~\eqref{Lambda}. Depending on the expressions of $k$ and $n$, we will divide the proof of~\eqref{56_3} into three steps:

\smallskip

\textbf{1.} Assume that $k= 2k_0 + j + 2 \widetilde k - 1$ and $n = 2k_0 + j + 2 \widetilde n - 1$, for $\widetilde k, \widetilde n \ge 0 $. Since 
	$$
\left| \left( k_0 + j + \widetilde k \right) - \left( k_0 + j + \widetilde n \right) \right| = \frac 12 \left|  k -  n \right|  \ge \frac q 2 \ge k_1,
	$$
from~\eqref{Lambda} and~\eqref{nueva}, we can write
	$$
	\begin{array}{l}
\displaystyle \left| \Lambda_k - \Lambda_n \right| = \left| \lambda^{(1)}_{k_0 + j + \widetilde k} - \lambda^{(1)}_{k_0 + j + \widetilde n} \right| \ge \frac \xi 2 \left| \left( k_0 + j + \widetilde k \right)^2 - \left( k_0 + j + \widetilde n \right)^2 \right|  \\
	\noalign{\smallskip}
\displaystyle \phantom{\displaystyle \left| \Lambda_k - \Lambda_n \right| } =  \frac \xi 8 \left| \left( k + 1 + j \right)^2 - \left( n + 1 + j \right)^2 \right| = \frac \xi 8 \left| k^2 - n^2 + 2 (k-n)(1+j) \right| \ge \frac \xi 8 \left| k^2 - n^2 \right| .
	\end{array}
	$$
We obtain thus the proof of~\eqref{56_3} for $\widetilde \delta = \xi /8$ and $q$ given by~\eqref{q}.

\smallskip

\textbf{2.} The case $k= 2k_0 + j + 2 \widetilde k$ and $n = 2k_0 + j + 2 \widetilde n$, with $\widetilde k, \widetilde n \in \N$, can be treated in the same way deducing~\eqref{56_3} for $\widetilde \delta = \xi /8$ and $q$ (see~\eqref{q}).

\smallskip

\textbf{3.} Let us analyze the last case $k= 2k_0 + j + 2 \widetilde k $ and $n = 2k_0 + j + 2 \widetilde n - 1$ (with $\widetilde k, \widetilde n \in \N$), $k, n \ge q$ and $\left| k - n \right| \ge q$, with $q$ satisfying~\eqref{q}.  In this case, one has
	$$
\left| \left( k_0 + \widetilde k \right) - \left( k_0 + j + \widetilde n \right) \right| = \left| \frac 12 \left( k - n \right) - j - \frac 12 \right| \ge \frac 12 \left| k - n \right| - \left( j + \frac 12 \right) \ge \frac 12 q - \left( j + \frac 12 \right)  \ge k_1,
	$$
whence
	$$
	\begin{array}{l}
\displaystyle \left| \Lambda_k - \Lambda_n \right| = \left| \lambda^{(2)}_{k_0 + \widetilde k} - \lambda^{(1)}_{k_0 + j + \widetilde n} \right| \ge \frac \xi 2 \left| \left( k_0 + \widetilde k \right)^2 - \left( k_0 + j + \widetilde n \right)^2 \right| =  \frac \xi 8 \left| \left( k - j \right)^2 - \left( n + 1 + j \right)^2 \right| \\
	\noalign{\smallskip}
\displaystyle \phantom{\displaystyle \left| \Lambda_k - \Lambda_n \right| } = \frac \xi 8 \left| k^2 - n^2 - \left[ 2 j (k+1) + 2n(1+j) + 1 \right] \right|  .
	\end{array}
	$$
Observe that if $k \le n$, from the previous inequality, we conclude~\eqref{56_3} for for $\widetilde \delta = \xi /8$ and $q$ given by~\eqref{q}. Let us now see the case $k > n$ (and then, $ k - n = \left| k - n \right| \ge q$). The previous inequality allows us to write
	$$
	\begin{array}{l}
\displaystyle \left| \Lambda_k - \Lambda_n \right| = \frac \xi 8 \left| k^2 - n^2 - \left[ 2 j (k+1) + 2n(1+j) + 1 \right] \right| \\
	\noalign{\smallskip}
\displaystyle \phantom{\displaystyle \left| \Lambda_k - \Lambda_n \right| } \ge \frac \xi 8 \left( k^2 - n^2 \right) - \left[ 2 j (k+1) + 2n(1+j) + 1 \right] \\
	\noalign{\smallskip}
\displaystyle \phantom{\displaystyle \left| \Lambda_k - \Lambda_n \right| } = \frac \xi 8 \left( k^2 - n^2 \right) \left[ 1- \frac{2 j (k+1) + 2n(1+j) + 1}{ k^2 - n^2 } \right]  \\
	\noalign{\smallskip}
\displaystyle \phantom{\displaystyle \left| \Lambda_k - \Lambda_n \right| }\ge \frac \xi 8 \left( k^2 - n^2 \right) \left[ 1- \frac{2 j (k+1) + 2n(1+j) + 1}{ q \left(k + n \right)} \right] \\
	\noalign{\smallskip}
\displaystyle \phantom{\displaystyle \left| \Lambda_k - \Lambda_n \right| } \ge \frac \xi 8 \left( k^2 - n^2 \right) \left[ 1- \frac {2j} q - \frac {1+j} q - \frac{1}{2q} \right] \ge \frac \xi {16} \left( k^2 - n^2 \right) .
	\end{array}
	$$
Let us remark that the last inequality is valid thanks to~\eqref{q}. 

\smallskip

In conclusion, we have proved the existence of a natural number $q \ge 1$, depending on the parameters in~\eqref{PFSylinear}, such that~\eqref{56_3} holds for $\widetilde \delta = \xi/16$ and $q$ provided by formula~\eqref{q}. As a consequence, one also has~\eqref{56} for a new $\delta >0$ and the same $q$.

\smallskip

Let us now show the estimate~\eqref{56_2} for the sequence~$\Lambda = \{ \Lambda_k\}_{k\geq1} = \{\lambda_k^{(1)}, \lambda_k^{(2)} \}_{k\geq1}$. From the definition of the sequence $\Lambda$, for any $ r >0$, we can write:
	$$
\mathcal{N}(r)= \# \left\{ k :  \Lambda_k\leq r \right\}= \# \left\{k : \lambda_k^{(1)} \le r \right\} + \# \left\{k : \lambda_k^{(2)} \le r \right\} = \# \mathcal A_1 (r) +  \# \mathcal A_2 (r) = n_1 + n_2, 
	$$
where $ \mathcal A_i (r)= \left\{k : \lambda_k^{(i)} \le r \right\}$ and $n_i = \# \mathcal A_i (r)$, $i=1,2$. Our next objective will be to give appropriate bounds for $n_1$ and $n_2$.

From the definition of  $ \mathcal A_1 (r)$ and $n_1$, we deduce that $n_1$ is a natural number which is characterized by $ \lambda_{n_1}^{(1)} \le r$ and $ \lambda_{n_1 + 1}^{(1)} > r $.
Let us first work with the inequality $ \lambda_{n_1}^{(1)} \le r$. From the definition of $\lambda_k^{(1)}$ (see~\eqref{lamk12}), one gets
	$$
\xi n_1^2 + \frac{\rho + 1}{2 \tau} \le r + \sqrt{ \dfrac{\xi\rho}{\tau} n_1^2 + \left(\dfrac{\rho+1}{2\tau}\right)^2 } \le r + \sqrt{ \dfrac{\xi\rho}{\tau}} \, n_1 + \frac{\rho + 1}{2 \tau}.
	$$
The previous inequality also implies
	$$
\xi n_1^2 - \sqrt{ \dfrac{\xi\rho}{\tau}} \, n_1 - r \le 0,
	$$
and
	$$
0 \le n_1 \le \frac 1{2\xi} \left( \sqrt{ \dfrac{\xi\rho}{\tau}} + \sqrt{ \dfrac{\xi\rho}{\tau} + 4 \xi r }\right) \le \frac 1{\sqrt \xi} \left( \sqrt{ \dfrac{\rho}{\tau}} + \sqrt{ r }\right) .
	$$
From the inequality $ \lambda_{n_1 + 1}^{(1)} > r $ we also deduce,
	$$
r < \xi \left( n_1 + 1 \right)^2 + \frac{\rho + 1}{2 \tau} - \sqrt{ \dfrac{\xi\rho}{\tau} \left( n_1 + 1 \right)^2 + \left(\dfrac{\rho+1}{2\tau}\right)^2 } \le \xi \left( n_1 + 1 \right)^2, 
	$$
that is to say, $n_1 > \sqrt r/ \sqrt \xi - 1$. Summarizing, $n_1$ is a nonnegative integer such that
	\begin{equation}\label{n1}
\frac{\sqrt r}{ \sqrt \xi} - 1 < n_1 \le \frac{\sqrt r}{ \sqrt \xi} + \sqrt{ \dfrac{\rho}{\xi \tau}} , \quad \forall r \ge 0.
	\end{equation}

We can repeat the previous arguments to obtain upper and lower bounds for $n_2$. Indeed, from the definition of  $ \mathcal A_2 (r)$ and $n_2$, we get that $n_2$ is a natural number that satisfies $ \lambda_{n_2}^{(2)} \le r$ and $ \lambda_{n_2 + 1}^{(2)} > r $. The first inequality provides the estimate
	$$
r \ge \lambda_{n_2}^{(2)} \ge \xi n_2^2, \quad \hbox{i.e.,} \quad n_2 \le \frac{\sqrt r}{ \sqrt \xi}. 
	$$
On the other hand, $n_2$ is such that
	$$
0 < \lambda_{n_2 + 1}^{(2)}- r \le \xi \left( n_2 + 1 \right)^2 + \sqrt{\dfrac{\xi\rho}{\tau}} \left( n_2 + 1 \right) + \dfrac{\rho+1}{\tau} - r,
	$$
whence
	$$
	\begin{array}{l}
\displaystyle n_2 + 1 > \frac 1{2 \xi} \left[-  \sqrt{\dfrac{\xi\rho}{\tau}} + \sqrt{\dfrac{\xi\rho}{\tau} + 4 \xi \left( r - \dfrac{\rho+1}{\tau} \right)}\right] = \frac 1{2 \sqrt \xi} \left(-  \sqrt{\dfrac{\rho}{\tau}} + \sqrt{ 4 r - \dfrac{3\rho+4}{\tau} }\right) \\
	\noalign{\smallskip}
\displaystyle \phantom{\displaystyle n_2 + 1} \ge  \frac 1{2 \sqrt \xi} \left( 2 \sqrt r -  \sqrt{\dfrac{\rho}{\tau}} -\sqrt{ \dfrac{3\rho+4}{\tau} } \right) .
	\end{array}
	$$
In the last inequality we have used that $\sqrt{a-b} \ge \sqrt a - \sqrt b$ provided $a,b>0$ and $a \ge b$. In conclusion, we have proved that $n_2$ is a nonnegative integer such that
	\begin{equation}\label{n2}
\frac {\sqrt r}{\sqrt \xi} - \frac 1{2 \sqrt \xi} \left( \sqrt{\dfrac{\rho}{\tau}} + \sqrt{ \dfrac{3\rho+4}{\tau} } \right) - 1 \le n_2 \le \frac{\sqrt r}{ \sqrt \xi} , \quad \forall r \ge 0.
	\end{equation}

Recall that $\mathcal N (r) = n_1 + n_2$. Thus, from inequalities~\eqref{n1} and~\eqref{n2}, we can write
	$$
\frac {2}{\sqrt \xi} \sqrt r - \frac 1{2 \sqrt \xi} \left( \sqrt{\dfrac{\rho}{\tau}} + \sqrt{ \dfrac{3\rho+4}{\tau} } \right) - 2 \le \mathcal N (r) \le \frac {2}{\sqrt \xi} \sqrt r + \sqrt{ \dfrac{\rho}{\xi \tau}}, \quad \forall r \ge 0,
	$$
and deduce~\eqref{56_2} with
	$$
p = \frac {2}{\sqrt \xi} \quad \hbox{and} \quad \alpha = \max \left\{ \frac 1{2 \sqrt \xi} \left( \sqrt{\dfrac{\rho}{\tau}} + \sqrt{ \dfrac{3\rho+4}{\tau} } \right) + 2, \sqrt{ \dfrac{\rho}{\xi \tau}}\right\}.
	$$
This ends the proof.
\end{proof}

We will finish this section giving a result on the set of eigenfunctions of the operators $L$ and $L^*$. It reads as follows:

%
%
	\begin{proposition}\label{SS3}
Let us consider the sequences $\mathcal{F} =\{\Psi_k^{(1)}, \Psi_k^{(2)} \}_{ k\geq1}$ and $\mathcal{F^*} = \{\Phi_k^{(1)}, \Phi_k^{(2)} \}_{k\geq1}$ given in Proposition~\ref{Lspec}. Then,
		\begin{enumerate}[i)]
			\item $\mathcal{F}$ and $\mathcal{F}^*$ are biorthogonal sequences.
			\item	$\hbox{span} \left(\mathcal{F} \right)$ and $\hbox{span} \left( \mathcal{F}^* \right)$ are dense in $H^{-1}(0,\pi;\mathbb{R}^2)$, $L^2(0,\pi;\mathbb{R}^2)$ and $H_0^1(0,\pi;\mathbb{R}^2)$.
			\item $\mathcal{F}$ and $\mathcal{F}^*$ are unconditional bases\footnote{A countable sequence $\{ x_n\}_{n \ge 1}$ in a Banach space $X$ is an unconditional basis for $X$ if for every $x \in X$ there exist unique scalars $a_n(x)$ such that $x = \sum_{n \ge 1} a_n(x) x_n$, where the series converges unconditionally for each $x \in X$.} 
for $H^{-1}(0,\pi;\mathbb{R}^2)$, $L^2(0,\pi;\mathbb{R}^2)$ and $H_0^1(0,\pi;\mathbb{R}^2)$.
		\end{enumerate}
	\end{proposition}
%

\begin{proof}
From the expressions of $\Psi_k^{(j)}$ and $\Phi_k^{(j)}$ (see~\eqref{pk12} and~\eqref{pk1_22}) we can write  
	$$
\Psi_k^{(j)} (\cdot) = V_{j,k} \eta_k (\cdot),\quad \mbox{and} \quad  \Phi_k^{(j)} = V_{j,k}^* \eta_k (\cdot), \quad  j=1,2, \quad k \geq 1,
	$$
where $V_{j,k}, V_{j,k}^* \in \mathbb{R}^2$ (the function $\eta_k$ is given in~\eqref{sin}).

Item $i)$ is simple to deduce, since $\{\eta_k\}_{k\geq 1}$ is an orthogonal basis for $H^{-1}(0,\pi)$, $H_0^1(0,\pi)$ and $L^2(0,\pi)$ (in this last case, an orthonormal basis) and $\{V_{1,k}, V_{2,k}\}_{k\geq 1}$ and $\{V_{1,k}^*, V_{2,k}^*\}_{k\geq 1}$ are biorthogonal basis of $\mathbb{R}^2$. Indeed, if $M_k = \left[V_{1,k}| V_{2,k}\right]$ and $N_k = \left[V_{1,k}^*| V_{2,k}^* \right]$, then,
\begin{equation*}
			M_k^{tr} N_k = M_k  N_k^{tr} = Id,\quad \forall k\geq1.
\end{equation*}
This proves item $i)$.

For showing item $ii)$ we only need to assure that $\hbox{span} \, (\mathcal{F})$ and $\hbox{span} \, (\mathcal{F}^*)$ are dense in $H_0^{1}(0,\pi;\mathbb{R}^2)$, since $H_0^{1}(0,\pi;\mathbb{R}^2)$ is dense in $L^2(0, \pi; \R^2)$ and in $H^{-1}(0,\pi;\mathbb{R}^2)$. Let us consider $f = (f_1, f_2)^{tr}\in H^{-1}(0,\pi;\mathbb{R}^2)$ such that
	\begin{equation*}
\left\langle f, \Psi_k^{(i)} \right\rangle =0, \quad \forall k\geq1, \quad i=1,2.
	\end{equation*}
(Recall that $\langle \cdot\,, \cdot \rangle$ stands for the usual duality pairing between $H^{-1} (0,1; \R^2)$ and $H_0^1 (0,1; \R^2)$). If we denote $f_{i,k}$ ($i=1,2$) the corresponding Fourier coefficients of the distribution $f_i \in H^{-1}(0, \pi)$ with respect to the sinus basis $\{ \eta_k (\cdot) \}_{k \ge 1}$, then the previous equality can be written under the form
	$$
(f_{1,k}, f_{2,k}) \, M_k = 0, \quad \forall k\geq 1.
	$$
Using that $\det M_k \neq 0$ for any $k \ge 1$, we deduce $f_{1,k}=f_{2,k}=0$, for all $k\geq1$ and, therefore, $f=0$. This proves the density of $\mathcal{F}$ in $H_0^{1}(0,\pi;\mathbb{R}^2)$. A similar argument can be used for $\mathcal{F}^*$. This shows item $ii)$.

Let us now prove item $iii)$. As before, we will only prove that $\mathcal F $ is an unconditional basis for $H_0^1 (0, \pi; \R^2)$. This amounts to prove that, for any $f = \left( f_1, f_2 \right)^{tr}\in H_0^1 (0, \pi; \R^2) $, the series 
	$$
S(f) := \sum_{k\geq1} \left( \left\langle \Phi_k^{(1)}, f \right\rangle \Psi_k^{(1)} + \left\langle \Phi_k^{(2)}, f \right\rangle \Psi_k^{(2)} \right)
	$$
is unconditionally convergent in $H_0^1 (0, \pi; \R^2)$. From the definition of the functions $\Psi_k^{(i)}$ and $\Phi_k^{(i)}$ (see~\eqref{pk12} and~\eqref{pk1_22}), it is easy to see that
	\begin{equation*}
	S(f)  = \sum_{k \ge 1} \left(
	\begin{array}{l}
f_{1,k}\\
f_{2,k}
	\end{array}
	\right) \eta_k,
	\end{equation*}
where $f_{i,k}$ is the Fourier coefficient of the function $f_i \in H_0^1(0, \pi)$ ($i=1,2$). Accordingly, this series converges unconditionally in $H_0^1(0,\pi;\mathbb{R}^2)$ (recall that $\{\eta_k\}_{k\geq1}$ is an orthogonal basis for $H_0^1(0,\pi)$ and $f_1,f_2\in H_0^1(0,\pi)$). This concludes the proof of the result.
\end{proof}
%
%
\begin{section}{Approximate and null controllability of the linear system~\eqref{PFSylinear}}\label{s4}

\setcounter{equation}{0}

We will devote this section to proving the approximate and null controllability at time $T>0$ of system~\eqref{PFSylinear}. To this aim, we will use in a fundamental way the properties of the spectrum of the operator $L$ (see~\eqref{LL*}) established in Propositions~\ref{Lspec},~\ref{propert} and~\ref{corh1h2}. Firstly, we will show the result on approximate controllability of the linear system (Theorem~\ref{CAproximada_}) and then the null controllability at time $T$ of the same system (Theorem~\ref{CNull_}).

%
%
\begin{subsection}{Approximate controllability: Proof of Theorem~\ref{CAproximada_}}\label{s4.1}
Let us fix $T>0$ and consider system~\eqref{PFSylinear} with $\xi,  \rho, \tau >0$ given. Let us first assume that system~\eqref{PFSylinear} is approximate controllable at time $T$. In this case, condition~\eqref{H2} holds. Indeed, otherwise, thanks to property~$(P2)$ of Proposition~\ref{propert}, the spectrum of the operator~$L$ is not simple, i.e., there exist $k,\ell \geq 1$  such that $\lambda_k^{(2)} = \lambda_\ell^{(1)}= \lambda_0$. Thus, if we take $a,b \in \R$, it is easy to see that the function
	$$
\varphi (x,t) =  \left(a\, \Phi_\ell^{(1)}(x)  + b\, \Phi_k^{(2)}(x)\right)e^{-\lambda_0 (T-t)}, \quad \forall (x,t) \in Q_T,
	$$
is the solution of the adjoint system~\eqref{adjoint} associated to the initial condition
	$$
\varphi_0 = a\Phi_\ell^{(1)} + b \Phi_k^{(2)} . 
	$$
This function satisfies (see~\eqref{ADB} and~\eqref{pk1_22}) 
	\begin{equation*}
B^*D^*\varphi_x(0,t)	= \xi \sqrt{\dfrac{2}{\pi\tau}} \left ( a \dfrac{\ell}{\sqrt{r_\ell}} -	b \dfrac{k}{\sqrt{r_k}} \right) e^{-\lambda_0 (T-t)}, \quad \forall t \in (0,T).
	\end{equation*}
Choosing
	$$
a = \frac k{\sqrt{r_k}} \quad \mbox{and} \quad b =\frac \ell{\sqrt{r_\ell}},
	$$
we have that $B^*D^*\varphi_x(0,\cdot)=0$ but $\varphi_0\neq 0$, contradicting the unique continuation property stated in the first point of Theorem~\ref{equivcontrol}. In conclusion, system~\eqref{PFSylinear} is not approximately controllable at time $T>0$. 

Let us now suppose that condition~\eqref{H2} holds and prove the unique continuation property for system~\eqref{adjoint}. Again, from the first point of Theorem~\ref{equivcontrol} we infer the approximate controllability property of system~\eqref{PFSylinear}.

Let us consider $\varphi_0 \in H_0^1(0,\pi)$ and assume that the corresponding solution $\varphi$ to the adjoint problem~\eqref{adjoint} satisfies
	$$
B^*D^*\varphi_x(0,t) = 0, \quad \forall t \in (0,T).
	$$
Observe that, thanks to Proposition~\ref{adjexist}
	$$
\varphi\in C^0([0,T];H_0^1(0,\pi;\mathbb{R}^2))\cap L^2(0,T;H^2(0,\pi;\mathbb{R}^2) \cap H_0^1(0,\pi;\mathbb{R}^2)),
	$$ 
and then, $B^*D^*\varphi_x(0,\cdot ) \in L^2 (0, T) $.

From Proposition~\ref{SS3} we deduce that $\mathcal{F}$ and $\mathcal{F}^*$ are biorthogonal bases for $H^{-1}(0,\pi;\mathbb{R}^2)$ and $H_0^1(0,\pi;\mathbb{R}^2)$. In particular, $\varphi_0 \in H_0^1(0,\pi;\mathbb{R}^2)$ can be written as
	$
\varphi_0 = \displaystyle \sum_{k\geq 1} \left( a_k \Phi_k^{(1)} + b_k \Phi_k^{(2)} \right) ,
	$
where
	$$ 
a_k = \left\langle \Psi_k^{(1)} , \varphi_0 \right\rangle, \quad b_k = \left\langle \Psi_k^{(2)} , \varphi_0 \right\rangle , \quad \forall k \ge 1.
	$$
Using Proposition~\ref{Lspec}, the corresponding solution $\varphi$ of system~\eqref{adjoint} associated to $\varphi_0$ is given by
	$$
\varphi(\cdot ,t) = \displaystyle \sum_{k\geq 1} \left( a_k \Phi_k^{(1)} e^{-\lambda_k^{(1)}(T-t)} + b_k \Phi_k^{(2)} e^{-\lambda_k^{(2)}(T-t)} \right) , \quad \forall t \in (0,T),
	$$
where $\lambda_k^{(i)}$, $\Psi_k^{(i)}$ and $\Phi_k^{(i)}$ ($k \ge 1$, $i=1,2$) are given in Proposition~\ref{Lspec}. Therefore,
	\begin{equation*}
0 = B^*D^*\varphi_x(0,t) = \displaystyle \sum_{k\geq 1} \sqrt{\frac 2\pi} \frac{k \xi}{ \sqrt{\tau r_k}} \left( a_k e^{-\lambda_k^{(1)}(T-t)} - b_k e^{-\lambda_k^{(2)}(T-t)} \right), \quad \forall t \in (0,T).
	\end{equation*}
From Proposition~\ref{corh1h2}, we can apply Lemma \ref{bioex} in order to deduce the existence of a biorthogonal family $\{ {q}^{(1)}_k,  {q}^{(2)}_k\}_{k\geq1}$ to $\{e^{-\lambda_k^{(1)}t)}, e^{-\lambda_k^{(2)}t)}\}_{k\geq1}$ in $L^2(0,T)$. Then, the previous identity, in particular, implies 
	\begin{equation*}
	\left\{
	\begin{array}{l}
\displaystyle\sqrt{\frac 2\pi} \frac{k \xi}{ \sqrt{\tau r_k}} \, a_k = \displaystyle\int_0^T B^*D^*\varphi_x(0,t) \, {q}^{(1)}_k (t) \, dt = 0, \quad \forall k \ge 1, \\
	\noalign{\medskip}
\displaystyle \sqrt{\frac 2\pi} \frac{k \xi}{ \sqrt{\tau r_k}} \, b_k = - \displaystyle\int_0^T B^*D^*\varphi_x(0,t) \, {q}^{(2)}_k (t) \, dt = 0, \quad \forall k \ge 1,
	\end{array}
	\right.
	\end{equation*}
and $a_k = b_k =0$ for any $k \ge 1$. In conclusion, $\varphi_0 = 0$ and we have proved the unique continuation property for the solutions of system~\eqref{adjoint}. This ends the proof of Theorem~\ref{CAproximada_}.
\end{subsection}
%

%
%

\begin{subsection}{Null controllability: Proof of Theorem~\ref{CNull_}}\label{s4.2}


Let us now prove the null controllability result stated in Theorem~\ref{CNull_}. To this aim, we consider $\xi$, $\rho$ and $\tau$ three positive real numbers satisfying assumptions~\eqref{H2} and~\eqref{H1}. We will obtain the proof writing the controllability problem for system~\eqref{PFSylinear} as a moment problem (see~\cite{FaRu}).

Let us take $y_0 = (\theta_0, \phi_0) \in H^{-1} (0, \pi; \R^2)$. As a consequence of Proposition~\ref{premomprob}, we have that the control $v \in L^2(0,T)$ is such that the solution $y=(\theta, \phi) \in C^0([0,T]; H^{-1} (0, \pi; \R^2))$ of system~\eqref{PFSylinear} satisfies $y (\cdot, T) =0$ if and only if $v \in L^2 (0,T)$ fulfills
	\begin{equation*}
\displaystyle \int_0^T B^*D^*\varphi_x(0,t) v(t) \, dt = - \langle y_0, \varphi(\cdot,0) \rangle , \quad \forall \varphi_0 \in H_0^{1}(0, \pi; \R^2),
	\end{equation*}
where $\varphi \in C^0([0,T]; H_0^{1} (0, \pi; \R^2))$ is the solution of the adjoint system~\eqref{adjoint} associated to $\varphi_0$. Using again Proposition~\ref{SS3}, we deduce that $\mathcal{F}^*$ is a basis of $H_0^1(0, \pi ; \R^2)$. In particular, we also deduce that the previous equality is equivalent to
	\begin{equation*}
\displaystyle \int_0^T B^*D^*\varphi^{(j)}_{k,x}(0,t) v(t) \, dt = - \left\langle y_0, \varphi^{(j)}_{k}(\cdot,0) \right\rangle , \quad \forall k \ge 1, \quad j=1,2, 
	\end{equation*}
where $\varphi^{(j)}_{k} (\cdot, t) = e^{-\lambda_k^{(j)}(T-t)}\Phi_k^{(j)}$ is the solution of system~\eqref{adjoint} corresponding to $\varphi_0 = \Phi_k^{(j)}$. Taking into account the expressions of $B$, $D$ and $\Phi_k^{(j)}$ (see~\eqref{ADB} and~\eqref{pk1_22}), we infer that $v \in L^2 (0,T)$ is a null control for system~\eqref{PFSylinear} associated to $y_0$ if and only if  

	\begin{equation*}
\displaystyle (-1)^{j+1} \sqrt{\frac{2}{\pi}} \frac{k \xi }{\sqrt{\tau r_k}} \int_0^T e^{-\lambda_k^{(j)}(T-t) } v(T-t) \,  dt = e^{-\lambda_k^{(j)}T}\left\langle y_0,  \Phi_k^{(j)} \right \rangle, \quad \forall k \ge 1, \quad j=1,2.
	\end{equation*}

Summarizing, we have transformed the null-controllability problem at time $T > 0$ for system~\eqref{PFSylinear} into the following moment problem: given $y_0 = (\theta_0, \phi_0) \in H^{-1} (0, \pi; \R^2)$, find $v \in L^2 (0,T)$ such that the function $u (t) := v(T-t) \in L^2 (0,T)$ satisfies
	\begin{equation}\label{tildev}
\displaystyle \int_0^T e^{-\lambda_k^{(j)} t} u(t) \, dt = c_{kj}, \quad \forall k\geq1,\quad j=1,2,
	\end{equation}
where $c_{kj} = c_{kj}(y_0)$ is given by
	\begin{equation}\label{ckj_}
c_{kj} = (-1)^{j+1} \sqrt{\frac{\pi}{2}} \, \frac{\sqrt{\tau r_k}}{k \xi } \, e^{-\lambda_k^{(j)}T}\left\langle y_0,  \Phi_k^{(j)} \right \rangle , \quad \forall k\geq1,\quad j=1,2 .
	\end{equation}

Our next task will be to solve problem~\eqref{tildev}. The assumptions~\eqref{H2} and~\eqref{H1}, Proposition~\ref{corh1h2} and Lemma~\ref{1.5M} guarantee the existence of $\widetilde T_0 >0$ such that for any $T \in (0, \widetilde T_0)$ there exists a biorthogonal family $\{ {q}^{(1)}_k, {q}^{(2)}_k \}_{k\geq1}$ to $\{e^{-\lambda_k^{(1)}t }, e^{-\lambda_k^{(2)}t }\}_{k\geq1}$ in $L^2(0,T)$ which satisfies
	\begin{equation}\label{bound}
\left\|q_k^{(j)}\right\|_{L^2(0,T)} \leq C e^{C\sqrt{\lambda_k^{(j)}} + \frac{C}{T}},\quad \forall k \geq 1, \quad j=1,2,
	\end{equation}
for a positive constant $C$ independent of $T$.

Let us first prove the result when $T \in (0,\widetilde T_0)$. Then, a formal solution to the moment problem~\eqref{tildev} is given by
	\begin{equation}\label{control}
		\begin{array}{lllll}
u (t) := v(T-t) = \sum_{k\geq 1} \left( c_{k1} {q}_k^{(1)} + c_{k2} {q}_k^{(2)} \right). 
		\end{array}
	\end{equation}

Let us now prove that $ u \in L^2(0,T)$ and, consequently, that $v \in L^2(0,T)$. From the expressions of $r_k$, $\lambda_k^{(j)}$ and $\Phi_k^{(j)}$ (see~\eqref{rkquant}, \eqref{lamk12} and~\eqref{pk1_22}) we can easily deduce the existence of constants $C,C_1,C_2>0$ such that
	\begin{equation*}
C_1k\leq r_k \leq C_2k,\quad C_1k^{2} \leq \left |\lambda_k^{(j)} \right| \leq C_2 k^2,\quad \left\|\Phi_k^{(j)} \right\|_{H_0^1} \leq C k^{3/2}, \quad \forall k\geq1, \quad j=1,2 ,
	\end{equation*}
and, from~\eqref{ckj_},
	\begin{equation*}
\left| c_{kj} \right| \leq  \frac{C}{\sqrt k}  \, e^{-\lambda_k^{(j)}T} \left\| y_0 \right\|_{H^{-1}} \left\|\Phi_k^{(j)} \right\|_{H_0^1} \leq C \,k \, e^{-\lambda_k^{(j)}T} \left\| y_0 \right\|_{H^{-1}} ,  \quad \forall k\geq1, \quad j=1,2.
	\end{equation*}

Coming back to the expression of the null control $v$ (see~\eqref{control}) and taking into account~\eqref{bound} and the previous inequality, we get
	\begin{equation}\label{cy0depcont}
	\begin{array}{l}
\displaystyle \|v\|_{L^2(0,T)} \leq C \, e^{\frac{C}{T}} \|y_0\|_{H^{-1}} \displaystyle \sum_{k\geq 1} \left( e^{C\sqrt{\lambda_k^{(1)}}}e^{-\lambda_k^{(1)}T} + e^{C\sqrt{\lambda_k^{(2)}}}e^{-\lambda_k^{(2)}T} \right)  \\
	\noalign{\smallskip}
\displaystyle \phantom{\displaystyle \|v\|_{L^2(0,T)}} \leq C \, e^{\frac{C}{T}}  \|y_0\|_{H^{-1}}\sum_{k\geq 1} \left( e^{\frac{C^2}{2T} + \frac{1}{2}\lambda_k^{(1)}T} e^{-\lambda_k^{(1)}T} + e^{\frac{C^2}{2T} + \frac{1}{2}\lambda_k^{(2)}T} e^{-\lambda_k^{(2)}T}\right) \\
	\noalign{\smallskip}
\displaystyle \phantom{\displaystyle \|v\|_{L^2(0,T)}} \leq C \, e^{\frac{C}{T}}  \|y_0\|_{H^{-1}} \sum_{k\geq 1}e^{-CTk^2} \leq C \, e^{\frac{C}{T}}  \|y_0\|_{H^{-1}} \int_0^{ \infty}e^{-CT s^2} \, ds = \frac{C}{2}\sqrt{\frac{\pi}{CT}} \, e^{\frac{C}{T}}  \|y_0\|_{H^{-1}} \\
	\noalign{\smallskip}
\displaystyle \phantom{\displaystyle \|v\|_{L^2(0,T)}} \leq C_0 \, e^{\frac{M}{T}}  \|y_0\|_{H^{-1}},
	\end{array}
	\end{equation}
for positive constants $C_0$ and $M$ independent of $T$. This inequality shows that $v \in L^2(0,T)$ and proves the first part of Theorem~\ref{CNull_}.

The second part is a direct consequence of the expression of the null control $v$ (see~\eqref{control}) and~\eqref{cy0depcont}. Indeed, if we define the operator $\mathcal{C}_T^{(0)}: H^{-1}(0, \pi ; \R^2) \rightarrow L^2(0,T)$ by
	\begin{equation*}
\mathcal{C}_T^{(0)} \left(y_0 \right) := \displaystyle \sum_{k\geq1}\left( c_{k1}(y_0) {q}_k^{(1)}(T- \cdot ) + c_{k2}(y_0) {q}_k^{(2)}(T- \cdot ) \right), \quad \forall y_0 \in H^{-1}(0, \pi ; \R^2),
	\end{equation*}
with $c_{kj}=c_{kj}(y_0)$ given by \eqref{ckj_}, it is not difficult to see that $\mathcal{C}_T^{(0)}$ is a linear operator which satisfies~\eqref{CT0} for a positive constants $C_0$ and $M$.  This ends the proof of Theorem \ref{CNull_} when $T \in (0, \widetilde T_0)$.

Let us now assume that $T \ge \widetilde T_0$. We will obtain the proof as a consequence of the previous case. Indeed, if $T \ge \widetilde T_0$ we can construct a null control at time $T$ for system~\eqref{PFSylinear} associated to $y_0 \in H^{-1} (0, \pi; \R^2)$ as
	$$
v (t) =  \mathcal{C}_{ T}^{(0)} \left( y_0 \right) (t) := \left\{
	\begin{array}{ll}
\displaystyle \mathcal{C}_{\widetilde T_0/2}^{(0)} \left( y_0 \right) (t) &\hbox{if } t \in \left[ 0 , \frac{\widetilde T_0}{2} \right], \\
	\noalign{\smallskip}
0 &\hbox{if } t \in \left[ \frac{\widetilde T_0}2, T\right]. 
	\end{array}
\right.
	$$
Clearly $\mathcal{C}_{ T}^{(0)} \in \mathcal{L} \left(H^{-1}(0,\pi;\mathbb{R}^2), L^2(0,T) \right)$
	$$
\left\| \mathcal{C}_{ T}^{(0)} \left( y_0 \right) \right\|_{L^2(0,T)} \le C_0 \, e^{2M/\widetilde T_0}  \|y_0\|_{H^{-1}} = C_1  \|y_0\|_{H^{-1}}
	$$
with $C_1$ a new positive constant independent of $T$. So, we can conclude~\eqref{CT0} for a new positive constant $C_0$ (only depending on the parameters in system~\eqref{PFSylinear}) and the same constant $M>0$ as before. This finishes the proof of Theorem~\ref{CNull_}.

\end{subsection}
\end{section}

%
%

\begin{section}{Boundary controllability of the phase-field system}\label{s5}

\setcounter{equation}{0}

In this section we will prove the exact controllability at time $T>0$ of the phase-field system~\eqref{PFSy} to the constant trajectory $(0,c)$, with $c= \pm 1$. To this end, we will perform a fixed-point strategy which will use in a fundamental way a null controllability result for the non-homogeneous linear system~\eqref{PFSylinear} ($f \in L^2(0,\pi; \R^2)$ is a given function in an appropriate weighted-Lebesgue space; see~\eqref{FVY}). 


\subsection{Null controllability of the non-homogeneous system~\eqref{PFSylinear}}\label{s5.1}

As said before, our next objective will be to show a null controllability result for non-homogeneous system~\eqref{PFSylinear} when $y_0=(\theta_0,\phi_0)\in H^{-1}(0,\pi;\mathbb{R}^2)$ and $f$ is a given function satisfying appropriate assumptions. To this end, we will follow some ideas from~\cite{tucsnack}.

Let us consider $\xi$, $\rho$ and $\tau$ three positive real numbers satisfying hypotheses~\eqref{H2} and~\eqref{H1}. The starting point is Theorem~\ref{CNull_} and Remark~\ref{r1.3}. As a consequence, we obtain an estimate for the cost of the null control of system~\eqref{PFSylinear}. With the notations of Remark~\ref{r1.3}, one has
	$$
\mathcal{K} (T) \leq C_0 \, e^{\tfrac{M}{T}},\quad \forall T > 0,
	$$
with $C_0$ and $M$ two positive constants only depending on $\xi$, $\rho$ and $\tau$.

%
%
%

In order to provide a null controllability result for the non-homogeneous problem \eqref{vectPFSy} at time $T>0$, we will introduce the functions $\gamma(t) := e^{\tfrac{M}{t}}$, $\forall t > 0 $, and, for $t \in [0,T]$,
	\begin{equation}\label{rhoF0}
\rho_{\mathcal{F}}(t) := e^{-\tfrac{b^2(a+1)M}{(b-1)(T-t)}}, \quad \rho_{0}(t) := e^{-\tfrac{aM}{(b-1)(T-t)}}, \quad \forall  t \in \left[ T\left( 1-\dfrac{1}{b^2} \right), T \right],
	\end{equation}
extended to $\left[ 0, T( 1- {1}/{b^2}) \right]$ in a constant way. Here $a,b>1$ are constants that will be chosen later. Observe that $\gamma$, $\rho_{\mathcal{F}}$ and $\rho_0$ are continuous and non increasing functions in $[0,T]$ and $\rho_{\mathcal{F}} (T) = \rho_0 (T) = 0$. 

With the previous functions, we also introduce the weighted normed spaces
	\begin{equation}\label{FVY}
	\begin{array}{c}
\displaystyle \mathcal{F} := \left\{ f \in L^2(Q_T;\mathbb{R}^2) : \frac{f}{\rho_\mathcal{F}} \in L^2(Q_T;\mathbb{R}^2) \right\}, \quad \mathcal{V} := \left\{ v \in L^2(0,T) : \frac{v}{\rho_0} \in L^2(0,T) \right\}, \\
	\noalign{\smallskip}
\displaystyle \mathcal{Y}_0 := \left\{ y \in  L^2(Q_T;\mathbb{R}^2) : \dfrac{y}{\rho_0} \in L^2(Q_T;\mathbb{R}^2) \right\},  \\
	\noalign{\smallskip}
\displaystyle \mathcal{Y} := \left\{ y \in  L^2(Q_T;\mathbb{R}^2) : \dfrac{y}{\rho_0} \in L^2(Q_T)\times C^0(\overline{Q}_T) \right\}.
	\end{array}
	\end{equation}
It is clear that $\mathcal{F}$, $\mathcal{V}$ and $\mathcal{Y}_0$ are Hilbert spaces. For instance, the inner product in $\mathcal F$ is given by
	$$
( f_1, f_2)_{\mathcal F} := \intdoble_{Q_T} \rho_{\mathcal F}^{-2} (t) f_1 (x,t) \cdot f_2 (x,t) \, dx \, dt, \quad \forall f_1,f_2 \in \mathcal F.
	$$
A similar definition can be made for $( \cdot , \cdot)_{\mathcal V}$ and $( \cdot , \cdot)_{\mathcal Y_0}$. On the other hand, $\mathcal{Y}$ is a Banach space with the norm
	$$
\| y \|_{\mathcal{Y}} := \left( \|y_1/ \rho_0\|^2_{L^2(Q_T)} + \|y_2/ \rho_0\|^2_{C^0( \overline Q_T)} \right)^{1/2}, \quad \forall y= (y_1, y_2) \in \mathcal{Y}.
	$$

With the previous notation, one has:

\begin{theorem}\label{boundsM5.1}
Let us consider $\xi$, $\rho$ and $\tau$ three positive real numbers satisfying~\eqref{H2} and~\eqref{H1}. Then, for every $T > 0$, there exist two bounded linear operators
	$$	
\mathcal{C}_T^{(1)}: H^{-1}(0,\pi;\mathbb{R}^2)\times\mathcal{F} \rightarrow \mathcal V
\quad \mbox{and} \quad
E_T^{(0)}: H^{-1}(0,\pi;\mathbb{R}^2)\times\mathcal{F} \rightarrow \mathcal Y_0
	$$
such that
\begin{enumerate}[(i)]
\item $\left\|\mathcal{C}_T^{(1)} \right\|_{\mathcal{L}(H^{-1}(0,\pi; \R^2)\times \mathcal{F},\mathcal V)}\leq C \, e^{C \left( T+ \frac 1T \right)} $ and $\left\|E_T^{(0)} \right\|_{\mathcal{L}(H^{-1}(0,\pi;\R^2) \times\mathcal{F},\mathcal Y_0)}\leq C \, e^{C \left( T+ \frac 1T \right)}$ for a positive constant $C$ independent of $T$.
\item $E_T^{(1)}:=E_T^{(0)}\left|_{H^{-1}(0,\pi)\times H_0^1(0,\pi)\times\mathcal{F}}\right. \in \mathcal{L}(H^{-1}(0,\pi)\times H_0^1(0,\pi)\times\mathcal{F},\mathcal Y)$ and, for a new constant $C>0$ independent of $T$, one has $\left\| E_T^{(1)}\right\|_{\mathcal{L}(H^{-1}(0,\pi)\times H_0^1(0,\pi)\times\mathcal{F},\mathcal Y)}\leq C \, e^{C \left( T+ \frac 1T \right)}$.
\item	For any $(y_0,f) \in H^{-1}(0,\pi;\mathbb{R}^2)\times\mathcal{F}$ (resp., $(y_0,f)\in H^{-1}(0,\pi)\times H_0^1(0,\pi)\times\mathcal{F}$), $y=	E_T^{(0)}(y_0,f)\in \mathcal Y_0$ (resp. $y = E_T^{(1)}(y_0,f)\in \mathcal Y$) is the solution of \eqref{vectPFSy} associated to $(y_0,f)$ and $v=\mathcal{C}_T^{(1)}(y_0,f)$.
\end{enumerate}
\end{theorem}


\begin{remark}
Before giving the proof of this result, let us underline that Proposition~\ref{boundsM5.1} provides a null controllability result for the non-homogeneous system~\eqref{vectPFSy} when $y_0 \in H^{-1}(0, \pi; \R^2)$ and $f \in \mathcal F$. Indeed, since $\rho_0$ is a continuous function on $[0,T]$ satisfying $\rho_0 (T)= 0$, it is clear that
	$$
y = E_T^{(0)}(y_0,f)\in \mathcal Y_0 \cap C^0([0,T];H^{-1}(0,\pi;\mathbb{R}^2)),
	$$  
solves~\eqref{vectPFSy} and satisfies $y(\cdot, T) = 0$ in $H^{-1}(0, \pi; \R^2)$.

\end{remark}


\begin{proof}[Proof of Theorem~\ref{boundsM5.1}]

In order to prove Theorem~\ref{boundsM5.1}, we adapt to system~\eqref{vectPFSy} a general technique developed in~\cite{tucsnack}  that permits to prove a null controllability result for a  non-homogenous linear problem from the corresponding controllability result for the homogenous problem. 

Let us consider $a,b>1$ and $T > 0$. With the previous definitions and notations, we define the sequence 
	\begin{equation*}
T_k = T-\dfrac{T}{b^k},\quad \forall k\geq0.
	\end{equation*}
From the definition of the functions $\rho_0$ and $\rho_{\mathcal F}$ (see~\eqref{rhoF0}) and the expression of $T_k$, one has
	\begin{equation}\label{est0}
\rho_0(T_{k+2})=\rho_{\mathcal F}(T_k) e^{\frac{M}{T_{k+2} - T_{k+1}}},\quad \forall k\geq 0.
	\end{equation}
This formula will be used in what follows.

%

Let us take $y_0 = (\theta_0, \phi_0) \in H^{-1}(0, \pi; \R^2)$ (resp., $y_0 \in H^{-1}(0, \pi) \times H_0^{1}(0, \pi)$) and $f \in \mathcal F$. Thus, we introduce the sequence $\{a_k\}_{k\geq 0} \subset H^{-1}(0,\pi;\mathbb{R}^2)$ (resp.~$\{a_k\}_{k\geq 0} \subset H^{-1}(0,\pi)\times H_0^1(0,\pi)$ if $y_0\in H^{-1}(0,\pi)\times H_0^1(0,\pi)$) defined by
	\begin{equation*}
a_0=y_0,\quad a_{k+1}=\tilde{y}_k (T_{k+1}^-),\quad \forall k\geq0,
	\end{equation*}
where $\tilde{y}_k$ is the solution to the linear system
	\begin{equation}\label{syst_tilde}
	\left\{
	\begin{array}{ll}
\tilde{y}_t - D \tilde{y}_{xx} + A \tilde{y} = f & \mbox{in } (0,\pi) \times  (T_k,T_{k+1}) := Q_k,\\
	\noalign{\smallskip}
\tilde{y}(0, \cdot) = \tilde{y}(\pi, \cdot) = 0 & \mbox{on } (T_k,T_{k+1})\\
	\noalign{\smallskip}
\tilde{y}(\cdot, T_k^+) = 0	& \mbox{in } (0,\pi),\\
	\end{array}
	\right. 
	\end{equation}
(the matrices $D$ and $A$ are given in~\eqref{ADB}). From Proposition~\ref{adjexist}, it is clear that this system admits a unique solution 
	$$
\tilde y_k \in L^2(T_k,T_{k+1} ;H^2(0,\pi) \cap H_0^1(0,\pi))\cap C^0([T_k,T_{k+1}];H_0^1(0,\pi; \R^2))
	$$
which satisfies~\eqref{acot1}. In particular, $\tilde y_k \in C^0(\overline Q_k; \R^2)$ and $a_{k+1} \in H_0^1 (0, \pi; \R^2)$, for any $k \ge 0$, and
	\begin{equation}\label{ak+1<f}
\| \tilde y_k \|_{ C^0(\overline Q_k; \R^2)} + \| a_{k+1} \|_{H_0^1} \leq e^{CT}  \|f\|_{L^2( Q_k ;\mathbb{R}^2)},\quad \forall k\geq 0,
	\end{equation}
where $C$ is a positive constant only depending on the coefficients of $D$ and $A$. 

For $k\geq 0$, we also consider the controlled autonomous problem
	\begin{equation}\label{syst_hat}
	\left\{
	\begin{array}{ll}
\hat{y}_t - D \hat{y}_{xx} + A \hat{y} = 0 & \mbox{in } Q_k,\\
	\noalign{\smallskip}
\hat{y}(0, \cdot) = Bv_k, \quad \hat{y}(\pi, \cdot) = 0	& \mbox{on } (T_k,T_{k+1})\\
	\noalign{\smallskip}
\hat{y}(\cdot, T_k^+) = a_k,\quad \hat{y}(\cdot, T_{k+1}^-) = 0 & \mbox{in } (0,\pi),\\
	\end{array}
	\right. 
	\end{equation}
where the control $v_k $ is given by $v_k = \mathcal{C}_{T_{k+1}-T_k}^{(0)}(a_k) \in L^2(T_k,T_{k+1})$ (the linear operator $ \mathcal{C}_{T_{k+1}-T_k}^{(0)}$ is given in Theorem~\ref{CNull_}). Thanks to Proposition~\ref{linearexist}, the solution $\hat y_k$ of the previous system satisfies 
	$$ 
	\left\{
	\begin{array}{l}
\hat y_0 \in L^2(Q_0;\mathbb{R}^2) \quad \hbox{(resp., }\hat y_0 \in L^2(Q_0;\mathbb{R}^2)\cap C^0([0 ,T_{1}];H^{-1}(0,\pi) \times H_0^{1}(0,\pi) \hbox{),} \\
\noalign{\smallskip}
\hat y_k \in L^2(Q_k;\mathbb{R}^2)\cap C^0([T_k,T_{k+1}];H^{-1}(0,\pi) \times H_0^{1}(0,\pi)), \quad \forall k \ge 1
	\end{array}
	\right.
	$$
and, from~\eqref{acot2} (resp.,~\eqref{acot4}), \eqref{ak+1<f} and Theorem~\ref{CNull_},
	\begin{equation*}
	\left\{
	\begin{array}{l}
\| \hat y_0 \|_{L^2(Q_0; \R^2)} \le e^{CT_1} \left( \| y_0 \|_{H^{-1}} + \| v_0 \|_{L^2(0,T_1)}\right) \le C_0 \, e^{CT} e^{\frac M{T_1}} \| y_0 \|_{H^{-1}} \\
	\noalign{\smallskip}
 \hbox{(resp., } \| \hat y_0 \|_{L^2(Q_0) \times C^0 (\overline Q_0)} \le C_0 \, e^{CT} e^{\frac M{T_1}} \| y_0 \|_{H^{-1} \times H_0^1} \hbox{)},
 	\end{array}
	\right.
	\end{equation*}
and, for any $k \ge 1$,
	$$
\| \hat y_k \|_{L^2(Q_k) \times C^0 (\overline Q_k)} \le e^{CT} \left( \| a_k \|_ {H^{-1} \times H_0^1} + \| v_k \|_{L^2(T_k,T_{k+1})}\right) \le C_0 e^{CT} e^{\frac{M}{T_{k+1}- T_k}} \| f \|_{L^2(Q_k; \R^2)}.
	$$

If we set $Y_k := \tilde y_k + \hat y_k$ in $Q_k = (0,\pi) \times  (T_k,T_{k+1})$, then 
	$$ 
	\left\{
	\begin{array}{l}
Y_0 \in L^2(Q_0;\mathbb{R}^2) \quad \hbox{(resp., }Y_0 \in L^2(Q_0;\mathbb{R}^2)\cap C^0([0 ,T_{1}];H^{-1}(0,\pi) \times H_0^{1}(0,\pi) \hbox{),} \\
\noalign{\smallskip}
Y_k \in L^2(Q_k;\mathbb{R}^2)\cap C^0([T_k,T_{k+1}];H^{-1}(0,\pi) \times H_0^{1}(0,\pi)), \quad \forall k \ge 1
	\end{array}
	\right.
	$$
and
	\begin{equation}\label{yhat}
	\left\{
	\begin{array}{l}
\| Y_0 \|_{L^2(Q_0; \R^2)} \le C \, e^{CT} e^{\frac M{T_1}} \left( \| y_0 \|_{H^{-1}} +  \|f\|_{L^2( Q_0 ;\mathbb{R}^2)}\right) \\
	\noalign{\smallskip}
 \hbox{(resp., } \| Y_0 \|_{L^2(Q_0) \times C^0 (\overline Q_0)} \le C \, e^{CT} e^{\frac M{T_1}} \left(\| y_0 \|_{H^{-1} \times H_0^1} +  \|f\|_{L^2( Q_0 ;\mathbb{R}^2)}\right) \hbox{)}, \\
	\noalign{\smallskip}
\| Y_k \|_{L^2(Q_k) \times C^0 (\overline Q_k)} \le C \, e^{CT} e^{\frac{M}{T_{k+1}- T_k}} \| f \|_{L^2(Q_k; \R^2)}, \quad \forall k \ge 1.
	\end{array}
	\right.
	\end{equation}

Let us divide the proof into two cases: the case $k=0$ and the case $k \ge 1$.

\smallskip

\textbf{Case $k=0$.} First, from Theorem~\ref{CNull_}, we can use that $bT_1=T(b-1)$ to obtain (recall that $v_0 = \mathcal{C}_{T_{1}}^{(0)}(y_0)$)
	\begin{equation*}
\| v_0 \|_{L^2(0,T_1)} \leq C_0 \, e^{\frac{M}{T_1}}\|y_0\|_{H^{-1}} = C_0 \, e^{\frac{Mb (a+1)}{(b-1)T}} \rho_0 (T_1) \|y_0\|_{H^{-1}}.
	\end{equation*}
Using now that $\rho_0$ is a positive continuous non-increasing function, from the previous estimate, we deduce the existence of a positive constant $C$ such that
	\begin{equation}\label{v0}
\left\| \frac{v_0}{\rho_0} \right\|_{L^2(0,T_1)} \leq C \, e^{\frac CT} \|y_0\|_{H^{-1}}.
	\end{equation}

On the other hand, from~\eqref{yhat},
	$$
	\left\{
	\begin{array}{l}
\| Y_0 \|_{L^2(Q_0; \R^2)} \le C \, e^{CT} e^{\frac M{T_1}} \left( \| y_0 \|_{H^{-1}} +  \|f\|_{L^2( Q_0 ;\mathbb{R}^2)}\right) \\
	\noalign{\smallskip}
\phantom{\| Y_0 \|_{L^2(Q_0; \R^2)}} = C \, e^{CT} \, e^{\frac{Mb (a+1)}{(b-1)T}} \rho_0 (T_1) \left( \| y_0 \|_{H^{-1}} +  \|f\|_{L^2( Q_0 ;\mathbb{R}^2)}\right), \\
	\noalign{\smallskip}
\hbox{(resp., } \| Y_0 \|_{L^2(Q_0) \times C^0(\overline Q_0)} \le C\, e^{CT} \, e^{\frac{Mb (a+1)}{(b-1)T}} \rho_0 (T_1) \left( \| y_0 \|_{H^{-1} \times H_0^1} +  \|f\|_{L^2( Q_0 ;\mathbb{R}^2)}\right)).
	\end{array}
	\right.
	$$
Observe that $ \|f\|_{L^2( Q ;\mathbb{R}^2)} \le \| f \|_{\mathcal F}$ (see the expression of $\rho_{\mathcal F}$ in~\eqref{rhoF0}). Hence, repeating the previous argument, we get
	\begin{equation}\label{Y0}
	\left\{
	\begin{array}{l}
\left\| \frac{Y_0}{\rho_0} \right\|_{L^2(Q_0; \R^2)} \le C \, e^{C \left( T+ \frac 1T \right)} \left( \| y_0 \|_{H^{-1}} +  \|f\|_{\mathcal F}\right) \\
	\noalign{\smallskip}
 \hbox{(resp., } \left\| \frac{Y_0}{\rho_0} \right\|_{L^2(Q_0) \times C^0(\overline Q_k)} \le C \, e^{C \left( T+ \frac 1T \right)} \left( \| y_0 \|_{H^{-1} \times H_0^1} +  \|f\|_{\mathcal F}\right)).
	\end{array}
	\right.
	\end{equation}

\smallskip

\textbf{Case $k \ge 1$.} Again, taking into account formula $v_k = \mathcal{C}_{T_{k+1}-T_k}^{(0)}(a_k) $, Theorem~\ref{CNull_},~\eqref{ak+1<f} and~\eqref{est0}, we infer
	$$
	\begin{array}{l}
\| v_{k} \|_{L^2(T_{k},T_{k+1})} \le C \, e^{\frac{M}{T_{k+1}-T_{k}}}\|a_{k}\|_{H^{-1}} \le C \, e^{CT} \, e^{\frac{M}{T_{k+1}-T_{k}}} \| f \|_{L^{2}(Q_{k-1}; \R^2)} \\
	\noalign{\smallskip}
\phantom{\| v_{k} \|_{L^2(T_{k},T_{k+1})}} = C \, e^{CT} \, \frac{\rho_0 (T_{k+1})}{\rho_{\mathcal F}(T_{k-1})} \| f \|_{L^{2}(Q_{k-1}; \R^2)} .
	\end{array}
	$$
As in the case $k=0$, using the fact that $\rho_0$ and $\rho_{\mathcal F}$ are non-increasing functions, from the previous inequality, we deduce 
	\begin{equation}\label{vk}
\left\| \frac{v_k}{\rho_0} \right\|_{L^2(T_{k},T_{k+1})} \leq C \, e^{ CT} \left\| \frac f{\rho_{\mathcal F}} \right\|_{L^{2} (Q_{k-1}; \R^2)} , \quad \forall k \ge 1.
	\end{equation}

We can also repeat the previous argument to obtain an estimate for $Y_k$ when $k \ge 1$. From~\eqref{yhat},
	$$
	\left\{
	\begin{array}{l}
\displaystyle \| Y_k \|_{L^2(Q_k) \times C^0 (\overline Q_k)} \le C \, e^{CT} e^{\frac{M}{T_{k+1}- T_k}} \| f \|_{L^2(Q_k; \R^2)} = C \, e^{CT} \, \frac{\rho_0 (T_{k+1})}{\rho_{\mathcal F}(T_{k-1})} \| f \|_{L^2(Q_k; \R^2)} \\
	\noalign{\smallskip}	
\displaystyle \phantom{\| Y_k \|_{L^2(Q_k) \times C^0 (\overline Q_k)}} \le C \, e^{CT} \, \frac{\rho_0 (T_{k+1})}{\rho_{\mathcal F}(T_{k})} \| f \|_{L^2(Q_k; \R^2)}, 
	\end{array}
	\right.
	$$
what implies
	\begin{equation}\label{Yk}
\left\| \frac{Y_k}{\rho_0} \right\|_{L^2(Q_k) \times C^0 (\overline Q_k)} \leq C \, e^{ CT} \left\| \frac f{\rho_{\mathcal F}} \right\|_{L^2(Q_k; \R^2)}, \quad \forall k \ge 1.
	\end{equation}

\smallskip

With the functions $v_k$ and $Y_k$, $k \ge 0 $, defined above, we define
	\begin{equation}\label{vsum}
\mathcal{C}_T^{(1)} (y_0, f) :=v = \displaystyle \sum_{k \geq 0} v_k 1_{[T_k,T_{k+1})} \quad \hbox{and} \quad E_T^{(0)} (y_0, f) := Y = \displaystyle\sum_{k \geq 0} Y_k 1_{[T_k,T_{k+1})},
	\end{equation}
where $1_I$ is the characteristic function on the set $I$. Let us first remark that, by construction, $\mathcal{C}_T^{(1)}$ and $E_T^{(0)}$ are linear operators. On the other hand, recall that $Y_k = \tilde y_k + \hat y_k$, $k \ge 0$, where $\tilde y_k $ and $\hat y_k$ are respectively the solution to systems~\eqref{syst_tilde} and~\eqref{syst_hat}. So, 
	$$
Y_k (T_{k+1}^-) = a_{k+1} = \hat y_{k+1} (T_{k+1}^+) = Y_{k+1} (T_{k+1}^+), \quad \forall k \ge 0,
	$$
which implies that the function $Y$ is continuous at time $T_k$, for any $k \ge 1$, and is the solution of system~\eqref{vectPFSy} associated to $(y_0, f, v)$.

Finally, thanks to~\eqref{v0}--\eqref{Yk}, we also deduce that $ \mathcal{C}_T^{(1)} (y_0, f) \in \mathcal V$ and $E_T^{(0)} (y_0, f) \in \mathcal Y_0$ (resp., $E_T^{(0)} (y_0, f) \in \mathcal Y$) for any $(y_0, f) \in H^{-1}(0, \pi; \R^2) \times \mathcal F$ (resp., for any $(y_0, f) \in H^{-1}(0, \pi) \times H_0^{1}(0, \pi) \times \mathcal F$) and
	$$
	\left\{
	\begin{array}{l}
\left\| \mathcal{C}_T^{(1)} (y_0, f) \right\|_{\mathcal V} = \left\| v \right\|_{\mathcal V} \le C \, e^{C \left( T+ \frac 1T \right)} \left( \| y_0 \|_{H^{-1}} +  \|f\|_{\mathcal F}\right), \\
	\noalign{\smallskip}
\left\| E_T^{(0)} (y_0, f) \right\|_{\mathcal Y_0} = \left\| Y \right\|_{\mathcal Y_0} \le C \, e^{C \left( T+ \frac 1T \right)} \left( \| y_0 \|_{H^{-1}} +  \|f\|_{\mathcal F}\right), \quad \forall (y_0, f) \in H^{-1}(0, \pi; \R^2) \times \mathcal F,
	\end{array}
	\right.
	$$
(resp.,
	$$
\left\| E_T^{(0)} (y_0, f) \right\|_{\mathcal Y} = \left\| Y \right\|_{\mathcal Y} \le C  e^{C \left( T+ \frac 1T \right)} \left( \| y_0 \|_{H^{-1}\times H_0^1} +  \|f\|_{\mathcal F}\right), \ \forall (y_0, f) \in H^{-1}(0, \pi) \times H_0^{1}(0, \pi) \times \mathcal F \hbox{)}.
	$$

The above estimates provide the proof of Proposition~\ref{boundsM5.1}. This ends the proof.
\end{proof}


\subsection{Proof of Theorem~\ref{NullC}}\label{s5.2}

We will devote this section to proving the local exact controllability at time $T>0$ of the phase-field system~\eqref{PFSy} stated in Theorem~\ref{NullC}. To this objective, let us take 
	$$
\tilde y_0 = (\tilde{\theta}_0,\tilde{\phi}_0)\in H^{-1}(0,\pi)\times (c + H_0^1(0,\pi))
	$$ 
($c= \pm 1$). As we saw in Section~\ref{s1}, the local exact controllability of system~\eqref{PFSy} at time $T$ to the constant trajectory $(0,c)$ is equivalent to the local null controllability of system~\eqref{bisPFSy} at time $T$ with $y_0 = (\theta_0, \phi_0) = (\tilde{\theta}_0, \tilde{\phi}_0 - c) \in H^{-1}(0,\pi)\times H_0^1(0,\pi)$ (the nonlinear functions $g_1$ and $g_2$ are given in~\eqref{g1g2}).

Let us take $a,b > 1$ (which will be determined below) and consider the functions $\rho_{\mathcal F}$ and $\rho_0$, defined in~\eqref{rhoF0}, and the spaces $\mathcal{F}$, $\mathcal{V}$ and $\mathcal{Y}$ given in~\eqref{FVY}. In order to prove the local null controllability result at time $T$ for system~\eqref{bisPFSy} we will perform a fixed-point strategy in the space $\mathcal Y$ which, in particular, will prove the existence of a control $v \in \mathcal V$ such that system~\eqref{bisPFSy} has a solution $y \in \mathcal Y$ associated to $(v, y_0)$. The condition $y \in \mathcal Y$ will imply the null controllability result for this system.

Let us fix $\varepsilon>0$ (to be determined bellow). With the previous data and notations, we consider the closed ball in the space $\mathcal F$
	\begin{equation*}
\overline{B}_{\varepsilon} = \left\{ f\in\mathcal{F}: \left\| f \right\|_{\mathcal F} \leq \varepsilon \right\}.
	\end{equation*}
Observe that if the initial datum $\tilde y_0 \in H^{-1}(0,\pi)\times (c + H_0^1(0,\pi))$ satisfies~\eqref{epsi}, then $y_0 = (\theta_0, \phi_0) = (\tilde{\theta}_0, \tilde{\phi}_0 - c) \in H^{-1}(0, \pi) \times H_0^1 (0, \pi)$ satisfies
	\begin{equation}\label{epsi2}
\| \theta_ 0 \|_{H^{-1}} + \| \phi_ 0 \|_{H_0^{1}} \le \varepsilon.
	\end{equation}

For each $f\in \overline B_\varepsilon \subset \mathcal{F}$, we denote $v_f = \mathcal{C}_T^{(1)}(y_0,f)\in \mathcal V$ and $y_f = (\theta_f,\phi_f):=E_T^{(1)}(y_0,f)\in \mathcal Y$, where the operators $\mathcal{C}_T^{(1)}$ and $E_T^{(1)}$ are given in Theorem~\ref{boundsM5.1}. As a consequence of this result and~\eqref{epsi2}, one has
	\begin{equation}\label{epsi3}
\left\| y_f \right\|_{\mathcal Y} + \left\| v_f \right\|_{\mathcal V} \leq C \, e^{C \left( T + \frac 1T \right)} \left( \|y_0\|_{H^{-1} \times H_0^1 } + \left\| f \right\|_{\mathcal{F}} \right) \le C \, e^{C \left( T + \frac 1T \right)} \, \varepsilon, \quad \forall f \in  \overline B_\varepsilon,
	\end{equation}
for a positive constant $C = C(\xi, \rho, \tau)$. Thus, we define the nonlinear operator $\mathcal{N} : \overline{B}_{\varepsilon} \rightarrow C^0(\overline Q_T; \R^2)$ given by (see~\eqref{g1g2})
	\begin{equation}\label{mthcalN_}
\mathcal{N}(f)=\left(
	\begin{array}{l}
\displaystyle \pm \frac{3\rho}{4\tau}\phi_f^2 + \frac{\rho}{4\tau}\phi_f^3 \\
	\noalign{\smallskip}
\displaystyle \mp \frac{3}{2\tau} \phi_f^2 - \frac{1}{2\tau}\phi_f^3
	\end{array}
	\right).
	\end{equation}
It is clear that the operator $\mathcal{N}$ is well-defined. On the other hand, if $\mathcal{N}$ admits a fixed point $f \in \mathcal F$, then $y_f \in \mathcal Y$, together with $v_f \in \mathcal V$, provides a solution of the system~\eqref{bisPFSy} associated to the initial datum $y_0 = (\theta_0, \phi_0)$. In fact, from Proposition~\ref{linearexist}, $y_f \in C^0([0,T]; H^{-1}(0,\pi; \R^2))$. Finally, condition $y_f \in C^0([0,T]; H^{-1}(0,\pi; \R^2)) \cap \mathcal Y $ in particular implies the null controllability result for system~\eqref{bisPFSy}. This would prove Theorem~\ref{NullC}.

The next task is to prove that the operator $\mathcal N$ has a fixed-point in the complete metric space $\overline B_\varepsilon \subset \mathcal{F}$. To this end, we will apply the Banach Fixed-Point Theorem. Before, let us select any $a > 1$ and $b$ such that
	$$
b^2 \in \left(1,\frac{2a}{a+1}\right).
	$$
With this choice, the functions $\rho_0^2/ \rho_\mathcal{F}$ and ${\rho_0^3}/{\rho_\mathcal{F}}$ are uniformly bounded in $[0,T]$, i.e., there exists a constant $C_T > 0$, depending on $T$, such that
	\begin{equation*}
\left\| \frac{\rho_0^2}{\rho_\mathcal{F}}\right\|_{C^0[0,T]} \leq C_T \quad \mbox{and} \quad \left\|\dfrac{\rho_0^3}{\rho_\mathcal{F}}\right\|_{C^0[0,T]} \leq C_T.
	\end{equation*}

Let us now check the assumptions of the Banach Fixed-Point Theorem:

\smallskip

\noindent
\textbf{1.} $\mathcal{N}(\overline{B}_{\varepsilon}) \subset \overline{B}_{\varepsilon}$: Indeed, if $f\in \overline{B}_{\varepsilon}$, then, from~\eqref{epsi3}, we obtain
	\begin{equation*}
	\begin{array}{l}
\displaystyle \left\| \mathcal{N}(f) \right\|_{\mathcal F} \leq C_T \left\|\dfrac{\mathcal{N}(f)}{\rho_\mathcal{F}}\right\|_{C^0(\overline{Q}_T ;\mathbb{R}^2)} \leq C_T \left(  \left\|\dfrac{\phi_f^2}{\rho_\mathcal{F}}\right\|_{C^0(\overline{Q}_T)} + \left\|\dfrac{\phi_f^3}{\rho_\mathcal{F}}\right\|_{C^0(\overline{Q_T})} \right) 
	\\
	\noalign{\smallskip}
\displaystyle \phantom{\displaystyle \left\| \mathcal{N}(f) \right\|_{\mathcal F}} \leq C_T \left( \left\|\frac{\rho_0^2}{\rho_\mathcal{F}}\right\|_{C^0(\overline{Q}_T)} \left\|\frac{\phi_f}{\rho_0}\right\|_{C^0 (\overline{Q}_T)}^2 + \left\|\dfrac{\rho_0^3}{\rho_\mathcal{F}}\right\|_{C^0(\overline{Q}_T)} \left\|\dfrac{\phi_f}{\rho_0}\right\|_{C^0 (\overline{Q_T})}^3 \right)
	\\
	\noalign{\smallskip}
\displaystyle \phantom{\displaystyle \left\| \mathcal{N}(f) \right\|_{\mathcal F}} \leq {C}_T\left( \left\| y_f \right\|_{ \mathcal Y}^2 + \left\| y_f \right\|_{ \mathcal Y}^3 \right)  \leq C_T \, e^{C \left( T + \frac 1T \right)} \left( \varepsilon^2 + \varepsilon^3 \right) \le \varepsilon 
	\end{array}
	\end{equation*}
for $\varepsilon=\varepsilon(T)$ small enough.

\smallskip

\noindent
\textbf{2.} $\mathcal{N}$ is a contraction mapping: Let us take $f_1, f_2 \in \overline B_\varepsilon \subset \mathcal F$ and denote $ y_i = (\theta_i, \phi_i) = E_T^{(1)}(y_0,f_i)\in \mathcal Y$, $i=1,2$. Firstly, observe that the non linearity $(g_1,g_2)$, given in~\eqref{g1g2}, satisfies
	\begin{equation*}
|g_j(s_1) - g_j(s_2)| \leq C( |s_1|^2 + |s_2|^2 + |s_1| + |s_2| )|s_1- s_2|, \quad \forall s_1, s_2 \in \R, \quad j=1,2.
	\end{equation*}
Thus, using again~\eqref{epsi3} and Theorem~\ref{boundsM5.1}, we have
	\begin{equation*}
	\begin{array}{l}
\displaystyle \left\| \mathcal{N}(f_1) - \mathcal{N}(f_2) \right\|_{\mathcal{F}} \leq C_T \sum_{j=1}^2 \left\|\frac{g_j(\phi_1) - g_j(\phi_2)}{\rho_\mathcal{F}}\right\|_{C^0(\overline{Q}_T)} 
	\\
	\noalign{\smallskip}
\displaystyle \phantom{\displaystyle \left\| \mathcal{N}(f_1) - \mathcal{N}(f_2) \right\|_{\mathcal{F}}} \leq C_T \left\| \frac{\rho_0}{\rho_\mathcal{F}} \left(|\phi_1|^2 + |\phi_2|^2 + |\phi_1| + |\phi_2| \right) \dfrac{|\phi_1-\phi_2|}{\rho_0} \right\|_{C^0(\overline{Q}_T)}  
	\\
	\noalign{\smallskip}
\displaystyle \phantom{\displaystyle \left\| \mathcal{N}(f_1) - \mathcal{N}(f_2) \right\|_{\mathcal{F}}} \leq C_T \left\| \left( \left| \dfrac{\phi_1}{\rho_0} \right|^2 + \left| \dfrac{\phi_2}{\rho_0} \right|^2 \right)\dfrac{\rho_0^3}{\rho_{\mathcal{F}}} + \left( \left| \dfrac{\phi_1}{\rho_0} \right| + \left| \dfrac{\phi_2}{\rho_0} \right| \right)\dfrac{\rho_0^2}{\rho_{\mathcal{F}}} \right\|_{C^0(\overline{Q}_T)} \left\| \dfrac{\phi_1-\phi_2}{\rho_0} \right\|_{C^0(\overline{Q}_T)}	\\
	\noalign{\smallskip}
\displaystyle \phantom{\displaystyle \left\| \mathcal{N}(f_1) - \mathcal{N}(f_2) \right\|_{\mathcal{F}}} \leq C_T \left( \left\| y_1 \right \|^2_{ \mathcal Y} + \left\| y_2 \right\|^2_{ \mathcal Y} + \left\| y_1 \right\|_{ \mathcal Y} + \left\| y_2\right\|_{ \mathcal Y} \right) \left\| E_T^{(1)}(y_0,f_1) - E_T^{(1)}(y_0,f_2) \right\|_{ \mathcal Y}
	\\
	\noalign{\smallskip}
\displaystyle \phantom{\displaystyle \left\| \mathcal{N}(f_1) - \mathcal{N}(f_2) \right\|_{\mathcal{F}}} \leq C_T \, e^{C \left( T + \frac 1T \right)} \left( \varepsilon^2 + \varepsilon \right) \left\| f_1-f_2 \right\|_{\mathcal{F}}.				
	\end{array}
	\end{equation*}
From this inequality it is clear that we can choose $\varepsilon = \varepsilon (T)$ (small enough) in such a way that $\mathcal{N}$ is a contraction mapping. 

\smallskip

In conclusion, we can apply the Banach Fixed-Point Theorem. This proves that the operator $\mathcal N$ has a fixed-point and provides the proof of Theorem~\ref{NullC}.

\end{section}
%

\begin{appendices}
	\appendix
\appendixpage
\numberwithin{equation}{section}

\begin{section}{}\label{APPX1_}
	
\setcounter{equation}{0}

This appendix will be devoted to dealing with the existence and uniqueness of solution of the linear systems~\eqref{AdjvectPFSy} and~\eqref{vectPFSy}. To be precise, we will prove Propositions~\ref{adjexist} and~\ref{linearexist}.

\begin{proof}[Proof of Proposition~\ref{adjexist}]
Let us assume that $\varphi_0 \in H_0^1(0,\pi;\mathbb{R}^2)$ and $g\in L^2(Q_T;\mathbb{R}^2)$. Let us denote  $\varphi_0 = (\theta_0,\phi_0)$ and $g=(g_1,g_2)$. Then the system \eqref{AdjvectPFSy} can be write as
	\begin{equation*}
	\left\{
	\begin{array}{ll}
-\theta_t - \xi\theta_{xx} + \dfrac{\rho}{\tau}\theta  - \dfrac{2}{\tau}\phi  =g_1 & \mbox{in } Q_T,  \\
	\noalign{\smallskip}
-\phi_t - \xi\phi_{xx} + \dfrac{1}{2}\rho\xi\theta_{xx}  - \dfrac{\rho}{2\tau}\theta + \dfrac{1}{\tau}\phi  = g_2	 & \mbox{in }  Q_T, \\
	\noalign{\smallskip}
\theta(0,\cdot) =  \phi(0,\cdot) = \theta(\pi,\cdot)=\phi(\pi,\cdot) = 0 & \mbox{on } (0,T), \\
	\noalign{\smallskip}
\theta(\cdot,T) = \theta_0, \quad \phi(\cdot,T) =\phi_0 & \mbox{in } (0,\pi),
	\end{array}
	\right.
	\end{equation*}
where $\varphi = (\theta,\phi)$. On the other hand, $ \xi\theta_{xx} = -\theta_t +\dfrac{\rho}{\tau}\theta - \dfrac{2}{\tau}\phi - g_1$. Thus, the previous system becomes
	\begin{equation}\label{PFSylinearbakApx_2}
	\left\{
	\begin{array}{ll}
-\theta_t - \xi\theta_{xx} - \dfrac{2}{\tau}\phi + \dfrac{\rho}{\tau}\theta  =g_1 & \mbox{in } Q_T,  \\
	\noalign{\smallskip}
-\phi_t - \xi\phi_{xx} - \dfrac{\rho}{2}\theta_t - \dfrac{\rho-1}{\tau}\phi + \dfrac{\rho(\rho-1)}{2\tau}\theta  = \dfrac{\rho}{2} g_1 + g_2	 & \mbox{in }  Q_T, \\
	\noalign{\smallskip}
\theta(0,\cdot) =  \phi(0,\cdot) = \theta(\pi,\cdot)=\phi(\pi,\cdot) = 0 & \mbox{on } (0,T), \\
	\noalign{\smallskip}
\theta(\cdot,T) = \theta_0, \quad \phi(\cdot,T) =\phi_0 & \mbox{in } (0,\pi),
	\end{array}
	\right.
	\end{equation}

Then, Proposition \ref{adjexist} is equivalent to prove that the system~\eqref{PFSylinearbakApx_2} has a unique strong solution $(\theta,\phi)$ satisfying
	$$
\theta,\phi\in C^0([0,T];H_0^1(0,\pi))\cap L^2(0,T;H^2(0,\pi) \cap H_0^1(0,\pi))
	$$ 
and
	\begin{equation}\label{acot1Apx_}
	\begin{array}{l}
\|\theta\|_{C^0(H_0^1)} + \|\phi\|_{C^0(H_0^1)} + \|\theta\|_{L^2(H^2\cap H_0^1)} + \|\phi\|_{L^2(H^2\cap H_0^1)} \\
	\noalign{\smallskip}
\phantom{\|\theta\|_{C^0(H_0^1)} + \|\phi\|_{C^0(H_0^1)}} \leq e^{CT}  \left( \|g_1\|_{L^2(L^2)} + \|g_2\|_{L^2(L^2)} +\|\theta_0\|_{H_0^1} + \|\phi_0\|_{H_0^1} \right).
	\end{array}
	\end{equation}
for a positive constant $C$, only depending on $\xi$, $\rho$ and $\tau$.

We will use the well-known Faedo-Galerkin method. First, let us consider the orthonormal basis $\{\eta_n\}_{n\in \N}$ of $L^2 (0,\pi)$ ($\eta_n$ is the normalized eigenfunction of the Dirichlet-Laplace operator, see~\eqref{sin}). For each $m\in\mathbb{N}$, we consider $V_{m}= [\eta_1,\eta_2,\cdots,\eta_m ]$, the subspace generated by the first $m$ vectors of $\{\eta_{n}\}_{n \in\mathbb{N}}$. Let us also consider $P_m$, the orthogonal projection operator onto the finite-dimensional space $V_m$ in $L^2(0, \pi)$. If we define 
	\begin{equation}\label{apx1}
\theta_0^m = P_m \theta_0, \quad \phi_0^m = P_m \phi_0, \quad g_1^m(\cdot, t) = P_m g_1(\cdot, t) \quad \hbox{and} \quad g_2^m(t, \cdot) = P_m g_2(t, \cdot), 
	\end{equation}
one has $\theta_0^m, \phi_0^m \in V_m$ and $g^m_1, g^m_2 \in L^2 (0,T; V_m)$, for any $m \in \N$, and
	\begin{equation}\label{apx2}
\theta_0^m\rightarrow \theta_0 , \quad \phi_0^m\rightarrow \phi_0 \mbox{ in } H_0^1(0,\pi),  \quad \mbox{and} \quad g_1^m\rightarrow g_1 , \quad g_2^m\rightarrow g_2 \mbox{ in } L^2(Q_T), \quad \mbox{as } m\rightarrow\infty .
	\end{equation}

We want an approximate solution $(\theta^m,\phi^m)\in C^0([0,T]; V_{m}^{2})$ of the approximate problem
	\begin{equation}\label{PFSylinearbakApx_2m}
	\left\{
	\begin{array}{ll}
-\theta^m_t - \xi\theta^m_{xx} - \dfrac{2}{\tau}\phi ^m + \dfrac{\rho}{\tau}\theta^m  =g_1^m & \mbox{in } Q_T,  \\
	\noalign{\smallskip}
-\phi ^m_t - \xi\phi ^m_{xx} - \dfrac{\rho}{2}\theta^m_t - \dfrac{\rho-1}{\tau}\phi ^m + \dfrac{\rho(\rho-1)}{2\tau}\theta^m  = \dfrac{\rho}{2} g_1^m + g_2^m & \mbox{in }  Q_T, \\
	\noalign{\smallskip}
\theta^m(0,\cdot) =  \phi ^m(0,\cdot) = \theta^m(\pi,\cdot)=\phi ^m(\pi,\cdot) = 0 & \mbox{on } (0,T), \\
	\noalign{\smallskip}
\theta^m(\cdot,T) = \theta^m_0, \quad \phi ^m(\cdot,T) =\phi ^m_0 & \mbox{in } (0,\pi),
	\end{array}
	\right.
	\end{equation}
under the form
	$$
\theta^m(x,t) = \sum_{j=1}^{m}\alpha_{jm}(t) \eta_{j}(x),\quad \phi^m(x,t)=\sum_{j=1}^{m}\beta_{jm}(t)\eta_{j}(x), \quad (x,t) \in Q_T.
	$$

It is clear that, for any $m \ge 1$, system~\eqref{PFSylinearbakApx_2m} is equivalent to a Cauchy problem for a linear ordinary differential system for the variables $\alpha_{jm}$ and $\beta_{jm}$, $1 \le j \le m$. In consequence, system~\eqref{PFSylinearbakApx_2m} admits a unique solution $(\theta^m,\phi^m)\in C^0([0,T]; V_{m}^{2})$ with $(\theta^m_t,\phi^m_t)\in L^2(0,T; V_{m}^{2})$.   

The proof of Proposition~\ref{adjexist} can be easily deduced from appropriate estimates of the approximate solution $(\theta^m,\phi^m)$ of system~\eqref{PFSylinearbakApx_2m}. 

If we multiply the first equation in~\eqref{PFSylinearbakApx_2m} by $-\tfrac{\rho}{2}\theta_t^m$, the second one by $\tfrac{2}{\tau}\phi^m$, we integrate on the interval $(0, \pi)$ and we add both equalities, we get,
	\begin{equation*}
	\begin{array}{cc}
\int_0^\pi \left( \dfrac{\rho}{2}|\theta_t^m|^2
- \dfrac{1}{\tau} \dfrac{d}{dt}(|\phi ^m|^2)
- \dfrac{\rho}{\tau}\theta^m_t \phi^m
\right) \, dx
+ \int_0^\pi \left ( 
-\dfrac{\rho\xi}{4}\dfrac{d}{dt}(|\theta^m_{x}|^2)
+\dfrac{2\xi}{\tau} |\phi ^m_x|^2 
\right) \, dx \\
	\noalign{\smallskip}
+ \int_0^\pi \left( 
\dfrac{\rho}{\tau} \phi ^m\theta_t^m
- \dfrac{2(\rho-1)}{\tau^2} |\phi ^m|^2
\right) \, dx
+ \int_0^\pi\left(
- \dfrac{\rho^2}{2\tau} \theta^m \theta^m_t
+ \dfrac{\rho(\rho-1)}{\tau^2} \theta^m \phi^m
\right) \,dx \\
	\noalign{\smallskip}
= -\dfrac{\rho}{2}\int_0^\pi g_1^m \theta_t^m dx
-\dfrac{2}{\tau}\int_0^\pi \left(\dfrac{\rho}{2}g_1^m + g_2^m \right)\phi^m \, dx.
	\end{array}
	\end{equation*}
Applying the Cauchy-Schwarz inequality in the previous equality, we obtain
	\begin{equation*}
	\begin{array}{c}
\dfrac{\rho}{2}\|\theta_t^m(\cdot,t)\|^2_{L^2 }
+\dfrac{2\xi}{\tau} \|\phi ^m_x(\cdot,t)\|^2_{L^2 }
- \dfrac{d}{dt}	\left(
\dfrac{1}{\tau} \|\phi ^m(\cdot,t)\|^2_{L^2 } 
+\dfrac{\rho\xi}{4}\|\theta^m_{x}(\cdot,t)\|^2_{L^2 }
\right) \leq \dfrac{\rho}{4}\|\theta_t^m(\cdot,t)\|^2_{L^2 } \\
	\noalign{\smallskip}
\phantom{}
+ C\left(
\|\theta^m(\cdot,t)\|^2_{L^2 }
+\|\phi ^m(\cdot,t)\|^2_{L^2 }
+ \|g_1^m(\cdot,t)\|^2_{L^2 }
+ \| g_2^m  (\cdot,t)\|^2_{L^2 }
\right), \hbox{ a.e.~}t \in (0,T), 
	\end{array}
	\end{equation*}
for a constant $C>0$ depending on the parameters $\xi$, $\rho$ and $\tau$. Using Poincar\'e inequality, it follows
	\begin{equation*}
	\begin{array}{ll}
\|\theta_t^m(\cdot,t)\|^2_{L^2 }
+\|\phi ^m_x(\cdot,t)\|^2_{L^2 }
-\dfrac{d}{dt}	\left(
\|\phi ^m(\cdot,t)\|^2_{L^2 } 
+\|\theta^m_{x}(\cdot,t)\|^2_{L^2 }
\right) \\
	\noalign{\smallskip}
\phantom{\|\theta_t^m(\cdot,t)\|^2}
\leq
C\left(
\|\phi ^m(\cdot,t)\|^2_{L^2 }
+\|\theta_x^m(\cdot,t)\|^2_{L^2 }
+ \|g_1^m (\cdot,t)\|^2_{L^2 }
+ \| g_2^m (\cdot,t)\|^2_{L^2 }
\right),
	\end{array}
	\end{equation*}
for a new constant $C>0$. Multiplying the previous inequality by $e^{-C(T-t)}$ and integrating in the interval $[t,T]$, with $t<T$, we have
	\begin{equation*}
	\begin{array}{ll}
\int_t^Te^{-C(T-s)}
\left(\|\theta_t^m(\cdot,s)\|^2_{L^2 }
+\|\phi ^m_x(\cdot,s)\|^2_{L^2 } \right) \, ds
+ e^{-C(T-t)}\left(
\|\phi ^m(\cdot,t)\|^2_{L^2 } 
+\|\theta^m_{x}(\cdot,t)\|^2_{L^2 }
\right)\\
	\noalign{\smallskip}
\phantom{\|\theta_t^m(\cdot,t)\|^2}
\leq
\|\phi ^m_0\|^2_{L^2 }			
+\|(\theta^m_0)_x\|^2_{L^2 }
+\int_t^Te^{-C(T-s)}\left(
\|g_1^m (\cdot,s)\|^2_{L^2 }
+ \| g_2^m  (\cdot,s)\|^2_{L^2 }
\right) \, ds.
	\end{array}
	\end{equation*}
Finally, multiplying the previous inequality by $e^{C(T-t)}$ and taking maximum with $t \in [0,T]$, we can deduce
	\begin{equation}\label{apxI02_}
	\left\{
	\begin{array}{ll}
\|\theta_t^m\|^2_{L^2(Q_T)}
+\|\phi ^m\|^2_{L^2(H_0^1 )}
+ \|\phi ^m\|^2_{C^0(L^2 )} 
+\|\theta^m\|^2_{C^0(H_0^1 )} \\
	\noalign{\smallskip}
\phantom{\|\theta_t^m\|^2_{L^2(Q_T)}}
\leq
e^{CT}\left(
\|\phi^m_0\|^2_{L^2 }			
+\|\theta^m_0\|^2_{H_0^1}
+\|g_1^m \|^2_{L^2(Q_T)}
+ \| g_2^m \|^2_{L^2(Q_T)}
\right)
\\
	\noalign{\smallskip}
\phantom{\|\theta_t^m\|^2_{L^2(Q_T)}}
\leq
e^{CT}\left(
\|\phi_0\|^2_{L^2 }			
+\|\theta_0\|^2_{H_0^1}
+\|g_1 \|^2_{L^2(Q_T)}
+ \| g_2 \|^2_{L^2(Q_T)}
\right).
	\end{array}
	\right.
	\end{equation}
Observe that in the previous inequalities we have used~\eqref{apx1}.

Let us notice that, from the first equation in~\eqref{PFSylinearbakApx_2m}, 
	\begin{equation}\label{apxI04_}
	\begin{array}{l}
\displaystyle \|\theta^m_{xx}\|_{L^2(Q_T)} = \frac 1\xi \left\| - \theta^m_t - \frac{2}{\tau}\phi ^m + \frac{\rho}{\tau}\theta^m  - g_1\right \|_{L^2(Q_T)}\\
	\noalign{\smallskip}
\displaystyle \phantom{\displaystyle \|\theta^m_{xx}\|_{L^2(Q_T)} } \leq
e^{CT}\left(
\|\phi_0\|^2_{L^2 }			
+\|\theta_0\|^2_{H_0^1}
+\|g_1 \|^2_{L^2(Q_T)}
+ \| g_2 \|^2_{L^2(Q_T)}
\right).
	\end{array}
	\end{equation}
From~\eqref{apxI02_} and~\eqref{apxI04_}, we get that the sequences $\{\theta^m\}_{m\in\N}$ and $\{\theta^m_t\}_{m\in\N }$ are respectively bounded in $ L^2(0,T;H^2(0,\pi)\cap H_0^1(0,\pi)) \cap C^0([0,T]; H_0^1(0,\pi))$ and $ L^2(Q_T)$. Then, there exist a subsequence, still denoted $\{\theta^m\}_{m\in\mathbb{N}}$, and a function $\theta\in L^\infty (0,T;H_0^1(0,\pi))\cap L^2(0,T;H^2(0,\pi)\cap H_0^1(0,\pi))$ such that $\theta_t \in L^2(Q_T)$ and
	\begin{equation}\label{apxI10_}
	\left\{
	\begin{array}{l}
\theta^m \stackrel{*}{\rightharpoonup}  \theta \quad \mbox{weakly-* in}\ L^\infty (0,T; H_0^1(0,\pi)), \quad  \theta^m_t \rightharpoonup \theta_t \quad \mbox{weakly in}\ L^2 (Q_T),\\
	\noalign{\smallskip}
\theta^m \rightharpoonup \theta \quad \mbox{weakly in}\ L^2(0,T;H^2(0,\pi)\cap H_0^1(0,\pi)) .
	\end{array}
	\right.
	\end{equation}
Observe that the previous regularity for function $\theta$ also implies $\theta \in C^0([0,T]; H_0^1(0, \pi))$.

In order to deal with $\phi^m$, let us multiply the second equation in~\eqref{PFSylinearbakApx_2m} by $-\phi^m_t$ and integrate on the interval $(0, \pi)$. After an integration by parts, we deduce
	$$
	\begin{array}{l}
\| \phi_t^m (\cdot,t) \|^2_{L^2} - \frac 1{2 \xi } \frac d {dt}\| \phi_x^m (\cdot,t) \|^2_{L^2} = \frac \rho 2 \int_0^\pi \theta_t^m \phi_t^m \,dx + \frac{\rho -1}\tau \int_0^\pi \phi^m \phi_t^m \,dx \\
	\noalign{\smallskip}
\phantom{\| \phi_t^m (\cdot,t) \|^2_{L^2}- \frac 1{2 \xi } \frac d {dt}\| \phi_t^m (\cdot,t) \|^2_{L^2}} - \frac{\rho (\rho -1)}{2\tau} \int_0^\pi \theta^m \phi_t^m \,dx +  \int_0^\pi \left( \frac \rho 2 g_1^m + g_2^m \right) \phi_t^m \,dx.
	\end{array}
	$$
Using again Cauchy-Scharwz inequality, we also obtain
	$$
	\left\{
	\begin{array}{l}
\| \phi_t^m (\cdot,t) \|^2_{L^2} - \frac d {dt} \| \phi_x^m (\cdot,t) \|^2_{L^2} \le C\left( \|\phi ^m(\cdot,t)\|^2_{L^2 } +  \|\theta^m(\cdot,t)\|^2_{L^2 } \right. \\
	\noalign{\smallskip}
\phantom{\| \phi_t^m (\cdot,t) \|^2_{L^2} - \frac d {dt} \| \phi_x^m (\cdot,t) \|^2_{L^2}} + \left. \|\theta_t^m(\cdot,t)\|^2_{L^2 } + \|g_1^m (\cdot,t)\|^2_{L^2 } + \| g_2^m (\cdot,t)\|^2_{L^2 }\right), \hbox{ a.e.~}t \in (0,T).
	\end{array}
	\right.
	$$
Reasoning as before, using inequality~\eqref{apxI02_} and again the second equation in~\eqref{PFSylinearbakApx_2m}, we deduce $\phi_t^m, \phi_{xx}^m \in L^2 (Q_T)$, $\phi^m \in C^0([0,T]; H_0^1(0, \pi))$ and 
	\begin{equation}\label{apxI05_}
\| \phi_t^m \|_{L^2 (Q_T)} + \| \phi^m \|_{L^2 (H^2 \cap H_0^1)} + \| \phi^m \|_{C^0 (H_0^1)} \le e^{CT}\left( \|\phi_0\|^2_{H_0^1 } + \|\theta_0\|^2_{H_0^1} + \|g_1 \|^2_{L^2(Q_T)} + \| g_2 \|^2_{L^2(Q_T)} \right)
	\end{equation}
As before, inequality~\eqref{apxI05_} allows us to extract a new subsequence (still denoted with the index $m$) and a function $\phi \in L^2(0,T; H^2(0,\pi) \cap H_0^1(0,\pi)) \cap C^0([0,T]; H_0^1(0, \pi))$ such that $\phi_t \in L^2 (Q_T)$ and
	\begin{equation}\label{apxI06_}
	\left\{
	\begin{array}{l}
\phi^m \stackrel{*}{\rightharpoonup}  \phi \quad \mbox{weakly-* in}\ L^\infty (0,T; H_0^1(0,\pi)), \quad  \phi^m_t \rightharpoonup \phi_t \quad \mbox{weakly in}\ L^2 (Q_T),\\
	\noalign{\smallskip}
\phi^m \rightharpoonup \phi \quad \mbox{weakly in}\ L^2(0,T;H^2(0,\pi)\cap H_0^1(0,\pi)) .
	\end{array}
	\right.
	\end{equation}

Finally, using the convergences in~\eqref{apx2},~\eqref{apxI10_} and~\eqref{apxI06_}, we can verify standardly that $(\theta,\phi)$ is a strong solution of the system~\eqref{PFSylinearbakApx_2}. In addition, inequality~\eqref{acot1Apx_} can be obtained combining the inequalities~\eqref{apxI02_},~\eqref{apxI04_} and~\eqref{apxI05_}. This proves the proposition. 
\end{proof}

\begin{proof}[Proof of Proposition~\ref{linearexist}] 

Let us take $y_0 \in H^{-1}(0, \pi)$, $v\in L^2(0,T)$ and $f \in L^2(Q_T; \R^2)$ and consider the functional $\mathcal{G} : L^2(Q_T;\mathbb{R}^2) \rightarrow \mathbb{R}$ given by
	\begin{equation*}
\mathcal{G} (g) = \left\langle y_0,\varphi(\cdot,0) \right\rangle - \int_0^T B^*D^* \varphi_x(0,t)v(t) \, dt + \intdoble_{Q_T} f \cdot \varphi \, dx\, dt.
	\end{equation*}
where $\varphi \in C^0([0,T];H_0^1(0,\pi;\mathbb{R}^2))\cap L^2(0,T;H^2 (0,\pi;\mathbb{R}^2)\cap H_0^1(0,\pi;\mathbb{R}^2))$ is the solution of \eqref{AdjvectPFSy} associated to $g$ and $\varphi_0=0$. 
From Proposition~\ref{adjexist}, we infer that $\mathcal{G}$ is bounded. In fact, from~\eqref{acot1} we can deduce the existence of a positive constant $C$, only depending on $D$ and $A$, such that
	\begin{equation*}
| \mathcal{G} (g)| \leq e^{CT}\left( \|y_0\|_{H^{-1}} + \|v\|_{L^2(0,T)} + \|f\|_{L^2(L^2)}\right) \|g\|_{L^2(L^2)},
	\end{equation*}
for all $g\in L^2(Q_T;\mathbb{R}^2)$. Then, by the Riesz Representation Theorem, there exists a unique function $y\in L^2(Q_T;\mathbb{R}^2)$ satisfying \eqref{transp}, i.e., a solution by transposition of \eqref{vectPFSy} in the sense of Definition~\ref{deft_}. Moreover, 
	\begin{equation*}
\|y\|_{L^2(L^2)}  = \|\mathcal{G}\| \leq e^{CT}\left( \|y_0\|_{H^{-1}}+ \|v\|_{L^2(0,T)} +\|f\|_{L^2(L^2)}\right),
	\end{equation*}
and $y$ satisfies the equality $ y_t - D y_{xx} + A y = f $ in $\mathcal{D}' (Q_T; \R^2)$.

Let us now see that the solution $y$ of the system~\eqref{vectPFSy} is more regular. To be precise, let us see that $ y_{xx} \in L^2 (0,T; ( H^2(0,\pi;\mathbb{R}^2) \cap H_0^1(0,\pi;\mathbb{R}^2) )' ) $ and 
	\begin{equation}\label{a60}
\| y_{xx} \|_{L^2 ((H^2 \cap H_0^1)')} \leq e^{CT}\left( \|y_0\|_{H^{-1}} + \|v\|_{L^2(0,T)} + \|f\|_{L^2(L^2)}\right) ,
	\end{equation}
for a new constant $C>0$ (only depending on $D$ and $A$). To this end, let us take two sequences $\{ y_0^n \}_{n \ge 1} \subset H_0^1(0, \pi; \R^2)$ and $\{ v^n \}_{n \ge 1} \in H_0^1(0,T)$ such that
	$$
y_0^n \to y_0 \hbox{ in } H^{-1}(0,\pi; \R^2) \quad \hbox{and} \quad v^n \to v \hbox{ in } L^2(0,T).
	$$
With the previous regularity assumption it is possible to show that system~\eqref{vectPFSy} for $y_0^n$, $v^n$ and $f$ has a unique strong solution $y_n \in C^0([0,T];H_0^1(0,\pi;\mathbb{R}^2))\cap L^2(0,T;H^2 (0,\pi;\mathbb{R}^2)\cap H_0^1(0,\pi;\mathbb{R}^2))$ which satisfies
	\begin{equation*}
\intdoble_{Q_T} y_n \cdot g \, dx\, dt = \left\langle y_0^n ,\varphi(\cdot,0) \right\rangle - \int_0^T B^*D^* \varphi_x(0,t)v^n (t) \, dt + \intdoble_{Q_T} f \cdot \varphi \, dx\, dt, \quad \forall n \ge 1,
	\end{equation*}
for any $g \in L^2(Q_T; \R^2)$, where $\varphi$ is the solution of the system~\eqref{AdjvectPFSy} associated to $g$ and $\varphi_0 = 0$. Indeed, if we take the new function $\widetilde y_n (\cdot, t) = y_n(\cdot, T-t) -(v^n(T-t), 0)$, one has that $\widetilde y_n$ satisfies a system like~\eqref{AdjvectPFSy} with regular data. Proposition~\ref{adjexist} provides the regularity and the previous formula. In fact, the previous equality and~\eqref{transp} also provide
	\begin{equation}\label{a61}
	\left\{
	\begin{array}{l}
\|y_n\|_{L^2(L^2)} \leq e^{CT}\left( \|y_0\|_{H^{-1}}+ \|v\|_{L^2(0,T)} +\|f\|_{L^2(L^2)}\right), \\
	\noalign{\smallskip}
y_n \to y \hbox{ in } L^2(Q_T; \R^2) \quad \hbox{and} \quad y_{n,xx} \to y_{xx} \hbox{ in } \mathcal{D}'(Q_T; \R^2),
	\end{array}
	\right.
	\end{equation}
for a new constant $C=C(D,A)>0$.

On the other hand, one has
	$$
\intdoble_{Q_T} y_{n,xx} \cdot \psi \, dx \, dt = \intdoble_{Q_T} y_n \cdot \psi_{xx} \, dx \, dt - \int_0^T B^* \psi_x(0,t) v^n (t) \, dt ,
	$$
for every $\psi \in L^2(0,T;H^2 (0,\pi;\mathbb{R}^2)\cap H_0^1(0,\pi;\mathbb{R}^2)) $. From this equality we deduce that the sequence $\{y_{n,xx} \}_{n \ge 1}$ is bounded in $L^2 (0,T; ( H^2(0,\pi;\mathbb{R}^2) \cap H_0^1(0,\pi;\mathbb{R}^2) )' )$. This property together with~\eqref{a61} gives $y_{xx} \in L^2 (0,T; ( H^2(0,\pi;\mathbb{R}^2) \cap H_0^1(0,\pi;\mathbb{R}^2) )' )$ and~\eqref{a60}.

Combining the identity $y_t = Dy_{xx} - Ay +f$ and the regularity property for $y_{xx}$, we also see that $y_t \in L^2 (0,T; ( H^2(0,\pi;\mathbb{R}^2) \cap H_0^1(0,\pi;\mathbb{R}^2) )' )$ and
	\begin{equation*}
\| y_{t} \|_{L^2 ((H^2 \cap H_0^1)')} \leq e^{CT}\left( \|y_0\|_{H^{-1}} + \|v\|_{L^2(0,T)} + \|f\|_{L^2(L^2)}\right) ,
	\end{equation*}
for a constant $C = C(D,A) > 0$. Therefore, $y \in C^0 ([0,T];X)$, where $X$ is the interpolation space 
	$$
X = \left[ L^2(0,\pi; \R^2), ( H^2(0,\pi;\mathbb{R}^2) \cap H_0^1(0,\pi;\mathbb{R}^2) )'  \right]_{1/2} \equiv H^{-1}(0,\pi; \R^2).
	$$
In conclusion, we have proved~\eqref{acot2}. Finally, it is not difficult to check that $y(\cdot , 0) = y_0$ in $H^{-1} (0, \pi; \R^2)$. This ends the proof.
\end{proof}
%
%
\end{section}
%
%

%
%
%
\begin{section}{}\label{APPX2_}

\setcounter{equation}{0}

In this appendix we will provide a positive answer on the null controllability of the phase-field system \eqref{PFSy} in the case  $c=0$. The computations and ideas used for obtaining this controllability result follow the ideas developed for the cases $c=1$ and $c=-1$.

Let us recall that $\tilde{\theta}=\tilde{\theta}(x,t)$ denotes the temperature of the material and the phase-field function $\tilde{\phi}=\tilde{\phi}(x,t)$ describes the phase transition of the material (solid or liquid) in such a way that $\tilde{\phi}=1$ means that the material is in solid state, $\tilde{\phi}=-1$ in liquid state and $\tilde \phi=0$ is an intermediate (mushy) phase.

In Theorem~\ref{NullC}, we proposed a local exact controllability result for the phase-field system~\eqref{PFSy} to the trajectories $(0,-1)$ or $(0,1)$. Our objective here is to prove a local null controllability result for the same system. 

Let us consider the phase-field system~\eqref{PFSy} with $c=0$, that is to say, the system 
	\begin{equation}\label{PFSyApx_}
	\left\{
	\begin{array}{ll}
\displaystyle \tilde{\theta}_t - \xi\tilde{\theta}_{xx} + \dfrac{1}{2}\rho\xi\tilde{\phi}_{xx} + \dfrac{\rho}{\tau}\tilde{\theta} = f_1(\tilde{\phi})	& \mbox{in } Q_T  , 
	\\
	\noalign{\smallskip}
\displaystyle \tilde{\phi}_t - \xi\tilde{\phi}_{xx} - \dfrac{2}{\tau}\tilde{\theta} = f_2(\tilde{\phi}) & \mbox{in }  Q_T, 
	\\
	\noalign{\smallskip}
\displaystyle \tilde{\theta}(0,\cdot) = v,\ \tilde{\phi}(0,\cdot) = 0,\ \tilde{\theta}(\pi,\cdot)=0 , \ \tilde{\phi}(\pi,\cdot) = 0 & \mbox{on }  (0,T),	 
	\\
	\noalign{\smallskip}
\displaystyle \tilde{\theta}(\cdot,0) = \tilde{\theta}_0, \  \tilde{\phi}(\cdot,0)=\tilde{\phi}_0 & \mbox{in }  (0,\pi).		\end{array}
	\right. 
	\end{equation}
where $\xi$, $\rho$ and $\tau$ are positive parameters and the nonlinear terms $f_1(\tilde{\phi})$ and $f_2(\tilde{\phi})$ are given by
	\[ 
\displaystyle f_1(\tilde{\phi}) = -\frac{\rho}{4\tau}\left( \tilde{\phi}-\tilde{\phi}^3 \right) \quad \mbox{and} \quad f_2(\tilde{\phi}) = \frac{1}{2\tau} \left( \tilde{\phi}-\tilde{\phi}^3 \right). 
	\]
For this system, a linearization around the equilibrium $(0,0)$ provides the following linear problem in vectorial form:
	\begin{equation}\label{vectPFSyApx_}
	\left\{
	\begin{array}{ll}
y_t - D y_{xx} + \hat{A} y = 0 & \mbox{in } Q_T,\\
	\noalign{\smallskip}
y(0,\cdot) = Bv, \quad y(\pi,\cdot)=0 & \mbox{on } (0,T), \\
	\noalign{\smallskip}
y (\cdot,0) = y_0, & \mbox{in }  (0,\pi), \\
	\end{array}
		\right. 
	\end{equation}
with $y_0 = ( \theta_0, \phi_0)$ ,  $y = ( \theta, \phi)$ and
	\begin{equation}\label{ADBApx_}
D = \left(
	\begin{array}{cc}
\xi & -\dfrac{1}{2} \rho \xi \\
	\noalign{\smallskip}
0 &  \xi
	\end{array}
	\right),
\quad
\hat{A} =
	\left(
	\begin{array}{cc}
\dfrac{\rho}{\tau} & \dfrac{\rho}{4\tau} \\
	\noalign{\smallskip}
-\dfrac{2}{\tau} & -\dfrac{1}{2\tau}
	\end{array}
\right),
\quad
B =
\left(
	\begin{array}{cc}
1 \\ 0
	\end{array}
	\right).
	\end{equation}

Following the same ideas used in the Appendix~\ref{APPX1_}, we can prove that, for every $y_0 \in H^{-1}(0,\pi; \R^2)$ and $v \in L^2(0,T)$, system~\eqref{vectPFSyApx_} has a unique solution by transposition (see Definition~\ref{deft_}) $y \in L^2(Q_T;\mathbb{R}^2)\cap C^0([0,T];H^{-1}(0,\pi;\mathbb{R}^2))$ which depends continuously on the data:
	\begin{equation*}
\| y \|_{L^2(L^2)} + \| y \|_{C^0(H^{-1})} \leq Ce^{CT} \left( \| y_0 \|_{H^{-1}} + \|v\|_{L^2 (0,T)} \right) ,
	\end{equation*}
for a constant $C>0$ only depending on the parameters $\xi$, $\rho$ and $\tau$ in system~\eqref{PFSyApx_}.

In order to state the null controllability result for systems~\eqref{PFSyApx_} and~\eqref{vectPFSyApx_}, let us consider the vectorial operators
	\begin{equation}\label{LL*Apx_}
\hat{L} = -D \partial_{xx} + \hat{A} \quad \mbox{and} \quad \hat{L}^* = -D^* \partial_{xx} + \hat{A}^*,
	\end{equation}
with domains $D(\hat{L})=D(\hat{L}^*)=H^2(0,\pi;\mathbb{R}^2) \cap H_0^1(0,\pi;\mathbb{R}^2)$. 

The first result in this appendix establishes the approximate controllability of system~\eqref{vectPFSyApx_} at time $T>0$. One has:
\begin{theorem} \label{CAproximada_Apx_}
	Let us consider $\xi$, $\rho$ and $\tau$ three positive real numbers and let us fix $T>0$. Then, system~\eqref{vectPFSyApx_} is approximately controllable in $H^{-1} (0, \pi; \R^2)$ at time $T $ if and only if the eigenvalues of the operators $\hat L$ and $\hat L^* $ are simple. Moreover, this equivalence amounts to the condition
	\begin{equation}\label{H2Apx_}
4\xi^2\tau^2(\ell^2 - k^2)^2 - 8\xi\rho\tau(\ell^2 + k^2)- 4 \rho-1 \neq 0, \quad \forall k,\ell\geq1, \quad \ell >k.
	\end{equation}
\end{theorem}

The second result in this appendix establishes the null controllability result at time $T>0$ of system~\eqref{vectPFSyApx_} and reads as follows:

\begin{theorem} \label{CNull_Apx_}
	Let us us fix $T>0$ and consider $\xi$, $\rho$ and $\tau$ positive real numbers satisfying~\eqref{H2Apx_} and 
	\begin{equation}\label{H1Apx_}
\xi \neq \dfrac{1}{j^2}\dfrac{\rho}{\tau},\quad \forall j\geq1.
	\end{equation}
Then, system~\eqref{vectPFSyApx_} is exactly controllable to zero in $H^{-1} (0, \pi; \R^2)$ at time $T>0$. Moreover, there exist  two positive constants $C_0$ and $M$, only depending on $\xi$, $\rho$ and $\tau$, such that for any $T  > 0$, there is a bounded linear operator $ \mathcal{C}_T^{(0)}: H^{-1}(0,\pi;\mathbb{R}^2) \rightarrow L^2(0,T) $ satisfying
	\begin{equation*}
\| \mathcal{C}_T^{(0)} \|_{\mathcal{L}(H^{-1}(0,\pi;\mathbb{R}^2),L^2(0,T))} \leq C_0 \, e^{M/T },
	\end{equation*}
and such that the solution 
	$$
y =( \theta , \phi) \in L^2 (Q_T; \R^2) \cap C^0([0,T];H^{-1}(0,\pi;\mathbb{R}^2)) 
	$$ 
of system~\eqref{vectPFSyApx_} associated to $y_0 = ( \theta_0 , \phi_0) \in H^{-1}(0,\pi;\mathbb{R}^2)$ and $v=\mathcal{C}_T^{(0)}(y_0)$ satisfies $y(\cdot , T)=0$. 
\end{theorem}

\begin{remark}
Observe that assumptions~\eqref{H2Apx_} and~\eqref{H1Apx_} play the role in Theorems~\ref{CAproximada_Apx_} and~\ref{CNull_Apx_} of conditions~\eqref{H2} and~\eqref{H1} in Theorems~\ref{CAproximada_} and~\ref{CNull_}.
\end{remark}

The local null controllability result for the nonlinear system~\eqref{PFSyApx_} is given in the next result:  

\begin{theorem}\label{NullCApx_}
Let us consider $\xi,\tau$ and $\rho$ three  positive numbers satisfying~\eqref{H2Apx_} and~\eqref{H1Apx_}, 
and let us fix $T>0$. Then, there exist $\varepsilon>0$ such that, for any $(\tilde{\theta}_0,\tilde{\phi}_0)\in H^{-1}(0,\pi)\times H_0^1(0,\pi)$ fulfilling
\begin{equation*}
	\|\tilde{\theta}_0\|_{H^{-1}} + \|\tilde{\phi}_0\|_{H_0^1} \le \varepsilon,
\end{equation*}
there exists $v\in L^2(0,T)$ for which system \eqref{PFSyApx_} has a unique solution 
	$$
(\tilde\theta,\tilde{\phi})\in \left[L^2(Q_T)\cap C^0([0,T]; H^{-1}(0,\pi; \R^2)) \right] \times C^0(\overline{Q}_T)
	$$
which satisfies
	\begin{equation*}
\tilde \theta (\cdot, T) = 0 \quad \hbox{and} \quad \tilde \phi (\cdot, T) = 0 \quad \hbox{in } (0, \pi).
	\end{equation*}
\end{theorem}

The proofs of Theorems~\ref{CAproximada_Apx_},~\ref{CNull_Apx_} and~\ref{NullCApx_} follow the same reasoning and ideas of the proofs of Theorems~\ref{CAproximada_},~\ref{CNull_} and~\ref{NullC}. They are based on an exhaustive study of the eigenvalues and eigenfunctions of the operators $\hat L$ and $\hat L^*$. In this sense, the properties of these eigenvalues and eigenfunctions are very close to the properties of the spectra of the operators $L$ and $L^*$ (see~\eqref{LL*}. Indeed, we have the following result.

\begin{proposition}\label{LspecApx_}
Let us consider the operators $\hat{L}$ and $\hat{L}^*$ given in~\eqref{LL*Apx_} (the matrices $D$ and $\hat{A}$ are given in~\eqref{ADBApx_}). Then, 
\begin{enumerate}
	\item The spectra of $\hat{L}$ and $\hat{L}^*$ are given by $\sigma(\hat{L}) = \sigma (\hat{L}^*)= \{\hat{\lambda}_k^{(1)}, \hat{\lambda}_k^{(2)} \}_{k\geq1}$ with
	\begin{equation}\label{lamk12Apx_}
\hat{\lambda}_k^{(1)} = \xi k^2 + \dfrac{2\rho + 1}{4\tau} - \hat{r}_k ,\quad \hat{\lambda}_k^{(2)} = \xi k^2 + \dfrac{2\rho + 1}{4\tau} + \hat{r}_k, \quad \forall k \ge 1,			
	\end{equation}
where
	\begin{equation*}
\hat{r}_k:= \sqrt{ \dfrac{\xi\rho}{\tau} k^2 + \left(\dfrac{2\rho+1}{4\tau}\right)^2 }.
	\end{equation*}
\item For each $k\geq1$, the eigenspaces of $\hat{L}$ (resp., $\hat{L}^*$) corresponding to $\hat{\lambda}_k^{(1)}$ and $\hat{\lambda}_k^{(2)}$ are respectively generated by
	\begin{equation*}
\hat{\Psi}_k^{(1)} = \frac{1}{8\sqrt{\tau \hat{r}_k}}\left(
	\begin{array}{c}
1-2\rho+4\tau \hat{r}_k \\
	\noalign{\smallskip}
8
	\end{array}
	\right)
	\eta_k,
\quad
\hat{\Psi}_k^{(2)} = \frac{1}{8\sqrt{\tau\hat{r}_k}}
	\left(
	\begin{array}{c}
1-2\rho-4\tau \hat{r}_k\\
	\noalign{\smallskip}
8
\end{array}
	\right)
	\eta_k, 
\end{equation*}
(resp., 
	\begin{equation*}
\hat{\Phi}_k^{(1)} = \frac{1}{8\sqrt{\tau \hat{r}_k}}
	\left(
	\begin{array}{c}
8 \\
	\noalign{\smallskip}
2\rho-1+4\tau \hat{r}_k
	\end{array}
	\right)
\eta_k,
	\quad
\hat{\Phi}_k^{(2)} = \frac{-1}{4\sqrt{\tau \hat{r}_k}}
	\left(
	\begin{array}{c}
8 \\
	\noalign{\smallskip}	
2\rho-1-4\tau \hat{r}_k
	\end{array}
	\right)
	\eta_k  \hbox{)}.
\end{equation*}
(The function $\eta_k$ is given in~\eqref{sin}).
\item The sequences $\hat{\mathcal{F}} =\{\hat{\Psi}_k^{(1)}, \hat{\Psi}_k^{(2)} \}_{ k\geq1}$ and $\hat{\mathcal{F^*}} = \{\hat{\Phi}_k^{(1)}, \hat{\Phi}_k^{(2)} \}_{k\geq1}$ are such that
\begin{enumerate}[i)]
	\item $\hat{\mathcal{F}}$ and $\hat{\mathcal{F}}^*$ are biorthogonal sequences.
	\item	$\hbox{span}  \left(\hat{\mathcal{F}} \right)$ and $\hbox{span}  \left(\hat{\mathcal{F}}^* \right)$ are dense in $H^{-1}(0,\pi;\mathbb{R}^2)$, $L^2(0,\pi;\mathbb{R}^2)$ and $H_0^1(0,\pi;\mathbb{R}^2)$.
	\item $\hat{\mathcal{F}}$ and $\hat{\mathcal{F}}^*$ are unconditional bases for $H^{-1}(0,\pi;\mathbb{R}^2)$, $L^2(0,\pi;\mathbb{R}^2)$ and $H_0^1(0,\pi;\mathbb{R}^2)$.
\end{enumerate}	
\end{enumerate}
\end{proposition}

The proof of Proposition~\ref{LspecApx_} follows the same ideas of the proofs of Propositions~\ref{Lspec} and~\ref{SS3}. The details are left to the reader.

Observe that the expressions of the eigenvalues of $\hat L$ and $\hat L^*$ (see~\eqref{lamk12Apx_}) are close to those of operators $L$ and $L^*$ (see~\eqref{lamk12}). In fact, replacing $(\rho,\tau)$ by $(2\rho,2\tau)$ in~\eqref{lamk12}, we obtain~\eqref{lamk12Apx_}. So, we can repeat the computations of the proof of Proposition~\ref{propert} in order to proof the following results concerning the spectral analysis for  $\sigma(\hat{L}) = \sigma (\hat{L}^*)= \{\hat{\lambda}_k^{(1)}, \hat{\lambda}_k^{(2)} \}_{k\geq1}$: 

\begin{proposition}\label{propertApx_}
Under the assumptions of Proposition~\ref{LspecApx_}, the following properties hold:
	\begin{enumerate}
\item[(P1)] $ \{\hat{\lambda}_k^{(1)}\}_{k\geq1}$ and $\{\hat{\lambda}_k^{(2)}\}_{k\geq1}$ (see~\eqref{lamk12Apx_}) are increasing sequences satisfying
	\begin{equation*}
0 < \hat{\lambda}_{k}^{(1)} < \hat{\lambda}_k^{(2)},\quad \forall k\geq1.		
	\end{equation*}
\item[(P2)] The spectrum of $\hat{L}$ and $\hat{L}^*$ is simple, i.e., $\hat{\lambda}_{k}^{(2)} \neq \hat{\lambda}_{\ell}^{(1)}$, for all $k,\ell\geq 1$ if and only if the parameters $\xi$, $\rho$ and $\tau$ satisfy condition~\eqref{H2Apx_}.
\item[(P3)] Assume that the parameters $\xi$, $\rho$ and $\tau$ satisfy~\eqref{H1Apx_}, i.e., there exists $j\geq0$  such that 
	\begin{equation*}
\frac{1}{(j+1)^2}\frac{\rho}{\tau}<\xi<\frac{1}{j^2}\frac{\rho}{\tau}.
	\end{equation*}
Then, there exists an integer $k_0=k_0(\xi,\rho,\tau, j) \geq 1$ and a constant $C= C(\xi,\rho,\tau, j) >0$  such that
	\begin{equation*}
\hat{\lambda}_{k+j}^{(1)} < \hat{\lambda}_{k}^{(2)} < \hat{\lambda}_{k+1+j}^{(1)} < \hat{\lambda}_{k+1}^{(2)}<\cdots, \ \forall k\geq k_0, \ \hbox{and} \ \min_{k \ge k_0} \left\{ \hat{\lambda}_k^{(2)} - \hat{\lambda}_{k+j}^{(1)} , \hat{\lambda}_{k+j+1}^{(1)} - \hat{\lambda}_k^{(2)} \right\} > C. 
	\end{equation*}
\item[(P4)]  Assume now that the parameters $\xi$, $\rho$ and $\tau$ satisfy~\eqref{H2Apx_} and~\eqref{H1Apx_}. Then, one has: 
	\begin{equation*}
\displaystyle \inf_{k,\ell\geq1}|\hat{\lambda}_{k}^{(2)}-\hat{\lambda}_{\ell}^{(1)}|>0 ,
	\end{equation*}
and there exists a positive integer $k_1 \in \N$, depending on $\xi$, $\rho$ and $\tau$, such that
	\begin{equation*}
\min \left\{ \left| \hat{\lambda}_k^{(1)} - \hat{\lambda}_\ell^{(1)} \right| , \left| \hat{\lambda}_k^{(2)} - \hat{\lambda}_\ell^{(2)} \right|, \left| \hat{\lambda}_k^{(2)} - \hat{\lambda}_\ell^{(1)} \right|   \right\} \ge \frac{\xi}{2} |k^2 - \ell^2| , \quad \forall k, \ell \ge 1, \quad |k - \ell | \ge k_1 .
	\end{equation*}
\end{enumerate}
\end{proposition}
%
%

We also have:

\begin{proposition}\label{corh1h2Apx_}
	Let us assume that the parameters $\xi,\rho,\tau$ satisfy \eqref{H2Apx_}. Then, the sequence $\{\hat{\lambda}_k^{(1)}, \hat{\lambda}_k^{(2)} \}_{k\geq1}$, given by~\eqref{lamk12Apx_}, can be rearranged into an increasing sequence $\hat{\Lambda}= \{ \hat{\Lambda}_k\}_{k\geq1}$ that satisfies~\eqref{serconv} and $\hat{\Lambda}_k\neq\Lambda_n$, for all $k,n \in \N$ with $k\neq n$. In addition, if~\eqref{H1Apx_} holds, the sequence $\{\hat{\Lambda}_k\}_{k\geq1}$ also satisfies \eqref{56} and \eqref{56_2}.
\end{proposition}

As said before, the proofs of Theorems~\ref{CAproximada_Apx_} and~\ref{CNull_Apx_} follow the same ideas of the proofs of Theorems~\ref{CAproximada_} and~\ref{CNull_}. To be precise, Theorem~\ref{CAproximada_Apx_} can be deduce from item~$(P2)$ in Proposition~\ref{propertApx_}. On the other hand, Theorem~\ref{CNull_Apx_} can be proved combining  Proposition~\ref{propertApx_}, Proposition~\ref{corh1h2Apx_} and Lemma~\ref{1.5M}, and following the same reasoning of Theorem~\ref{CNull_}.

Finally, the same proof presented in Section~\ref{s5.2} can be easily adapted to Theorem~\ref{NullCApx_}, with the observation that the operator $\mathcal{N}$ in~\eqref{mthcalN_}, that represents the nonlinearity of the system \eqref{PFSy} with $c=\pm 1$, can be defined as follows:
	\begin{equation*}
\mathcal{N}(f) =
	\left(
	\begin{array}{l}
\frac{\rho}{4\tau}\tilde{\phi}^3 \\
	\noalign{\smallskip}
-\frac{1}{2\tau}\tilde{\phi}^3\\
	\end{array}
	\right).
	\end{equation*}
The proof of Theorem~\ref{NullCApx_} can be deduced applying the Banach Fixed-Point Theorem. The details are left to the reader.

\end{section}

\end{appendices}	
	

%
%
%

%

\begin{thebibliography}{99} 
	%
	\bibitem{AB} \textsc{F.~Alabau-Boussouira}, \emph{Controllability of cascade coupled systems of multi-dimensional evolution PDEs by a reduced number of controls}, C.~R.~Math.\ Acad.\ Sci.\ Paris~\textbf{350} (2012), no.~11-12, 577--582.
	%
	\bibitem{AB-L} \textsc{F.~Alabau-Boussouira, M.~L\'eautaud}, \emph{Indirect controllability of locally coupled wave-type systems and applications}, J.~Math.\ Pures Appl.~(9) \textbf{99} (2013), no.~5, 544--576.
	%
	\bibitem{AK-B-D-K} \textsc{F.~Ammar Khodja, A.~Benabdallah, C.~Dupaix, I.~Kostin}, \emph{Controllability to the trajectories of phase-field models by one control force}, SIAM J.~Control Optim.\ \textbf{42} (2003), no.~5, 1661--1680. 
	%
	\bibitem{ACML4} \textsc{F.~Ammar Khodja, A.~Benabdallah, M.~Gonz\'alez-Burgos, L.~de Teresa}, \emph{Recent results on the controllability of linear coupled parabolic problems: a survey}, Math.\ Control Relat.\ Fields~\textbf{1} (2011), no.~3, 267--306.
	%
	\bibitem {AK-B-GB-dT-2011} \textsc{F.~Ammar Khodja, A.~Benabdallah, M.~Gonz\'alez-Burgos, L.~de Teresa}, \emph{The Kalman condition for the boundary controllability of coupled parabolic systems. Bounds on biorthogonal families to complex matrix exponentials}, J.~Math.\ Pures Appl.~(9) \textbf{96} (2011), no.~6, 555--590.
	%
	\bibitem {MGBminT} \textsc{F.~Ammar Khodja, A.~Benabdallah, M.~Gonz\'alez-Burgos, L.~de Teresa}, \emph{Minimal time for the null controllability of parabolic systems: the effect of the condensation index of complex sequences}, J.~Funct.\ Anal.~\textbf{267} (2014), no.~7, 2077--2151.
	%
	\bibitem {FAMLnewminim}\textsc{F.~Ammar Khodja, A.~Benabdallah, M.~Gonz\'alez-Burgos, L.~de Teresa}, \emph{New phenomena for the null controllability of parabolic systems: Minimal time and geometrical dependence}, J.~Math.\ Anal.\ Appl.~\textbf{444} (2016), no.~2, 1071--1113.
	%
	\bibitem{BBGB} \textsc{A.~Benabdallah, F.~Boyer, M.~Gonz\'alez-Burgos, G.~Olive}, \emph{Sharp estimates of the one-dimensional boundary control cost for parabolic systems and application to the $N$-dimen\-sional boundary null controllability in cylindrical domains}, SIAM J.~Control Optim.~\textbf{52} (2014), no.~5, 2970--3001.
	%
	\bibitem{cag} \textsc{G.~Caginalp}, \emph{An analysis of a phase field model of a free boundary}, Arch.~Rational Mech.\ Anal.~\textbf{92} (1986), no.~3, 205--245.
	%
	\bibitem {C}  \textsc{J.-M.~Coron}, \emph{Control and Nonlinearity}, Mathematical Surveys and Monographs, 136, American Mathematical Society, Providence, RI (2007).
	%
	\bibitem{D-FC-GB-Z} \textsc{A.~Doubova, E.~Fern\'andez-Cara, M.~Gonz\'alez-Burgos, E.~Zuazua}, \emph{On the controllability of parabolic systems with a nonlinear term involving the state and the gradient}, SIAM J.~Control Optim.~\textbf{41} (2002), no.~3, 798--819. 
	%
	\bibitem{Fattorini} \textsc{H.O.~Fattorini}, \emph{Some remarks on complete controllability}, SIAM J.~Control~\textbf{4} (1966), 686--694.
	%
	\bibitem{FaRu} \textsc{H.O.~Fattorini, D.L.~Russel}, \emph{Exact controllability theorems for linear parabolic equations in one space dimension}, Arch.~Rational Mech.\ Anal.~\textbf{43} (1971), 272--292. 
	%
	\bibitem{FC-E}  \textsc{E.~Fern{\'a}ndez-Cara, E.~Zuazua}, \emph{The cost of approximate controllability for heat equations: the linear case}, Adv.~Differential Equations \textbf{5} (2000), no.~4-6, 465--514.
	%
	\bibitem{MGBBoundC} \textsc{E.~Fern\' andez-Cara, M. Gonz\' alez-Burgos, L. de Teresa}, \emph{Boundary controllability of parabolic coupled equations}, J.~Funct.\ Anal.~\textbf{259} (2010), no.~7, 1720--1758.
	%
	\bibitem{FC-GB-dT} \textsc{E.~Fern\' andez-Cara, M. Gonz\' alez-Burgos, L. de Teresa}, \emph{Controllability of linear and semilinear non-diagonalizable parabolic systems}, ESAIM Control Optim.\ Calc.\ Var.~\textbf{21} (2015), no.~4, 1178--1204. 
	%
	\bibitem{FC-Z}\textsc{E.~Fern\' andez-Cara, E.~Zuazua}, \emph{Null and approximate controllability for weakly blowing up semilinear heat equations}, Ann.~Inst.\ H.~Poincar\'e Anal.~Non Lin\'eaire \textbf{17} (2000), no.~5, 583--616. 
	%
	\bibitem{GB-PG} \textsc{M.~Gonz\'alez-Burgos, R.~P\'erez-Garc\'ia}, \emph{Controllability results for some nonlinear coupled parabolic systems by one control force}, Asymptot.~Anal. \textbf{46} (2006), no.~2, 123--162.
	%
	%
	\bibitem{tucsnack} \textsc{Y.~Liu, T.~Takahashi, M.~Tucsnak}, \emph{Single input controllability of a simplified fluid-structure interaction model}, ESAIM Control Optim.\ Calc.\ Var.~\textbf{19} (2013), no.~1, 20--42.
	%
	\bibitem{Mil04} \textsc{L.~Miller}, \emph{Geometric bounds on the growth rate of null-controllability cost for the heat equation in small time}, J.~Differential Equations~\textbf{204} (2004), no.~1, 202--226.
	%
	\bibitem{Seidman} \textsc{T.I. Seidman}, \emph{Two results on exact boundary control of parabolic equations}, Appl.~Math.\ Optim.~\textbf{11} (1984), no.~2, 145--152.
	%
	%
	\bibitem{TW} \textsc{M.~Tucsnak, G.~Weiss}, \emph{Observation and Control for Operator Semigroups}, Birkh\"auser Advanced Texts: Basler Lehrb\"ucher, Birkh\"auser Verlag, Basel (2009).
	%
	\bibitem{Z} \textsc{J. Zabczyk,} \emph{Mathematical Control Theory: An Introduction}, Systems \& Control: Foundations \& Applications, Birkh\"{a}user Boston, Inc., Boston, MA (1992).
	%
\end{thebibliography}
\end{document}